\documentclass{amsart}
\usepackage{amsmath,amsfonts,amssymb,amsthm,epsfig,amscd,epsfig,psfrag,comment}
\usepackage[all]{xy}

\newtheorem{Theorem}{Theorem}[section]
\newtheorem{Proposition}[Theorem]{Proposition} 
\newtheorem{Lemma}[Theorem]{Lemma}

\newtheorem{Corollary}[Theorem]{Corollary}

\newtheorem{Specialthm}{Theorem}

\newcommand{\puteps}[2][0.5]
{\includegraphics[scale=#1]{#2.eps}}

\theoremstyle{definition}

\newtheorem{Remark}[Theorem]{Remark}

\newcommand{\p}{{\mathbb{P}^1}}
\newcommand{\C}{\mathbb{C}}
\def\F{{\sf{F}}}
\def\sF{{\mathcal{F}}}
\def\G{{\sf{G}}}
\def\sG{{\mathcal{G}}}
\def\rH{{\widetilde{\mathrm{H}}}}
\def\T{{\sf{T}}}
\def\sT{{\mathcal{T}}}

\def\P{{\mathcal{P}}}
\def\Q{{\mathcal{Q}}}
\newcommand{\Kh}{\mathrm{H}_{\mathrm{alg}}}
\newcommand{\Khh}{\mathrm{H}_{\mathrm{Kh}}}
\newcommand{\Cone}{\mathrm{Cone}}
\newcommand{\Con}{\mathrm{Con}}
\def\O{{\mathcal O}}
\def\E{{\mathcal{E}}}
\def\H{{\mathcal{H}}}
\def\R{{\mathcal{R}}}
\def\dim{\mbox{dim}}
\newcommand{\othercomp}{V_n^i}
\newcommand{\uq}{U_q(\mathfrak{sl}_2)}
\newcommand{\base}{\mathbb{C}[q,q^{-1}]}
\newcommand{\oH}{\overline{\mathrm{H}}_{\mathrm{alg}}}
\newcommand{\uoker}{\Q}

\newcommand{\ver}{\shortparallel}
\newcommand{\hor}{=}
\newcommand{\sh}{M}
\newcommand{\triv}{\bigcirc}

\DeclareMathOperator{\spn}{span}
\DeclareMathOperator{\Hom}{Hom}
\DeclareMathOperator{\Ext}{Ext}

\begin{document}

\title[Knot homology via derived categories of coherent sheaves I, sl(2) case]{Knot homology via derived categories of coherent sheaves I, sl(2) case}

\author{Sabin Cautis}
\email{scautis@math.harvard.edu}
\address{Department of Mathematics\\ Rice University \\ Houston, TX}

\author{Joel Kamnitzer}
\email{jkamnitz@math.berkeley.edu}
\address{Department of Mathematics\\ UC Berkeley \\ Berkeley, CA}

\begin{abstract}
Using derived categories of equivariant coherent sheaves, we construct a categorification of the tangle calculus associated to sl(2) and its standard representation.  Our construction is related to that of Seidel-Smith by homological mirror symmetry.  We show that the resulting doubly graded knot homology agrees with Khovanov homology.
\end{abstract}

\date{\today}
\maketitle
\tableofcontents
\newpage

\section{Introduction}
\subsection{Categorification}
There is a diagrammatic calculus involving tangles and tensor products $V^{\otimes n} $ of the standard representation of the quantum group $ \uq $.  Suppose that one is given a planar projection of a tangle $ T $ with $ n $ free endpoints at the top and $ m $ free endpoints at the bottom.  Following Reshetikhin-Turaev \cite{RT}, one can associate to this tangle a map of $ \uq $ representations 
\begin{equation*}
\psi(T) : V^{\otimes n} \rightarrow V^{\otimes m}.
\end{equation*}
called the Reshtikhin-Turaev invariant of the tangle.  This is done by analyzing the tangle projection from top to bottom and considering each cap, cup, and crossing in turn.  To each cap we associate the standard map $ \base \rightarrow V \otimes V $, to each cup we associate the standard map $ V \otimes V \rightarrow \base $, to each right handed crossing we associate the braiding $ V \otimes V \rightarrow V \otimes V $, which is defined using the $R$-matrix of $ \uq $, and to each left handed crossing, the inverse of the braiding (up to a scaling factor).  The overall map $ \psi(T) $ is defined to be the composition of these maps.  Reshetikhin-Turaev \cite{RT} showed that $ \psi(T) $ does not depend on the planar projection of the tangle.  Thus we have a map
\begin{equation*}
\psi : \biggl\{ (n,m) \text{ tangles } \biggr\} \rightarrow \Hom_{\uq}(V^{\otimes n}, V^{\otimes m}).
\end{equation*}
If $ T = K $ is link, then $ \psi(K) $ is a map $ \base \rightarrow \base $ and $ \psi(K)(1) $ is the Jones polynomial.

Khovanov has proposed the idea of categorifiying this calculus.  We define a \textbf{graded triangulated category} to be a triangulated category with a automorphism, denoted $ A \mapsto A\{1\} $, which commutes with all of the triangulated structure.  Note that the Grothendieck group of a graded triangulated category is naturally a module over $ \base $ (where $ q $ acts by $ \{-1\} $).  

A \textbf{weak categorification} of the above calculus is a choice of graded triangulated category $ D_n $ for each $ n$ and a map
\begin{equation*}
\Psi : \biggl\{ (n,m) \text{ tangles } \biggr\} \rightarrow \biggl\{\text{ isomorphism classes of exact functors } D_n \rightarrow D_m \biggr\}
\end{equation*}
such that we recover the original calculus on the Grothendieck group level, ie we have $ K(D_n) \cong V^{\otimes n} $ as $ \base $ modules and $ [\Psi(T)] = \psi(T) $.  We will also insist that $ D_0 $ be the derived category of graded vector spaces.

One categorification was conjectured by Bernstein-Frenkel-Khovanov in \cite{BFK} and proven by Stroppel in \cite{S}.  In this categorification $ D_n $ was the derived category of a graded version of the direct sum of various parabolic category $ \mathcal{O} $ for the Lie algebra $ \mathfrak{gl}_n$.  Another categorification was constructed by Khovanov in \cite{Ktangle} and extended in \cite{CK}.  In this categorification $ D_n $ was the derived category of graded modules over a certain combinatorially defined graded algebra. 

One source of interest in categorification is link invariants.  Suppose that $ K $ is a link.  Then $ \Psi(K) $ is a functor from the derived category of graded vector spaces to itself.  We define $ H^{i,j}(K) $ to be the cohomology of the complex of graded vector spaces $ \Psi(K)(\mathbb{C}) $.  Thus, $ H^{i,j} $ is a bigraded knot homology theory whose graded Euler characteristic is the Jones polynomial.  In the Khovanov and Stroppel categorifications, this knot homology theory is the celebrated Khovanov homology \cite{Kknot}.

\subsection{Categorification via derived categories of coherent sheaves}
In this paper, we present a categorification where the categories $ D_n = D(Y_n) $ are the derived categories of $ \mathbb{C}^\times $-equivariant coherent sheaves on certain smooth projective varieties $ Y_n $ (see section \ref{se:varieties} for their definition). 

Given a planar projection of a $(n,m) $ tangle $ T $, we construct a functor $ \Psi(T) : D(Y_n) \rightarrow D(Y_m) $ by associating to each cup, cap, and crossing certain basic functors $$ \G_n^i : D(Y_{n-2}) \rightarrow D(Y_n),\quad \F_n^i : D(Y_n) \rightarrow D(Y_{n-2}),\quad  \T_n^i(l) : D(Y_n) \rightarrow D(Y_n). $$ 
The functors $ \G_n^i, \F_n^i $ are defined using a correspondence $ X_n^i $ between $ Y_{n-2} $ and $ Y_n $ (see section \ref{se:capcup}).  The functors $ \T_n^i(l) $ are defined using a correspondence $ Z_n^i $ between $ Y_n $ and itself.  They also admit other descriptions as twists in the spherical functors $ \G_n^i $ (see section \ref{se:crossing}).

In section \ref{se:invariance} (which is the bulk of this paper), we check certain relations among these basic functors corresponding to Reidemeister and isotopy moves between different projections of the same tangle.  We obtain the following result.
\begin{Specialthm}[Theorem \ref{thm:main}]
The isomorphism class of $\Psi(T) $ is an invariant of the tangle $ T $.
\end{Specialthm}
Among the relations that we check, the most basic is that $ \T_n^i $ is an equivalence.  Using a result of Horja \cite{Ho}, this follows from the fact that it is a spherical twist (see section \ref{se:twists}).  The other important relation is that $ \T_n^i $ and $ \T_n^{i+1} $ braid.  Here our proof uses ideas of Anno \cite{A} which generalize the method of Seidel-Thomas \cite{ST} (see section \ref{se:twistbraid}).

In section \ref{se:Kgroups}, we consider $ \Psi(T) $ on the level of Grothendieck group.  Via explicit computations we prove the following result.
\begin{Specialthm}[Theorem \ref{th:Kgroup}]
On the Grothendieck groups, $ \Psi(T) $ induces the representation theoretic map $ \psi(T) $.
\end{Specialthm}

Thus, we have constructed a weak categorification of the above tangle calculus.

In section \ref{se:cobordism}, using a method introduced by Khovanov \cite{Kcobordism}, we associate to each tangle cobordism $T \rightarrow T'$ a natural transformation $\Psi(T) \rightarrow \Psi(T')$ thus showing $\Psi$ is an invariant of tangle cobordisms.

A natural question to ask is how this categorification compares with those mentioned above due to Khovanov and Stroppel.  For general tangles this is a difficult question because one must compare very differently defined triangulated categories.  However, for links the situation is simpler.  In section \ref{se:unorient} we show that our knot homology satisfies the same long exact sequence as Khovanov homology (Corollary \ref{th:longexact}) and we show that:
\begin{Specialthm}[Theorem \ref{th:Khov}]
Let $ K $ be a link and let $ H_{\text{alg}}^{i,j}(K) $ denote the knot homology obtained from our categorification.  We have $ H_{\text{alg}}^{i,j}(K) \cong H_{\text{Kh}}^{i+j,j}(K) $. 
\end{Specialthm}

Finally, in section \ref{se:reduced} we describe a reduced version of $H_{\text{alg}}$ for marked links. 

\subsection{Motivation and future work}
In this work, we were free to make many choices.  Namely we chose the varieties $ Y_n $ and we chose how to define a functor for each cap, cup, and crossing.  These choices were motivated by two considerations.  

First, Seidel-Smith \cite{SS} constructed a knot homology theory defined using Floer homology of Lagrangians in a certain sequence of symplectic manifolds $ M_n$.  Our variety $ Y_{2n} $ is a compactification of $ M_n $, after a change of complex structure.  Our knot homology theory is related to theirs by homological mirror symmetry (or more precisely hyperK\"ahler rotation).  The details of this connection are presented in section \ref{se:SS}.

A second motivation is that the variety $ Y_n $ arises in the geometric Langlands program as the convolution product
\begin{equation*}
Gr^{\omega} \tilde{\times} \cdots \tilde{\times} Gr^{\omega}.
\end{equation*}
Here $ \omega $ is the minuscule coweight of $ PGL_2 $ and $ Gr^{\omega} $ is the corresponding $ PGL_2(\mathbb{O}) $ orbit in the affine Grassmannian of $ PGL_2 $.  Via the geometric Satake correspondence \cite{MV}, this convolution product corresponds to the representation $ V^{\otimes n} $ of $ SL_2$.  The correspondences $ X_n^i, Z_n^i $ also have natural interpretations in this context. More details of this relationship will be explained in \cite{CK2}.

This perspective suggests how to generalize the foregoing construction to Lie algebras and representations other that $ \mathfrak{sl}_2 $ and its standard representation.  In a future work \cite{CK2}, we will give this construction in the next simplest case, the categorification of the diagrammatic calculus associated to $ \mathfrak{sl}_m $ and its standard representation.

\subsection{Acknowledgements}
In the initial stages, we benefited greatly from conversations with Denis Auroux and Roman Bezrukavnikov.  We also thank them for organizing the mirror symmetry and representation theory seminar at MIT in 2005-2006 where we began this project.  We thank Rina Anno, David Ben-Zvi, Sasha Braverman, Chris Douglas, Dennis Gaitsgory, Joe Harris, Brendan Hassett, Reimundo Heluani, Mikhail Khovanov, Allen Knutson, Kobi Kremnizer, Andrew Lobb, Scott Morrison, Nick Proudfoot, Ivan Smith, and Edward Witten for helpful conversations.

The second author thanks the American Institute of Mathematics for support and the mathematics departments of MIT and UC Berkeley for hospitality. 

\section{Geometric background}
We begin by describing a certain geometric setup upon which our work is based.

\subsection{The varieties} \label{se:varieties} 
Let $ N $ be a fixed large integer.  Fix a vector space $ \C^{2N} $ of dimension $2N$ and a nilpotent linear operator $ z : \C^{2N} \rightarrow \C^{2N} $ of Jordan type $(N, N) $.  More explicitly, we choose a basis $ e_1, \dots, e_N, f_1, \dots, f_N $ for $ \C^{2N} $ and define $ z $ by $ z e_i = e_{i-1}, z f_i = f_{i-1}, zf_1 = 0 = z e_1 $.

We consider now the following variety, denoted $ Y_n $, defined by
\begin{equation*}
Y_n := \{ (L_1, \dots, L_{n}) : L_i \subset \C^{2N} \text{ has dimension $ i $}, L_1 \subset L_2 \subset \dots \subset L_{n}, \text{ and } z L_i \subset L_{i-1} \}.
\end{equation*}
Note that $Y_n $ is independent of $N$ as long as $ N \ge 2n $, which we will always assume.  Alternatively, we can set $ N = +\infty$, in which case we can think of $ \C^{2N} $ as $ \C^2 \otimes \C[z^{-1}] $.

Note that $ Y_n $ is a smooth projective variety of dimension $ n $.  In fact, $ Y_{k+1} $ is a $\p $ bundle over $ Y_k $.  To see this, suppose that we have $ (L_1, \dots, L_k) \in Y_k $ and are considering possible choices of $ L_{k+1} $.  It is easy to see that we must have $ L_k \subset L_{k+1} \subset z^{-1}(L_k) $.  Since $ z^{-1}(L_k) / L_k $ is always two dimensional, this fibre is a $ \p $.  Hence the map $ Y_{k+1} \rightarrow Y_k $ is a $ \p $ bundle.  

There is an action of $\C^\times$ on $\C^{2N}$ given by $t \cdot e_i = t^{-2i}e_i$ and $t \cdot f_i = t^{-2i}f_i$. The reason we choose this action rather than $t \cdot e_i = t^{-i} e_i$ is to avoid having to shift the grading by half-weights later on. Notice that 
$$t \cdot (ze_i) = t \cdot e_{i-1} = t^{-2i+2}e_{i-1} = t^{2} t^{-2i}e_{i-1} = t^{2}z(t^{-2i} e_i) = t^{2}z(t \cdot e_i)$$
and similarly $t \cdot (zf_i) = t^{2}z(t \cdot f_i)$. So for any $v \in \C^{2N}$ we have 
$$t \cdot (zv) = t^{2}z(t \cdot v).$$
Thus $t \cdot (zL_i) = z(t \cdot L_i)$ so if $zL_i \subset L_{i-1}$ then $t \cdot zL_i \subset t \cdot L_{i-1}$ which means $z(t \cdot L_i) \subset t \cdot L_{i-1}$. Consequently, the $\C^\times$ action on $\C^{2N}$ induces a $\C^\times$ action on $Y_n$ by
$$t \cdot (L_1, \dots, L_{n}) = (t \cdot L_1, \dots, t \cdot L_{n}).$$

\subsection{Some diagrams}
The various spaces $ Y_n $ are related by the following diagrams.  For each $ 1 \le i \le n $, define the following subvariety of $ Y_n $, 
\begin{equation*}
X_n^i := \{ (L_1, \dots, L_{n}) \in Y_n : L_{i+1} = z^{-1} (L_{i-1}) \}.
\end{equation*}
Note that $ X_n^i $ is a $\C^\times$-equivariant divisor in $ Y_n $ and thus inherits a $\C^\times$ action from $Y_n$.

We also have a $\C^\times$-equivariant map
\begin{align*}
X_n^i &\xrightarrow{q} Y_{n-2} \\
(L_1, \dots, L_{n}) &\mapsto (L_1, \dots, L_{i-1}, zL_{i+2}, \dots, zL_{n}).
\end{align*}
That $q$ is $\C^\times$-equivariant is a consequence of $t \cdot (zL_i) = z(t \cdot L_i)$. 

Note that if $(M_1, \dots, M_{n-2}) \in Y_{n-2} $, then
\begin{equation*}
q^{-1}(M_1, \dots, M_{n-2}) = \{ (M_1, \dots, M_{i-1}, L_i, z^{-1}(M_{i-1}), z^{-1}(M_i), \dots, z^{-1}(M_{n-2}) \},
\end{equation*}
where $ L_i $ can be any subspace of $ \C^{2N} $ which lies between $ M_{i-1} $ and $ z^{-1}(M_{i-1}) $.  In particular there is a $ \p $ worth of choice for $ L_i $, and so the map $ q $ is a $ \p $ bundle.

To summarize we have the diagram of spaces
\begin{equation*}
\begin{CD}
X_n^i @>i>> Y_n \\
@VqVV \\
Y_{n-2} 
\end{CD}
\end{equation*}
where $ i $ is the $\C^\times$-equivariant inclusion of a divisor and $ q $ is an $\C^\times$-equivariant $ \p $ bundle.

In particular, we may view $ X_n^i $ as a $\C^\times$-equivariant subvariety of $Y_{n-2} \times Y_n $ via these maps.  Explicitly, we have
\begin{equation*}
X_n^i = \{ (L_\cdot, L'_\cdot) : L_j = L'_j \text{ for } j \le i-1, L_j = zL_{j+2} \text{ for } j \ge i-1 \}.
\end{equation*}

\subsection{Vector Bundles on $Y_n$ and $X_n^i$}

For $1 \le k \le n$ we have a $\C^\times$-equivariant vector bundle of dimension $k$ on $Y_n$ whose fibre over the point $(L_1, \dots, L_{n}) \in Y_n$ is $L_k$. Abusing notation a little, we will denote this vector bundle by $L_k$. Since $L_{k-1} \subset L_k$ we can consider the quotient $\E_k = L_k/L_{k-1}$ which is a $\C^\times$-equivariant line bundle on $Y_n$. Similarly, we get equivariant vector bundles $L_k$ on $X_n^i$ as well as the corresponding quotient line bundles $\E_k$ for any $1 \le k \le n$. Since $L_k$ on $Y_n$ restricted to $X_n^i$ is isomorphic to $L_k$ on $X_n^i$ we can omit subscripts and superscripts telling us where each $L_k$ and $\E_k$ lives. 

Another notational convention will be useful to us.  Often in this paper, we will deal with product of these spaces, such as $ Y_a \times Y_b \times Y_c $.  In this case, we will use the notation $ \E_i $ for $ \pi_1^*(\E_i)$, $\E'_i $ for $ \pi_2^*(\E_i) $ and $\E''_i $ for $\pi_3^*(\E_i)$.  In this way, the fibre of the line bundle $ \E'_i$ at a point $ (L_\cdot, L'_\cdot, L''_\cdot) \in Y_a \times Y_b \times Y_c $ is the vector space $ L'_i/L'_{i-1} $.

The map $z:L_i \rightarrow L_{i-1}$ has weight $2$ since 
$$(t \cdot z)(v) = t \cdot (z(t^{-1} \cdot v)) = t^{2}z(t \cdot t^{-1} \cdot v) = t^{2}zv.$$
If $\sF$ is a $\C^\times$-equivariant sheaf on $Y_n$ then we denote by $\sF\{m\}$ the same sheaf but shifted with respect to the $\C^\times$ action so that if $f \in \sF(U)$ is a local section of $\sF$ then viewed as a section $f' \in \sF\{m\}(U)$ we have $t \cdot f' = t^{-m}(t \cdot f)$. Using this notation we obtain $\C^\times$-equivariant maps 
$$z: L_i \rightarrow L_{i-1}\{2\} \ \text{ and } \ z: \E_i \rightarrow \E_{i-1}\{2\}.$$
We will repeatedly use these maps in the future. 

\subsection{A different description of $Y_n$}
Earlier, we saw that $ Y_n $ is a iterated $ \p $ bundle.  In fact these bundles are topologically trivial.  The following is a nice way of seeing this and will be useful in what follows.

Fix a vector space $ \C^2 $ with basis $ e, f $ and choose the Hermitian inner product on $ \C^2 $ such that $ e,f $ is an orthonormal basis.  Let $ \p $ denote the manifold of lines in $ \C^2 $.  Let $ a : \p \rightarrow \p $ denote the map which takes a line to its orthogonal complement.  In the usual coordinates on $ \p$, this map is given by $ z \mapsto -1/\bar{z} $.

\begin{Theorem} \label{th:isoman}
There is a diffeomorphism $Y_n \rightarrow ({\p})^{n} $.  Moreover under this diffeomorphism $ X_n^i $ is taken to the submanifold $ A_n^i := \{(l_1, \dots, l_i, a(l_i), l_{i+2}, \dots l_{n}) \} \subset ({\p})^{n} $.
\end{Theorem}

First, let $ C : \C^{2N} \rightarrow \C^2 $ denote the linear map which takes every $ e_i $ to $ e $ and every $ f_i $ to $ f $.  Next introduce the Hermitian inner product of $ \C^{2N} $ with orthonormal basis $ e_1, \dots, e_N, f_1, \dots, f_N $.

\begin{Lemma}
Let $ W \subset \C^{2N}$ be a subspace (of codimension at least two) such that $ zW \subset W $.  Then $ C $ restricts to an unitary isomorphism $ z^{-1} W \cap W^{\bot} \rightarrow \C^2 $.
\end{Lemma}

\begin{proof}
First, note that $ \dim(z^{-1} W \cap W^\bot) = 2 $, since $ W \subset z^{-1} W $ and $z$ has rank two.

So it suffices to show that $ C $ restricts to a unitary map.  Let $ v, w \in z^{-1}W \cap W^\bot $.

Let us expand $ v = v_1 + \dots + v_N, w = w_1 + \dots + w_N $, where $ v_i, w_i \in \spn(e_i, f_i) $. This is an orthogonal decomposition.  Note also that $ C $ does restrict to a unitary isomorphism $ \spn(e_i, f_i) \rightarrow \C^2 $.

Hence
\begin{equation*} 
\langle C(v), C(w) \rangle = \sum_{i,j} \langle C(v_i), C(w_j) \rangle 
\end{equation*}
while
\begin{equation*}
\langle v, w \rangle = \sum_i \langle v_i, w_i \rangle = \sum_i \langle C(v_i), C(w_i) \rangle.
\end{equation*}

So it suffices to show that 
\begin{equation} \label{eq:need}
\sum_{i,j, i \ne j} \langle C(v_i), C(w_j) \rangle = 0.
\end{equation}

Now, since $ w \in z^{-1} W $ and $v \in W^\bot $, we see that $ \langle v, z^k w \rangle = 0 $ for all $ k \ge 1$.  Note that 
\begin{equation*}
\langle v, z^k w \rangle = \langle C(v_1), C(w_{k+1}) \rangle + \dots + \langle C(v_{N-k}), C(w_N) \rangle.
\end{equation*} 
So adding up these equations for all possible $ k $, we deduce that
\begin{equation*}
\sum_{i < j} \langle C(v_i), C(w_j) \rangle = 0.
\end{equation*}

Similarly, we deduce that 
\begin{equation*}
\sum_{i < j} \langle C(w_i), C(v_j) \rangle = 0.
\end{equation*}
Taking the complex conjugate of this equation and adding it to the previous equation gives (\ref{eq:need}).
\end{proof}

\begin{proof}[Proof of Theorem \ref{th:isoman}]
Given $ (L_1, \dots, L_{n}) \in Y_n $, let $ M_1, \dots, M_{n} $ be the sequence of lines in $ \C^{2N} $ such that 
\begin{equation*}
L_k = M_1 \oplus \dots \oplus M_k \quad L_{k-1} \bot M_k.
\end{equation*}
Hence $ M_k $ is a one dimensional subspace of the two dimensional vector space $ z^{-1} L_{k-1} \cap L_{k-1}^{\bot} $.  By the lemma $ C(M_k) $ is a one dimensional subspace of $ \C^2 $.

Thus we have a map 
\begin{align*}
Y_n &\rightarrow (\p)^{n} \\
(L_1, \dots, L_{n}) &\mapsto (C(M_1), \dots, C(M_{n}))
\end{align*}

By induction on $ n $, it is easy to see that this map is an isomorphism.

Now, let us consider $ (L_1, \dots, L_{n}) \in X_n^i $.  Since $ L_{i+1} = z^{-1} L_{i-1} $, we see that $ M_i, M_{i+1} $ are both subspaces of $ z^{-1} L_{i-1} \cap L_{i-1}^\bot$.  As they are perpendicular, they are sent by $ C $ to perpendicular subspaces of $ \C^2 $.  Hence $ a(C(M_i)) = C(M_{i+1}) $ as desired.
\end{proof}

\subsection{Comparison with resolution of slices} \label{se:SS}
The main purpose of this subsection is to make contact with the work of Seidel-Smith.  For this section, we will work only with the even spaces $ Y_{2n} $.

\subsubsection{The Springer fibre}
We define the subspace $ \C^{2n} \subset \C^{2N} $ as the span of $ e_1, \dots, e_n, f_1, \dots, f_n $ and now we consider the subvariety $F_n $ of $ Y_{2n} $,
\begin{equation*}
F_n := \{ (L_1, \dots, L_{2n}) \in Y_{2n} : L_{2n} = \C^{2n} \}.
\end{equation*}
Note that $ z $ restricted to $ \C^{2n} $ is a nilpotent operator of Jordan type $ (n,n) $ and $ F_n $ is its Springer fibre.  So from the general theory of Springer fibres, $ F_n $ is a reducible connected projective variety of dimension $ n $.

Now, we define a certain open neighbourhood of $ F_n $, inside $ Y_{2n} $.  Define a linear operator $ P : \C^{2N} \rightarrow \C^{2n} $ by $ Pe_i = e_i $ if $ i \le n $ and $ P e_i = 0 $ if $ i > n $, and similarly for $ f_i$.  So $ P $ is a projection onto the subspace $ \C^{2n} $.  Define
\begin{equation*}
U_n := \{ (L_1, \dots, L_{2n} \in Y_{2n} : P(L_{2n}) = \C^{2n} \}.
\end{equation*}

\subsubsection{Nilpotent slices}
Let $ e_n $ be a nilpotent operator in $ \C^{2n} $ of type $(n,n) $, so in a particular basis $ e_n $ has the form
\begin{equation*}
\begin{bmatrix}
0& I&  & &  &\\
&  0& I&  & &\\
&& \dots& && \\
& & & 0& I& \\
&&&& 0&
\end{bmatrix},
\end{equation*}
where the entries are $ 2 \times 2 $ blocks.

Let $ f_n $ be the matrix which completes $ e_n $ to a Jacobson-Morozov triple, so
\begin{equation*}
f_n = \begin{bmatrix}
0& &&&& \\
(n-1)I&  0& &&& \\
& 2(n-2)I& 0 && \\
& & \dots& && \\
&&& (n-1) I& 0 &\\
\end{bmatrix},
\end{equation*}

Now, consider the set $ S_n := e_n + ker(f_n\cdot) \subset \mathfrak{sl}(\C^{2n}) $ where $ ker(f_n\cdot) $ is the kernel of left multiplication by $ F_n $.  It is easy to see that $ S_n $ is the set of traceless matrices of the form
\begin{equation*}
\begin{bmatrix}
0& I&  & &  &\\
&  0& I&  & &\\
& & \dots& && \\
& & & 0& I& \\
*& *& \dots& *& *&
\end{bmatrix}.
\end{equation*}

A small modification of $ S_n $ was considered by Seidel-Smith in \cite{SS}.  They considered $ e_n + ker(\cdot f_n) $, so the unknown entries in the matrix were on the left side.  In any case, we have the following result.

\begin{Lemma}[{\cite[Lemma 17]{SS}}]
$ S_n $ is a slice to the adjoint orbit through $ e_n $.
\end{Lemma}

Now we consider $ S_n \cap \mathcal{N} := \{ X \in S_n : X \text{ is nilpotent }\} $.  This is a singular variety.  It has a natural Grothendieck-Springer resolution, $ \widetilde{S_n \cap \mathcal{N}} $ defined by
\begin{equation*}
\widetilde{S_n \cap \mathcal{N}} := \{ (X, (V_1, \dots, V_{2n})) : X \in S_n \cap \mathcal{N}, (V_1, \dots, V_{2n}) \text{ is a complete flag in } \C^{2n}, \text{ and } X V_i \subset V_{i-1} \}.
\end{equation*}

\subsubsection{Relation to previous definition}
The following observation is originally due to Lusztig (\cite{L}).
\begin{Proposition} \label{prop:Unslice}
There is an isomorphism $ U_n \rightarrow \widetilde{S_n \cap \mathcal{N}} $ which is the identity on $ F_n $.
\end{Proposition}

\begin{proof}
Note that if $ (L_1, \dots, L_{2n}) $ is an element of $ U_n $, the projection map $ P $ restricts to an isomorphism between $ L_{2n} $ and $ \C^{2n}$.  We let $ P^{-1} $ denote its inverse.  

Note that on $L_{2n} $ we have a linear operator $ z $ which is the restriction of $ z $ from $ \C^{2N} $.  Now, we can compose $ z $ by the isomorphism $ P $ to get a linear operator $ X := P z P^{-1} $ on $ \C^{2n} $.  So, $ X(v) = P( z w ) $ where $ w $ in the unique element of $ L_{2n} $ such that $ P(w) = v $.  In particular, $ w = v + v' $ where $ v' $ lies in the span of $ e_{n+1}, \dots, e_{N}, f_{n+1}, \dots, f_{N} $.  Hence we see that $ X(e_i) = e_{i-1} + v'' $ where $ v'' $ lies in the span of $ e_n $ and $f_n$.  In particular such an $ X $ lies in $ S_n$.

We define the map
\begin{align*}
U_n &\rightarrow \widetilde{S_n \cap \mathcal{N}} \\
(L_1, \dots, L_{2n}) &\mapsto (P z P^{-1}, (P(L_1), \dots, P(L_{2n}))).
\end{align*}

We omit the construction of the inverse map.
\end{proof}

\begin{Remark}
It is known that $ \widetilde{S_n \cap \mathcal{N}} $ admits a hyperK\"akler structure (for example, since it is isomorphic to a Nakajima quiver variety by the work of Maffei \cite{M}) and thus so does $ U_n $.  This fact will be used below, but only for motivational purposes. 
\end{Remark}

\subsubsection{Seidel-Smith construction}
We can now explain one motivation for our work.

In their paper \cite{SS}, Seidel-Smith give a construction of a link invariant using symplectic geometry.  In particular, they consider the map $ \chi : S_n \rightarrow \C^{2n-1} / {\Sigma}_{2n} $, where the map $ \chi $ takes the eigenvalues of a matrix (here $ \Sigma_{2n} $ denotes the symmetric group and $ \C^{2n-1} $ denotes the subset of $ \C^{2n} $ of numbers whose sum is 0).  Over $ \C^{2n-1} \smallsetminus \Delta / {\Sigma}_{2n} $, the map $ \chi $ is a fibre bundle.  They use this fibre bundle to construct an action of the braid group $ B_{2n} = \pi_1(\C^{2n-1} \smallsetminus \Delta / {\Sigma}_{2n}) $ on the set of Lagrangian submanifolds of a given fibre $ M_n = \chi^{-1} (\lambda_1, \dots, \lambda_{2n}) $. 

Using this action they construct a link invariant by associating to each link $ K $, $ HF(L, \beta(L)) $, where $ HF(,) $ denotes Floer cohomology, $ \beta \in B_{2n} $ is a braid whose closure is $ K $, $ L $ is certain chosen Lagrangian in $ M_n $ not depending on $ K $, and $ \beta(L) $ denotes the above action of $ \beta $ on $ L $.

The homological mirror symmetry principle suggests that there should exist some derived category of coherent sheaves equivalent to the Fukaya category of the affine K\"ahler manifold $ M_n $.   In particular, the braid group should act on this derived category and it should be possible to construct a link invariant in the analogous manner, namely as $ \Ext(L, \beta(L))$ for an appropriate chosen object $ L $ in a certain derived category of coherent sheaves.  

By the general theory of nilpotent slices (see for example Lemma 21 of \cite{SS}), $ M_n $ is diffeomorphic to $ \widetilde{S_n \cap \mathcal{N}} $, but has a different complex structure (recall that $ \widetilde{S_n \cap \mathcal{N}} $ is hyperK\"ahler).  In particular, they are related by a classic ``deformation vs. resolution'' picture for hyperK\"ahler singularities.  Hence string theory suggests that the Fukaya category of $ M_n $ should be related to the derived category of coherent sheaves on $ U_n \cong \widetilde{S_n \cap \mathcal{N}} $ or perhaps the subcategory of complexes of coherent sheaves whose cohomologies are supported on the Springer fibre $ F_n$.  

In our paper, we found it more convenient to work with the full derived category of the compactification $Y_{2n} $ of $ U_n$ and to work with tangles, rather than just links presented as the closures of braids.  However, suppose that a link $ K $ is presented as the closure of braid $ \beta \in B_{2n} $.  Tracing through our definitions from section \ref{se:funtangle}, we see that (ignoring the bigrading)
\begin{equation*}
\Kh(K) = H(\Psi(K)(\C)) \cong \Ext_{D(Y_{2n})}(L, \beta(L)) = \Ext_{D(U_n)}(L, \beta(L))
\end{equation*}
where $ L $ is the structure sheaf of a certain component of the Springer fibre $ F_n $ (tensored with a line bundle).  Thus, our link invariant has the same form as would be expected from homological mirror symmetry.

\section{Background to FM transforms and twists}
\subsection{Fourier-Mukai transforms}
In this paper, we will define functors using Fourier-Mukai transforms, also called integral transforms.  We begin with some background, following \cite[section 5.1]{H}.

Let $ X $ be a smooth projective variety with a $\C^\times$-action. We will be working  with the bounded derived category of $\C^\times$-equivariant coherent sheaves on $X$, which we denote $ D(X) $.  All pullbacks, pushforwards, Homs, and tensor products of sheaves will be derived functors.  Also, we assume all spaces come equipped with a $\C^\times$ action while all maps and sheaves are assumed to be $\C^\times$-equivariant with respect to this action. 

Let $ X, Y $ be two smooth projective varieties.  A \textbf{Fourier-Mukai kernel} is any object $ \P $ of the derived category of $\C^\times$-equivariant coherent sheaves on $ X \times Y $.  Given $\P \in D(X \times Y)$, we may define the associated \textbf{Fourier-Mukai transform}, which is the functor
\begin{equation*}
\begin{aligned}
\Phi_\P : D(X) &\rightarrow D(Y) \\
\sF &\mapsto {\pi_2}_* (\pi_1^* (\sF) \otimes \P) 
\end{aligned}
\end{equation*}

Fourier-Mukai transforms have right and left adjoints which are themselves Fourier-Mukai transforms.  In particular, the right adjoint of $ \Phi_\P $ is the FM transform with respect to $ \P_R := \P^\vee \otimes \pi_2^* \omega_X [\dim(X)] \in D(Y \times X)$.  Similarly, the left adjoint of $ \Phi_\P $ is the FM transform with respect to $ \P_L := \P^\vee \otimes \pi_1^* \omega_Y [\dim(Y)] $, also viewed as a sheaf on $ Y \times X $.

We can express composition of Fourier-Mukai transform in terms of their kernel.  If $ X, Y, Z $ are varieties and $\Phi_\P : D(X) \rightarrow D(Y),  \Phi_\Q : D(Y) \rightarrow D(Z) $ are Fourier-Mukai transforms, then $ \Phi_\Q \circ \Phi_\P $ is a FM transform with respect to the kernel
\begin{equation*}
\Q \ast \P := {\pi_{13}}_*(\pi^*_{12}(\P) \otimes \pi^*_{23}(\Q)).
\end{equation*}
The operation $ \ast $ is associative.  Moreover by \cite{H} remark 5.11, we have $ (\Q * \P)_R \cong \P_R * \Q_R $.

If $\Q \in D(Y \times Z)$ then we define the functor $\Q \ast: D(X \times Y) \rightarrow D(X \times Z)$ by $\P \mapsto \Q \ast \P$. 

\begin{Lemma}\label{kerneladjoint} If $\Q \in D(Y \times Z)$ then the left adjoint of $\Q \ast$ is $(\Q \ast)^L = \Q_L \ast$ and the right adjoint $(\Q \ast)^R = \Q_R \ast$. 
\end{Lemma}
\begin{proof}
Let $\P \in D(X \times Y)$ and $\P' \in D(X \times Z)$. We need to show that there is a natural isomorphism
$$\Hom(\Q_L \ast \P', \P) \cong \Hom(\P', \Q \ast \P).$$
Using the definition of $ * $ and Grothendieck duality we have
\begin{eqnarray*}
\Hom(\Q_L \ast \P', \P) &=& \Hom( \pi_{12\ast}(\pi_{13}^\ast \P' \otimes \pi_{23}^\ast \Q^\vee \otimes \pi_3^\ast \omega_Z) [\dim(Z)], \P) \\
&\cong& \Hom( \pi_{13}^\ast \P' \otimes \pi_{23}^\ast \Q^\vee \otimes \pi_3^\ast \omega_Z [\dim(Z)]), \pi_{12}^\ast \P \otimes \pi_3^\ast \omega_Z [\dim(Z)]) \\
&=& \Hom( \pi_{13}^\ast \P', \pi_{23}^\ast \Q \otimes \pi_{12}^\ast \P) \\
&\cong& \Hom( \P', \pi_{13\ast}(\pi_{23}^\ast \Q \otimes \pi_{12}^\ast \P)) \cong \Hom( \P', \Q \ast \P)
\end{eqnarray*}
Since each isomorphism above is natural so is the composition. Notice $\Q_R^\vee = \Q \otimes \pi_1^\ast \omega_X^\vee [-\dim(X)]$ so $(\Q_R)_L = (\Q_R)^\vee \otimes \pi_1^\ast \omega_X [\dim(X)] = \Q.$ Replacing $\Q$ by $\Q_R$ in the above calculation shows 
\begin{eqnarray*}
\Hom(\Q \ast \P', \P) = \Hom((\Q_R)_L \ast \P', \P)
= \Hom( \P', \Q_R \ast \P)
\end{eqnarray*}
which means that $\Q_R \ast$ is the right adjoint of $\Q \ast$. 
\end{proof}

\begin{Remark} One can also define the functor $\ast \Q$ by $\P \mapsto \P \ast \Q$. Then the right adjoint becomes  $(\ast \Q)_R = \ast \Q_L$ since 
$$\Hom(\P' \ast \Q, \P) \cong \Hom(\P_R, \Q_R \ast \P'_R) \cong \Hom(\Q \ast \P_R, \P'_R) \cong \Hom(\P', \P \ast \Q_L).$$
Similarly one can show that the left adjoint is $(\ast \Q)_L = \ast \Q_R$. 
\end{Remark}

\subsection{Twists in spherical functors} \label{se:twists}
Let $ X, Y $ be smooth projective varieties and let $ \Phi_\P : D(X) \rightarrow D(Y) $ be a Fourier-Mukai transform with respect to $ \P $.  

There is a natural map $ \beta_\P : \P \ast \P_R \rightarrow \O_\Delta $ in $ D(Y \times Y) $.  By Lemma \ref{kerneladjoint}, $ (\P *)^R \cong \P_R *  $.  Hence we have a natural transformation, $ \gamma_\P : (\P *) \circ (\P_R * ) \rightarrow id $.  By definition $$ \beta_\P := \gamma_\P(\O_\Delta) : \P * \P_R \rightarrow \O_\Delta. $$  Moreover, $ \gamma_\P $ is actually determined completely by $ \beta_\P $ in the sense that $\gamma_\P (\sF) = \beta_P * \sF $ for any $ \sF \in D(A \times X) $.

Similarly, there is a natural transformation $ \tau_\P : id \rightarrow (P_R * ) \circ (\P * ) $ which leads to a morphism $ \sigma_\P : \O_\Delta \rightarrow \P_R * \P $.

The maps $ \beta_\P $ are compatible with convolution in the following sense.

\begin{Lemma} \label{th:convbeta}
Suppose that $\P \in D(X \times Y) $ and $ \Q \in D(Y \times Z) $ are Fourier-Mukai kernels. 
Then the following diagram commutes:
\begin{equation*}
\begin{CD}
\Q * \P * \P_R * \Q_R @>\Q * \beta_\P * \Q_R>> \Q * \Q_R \\
@VVV @V{\beta_\Q}VV \\
(\Q * \P) * (\Q * \P)_R @>\beta_{\Q*\P}>> \O_\Delta
\end{CD}
\end{equation*}
where the left hand arrow is the natural isomorphism.
\end{Lemma}

\begin{proof}
Since $ (\Q * \P) * \cong (\Q * ) \circ (\P * ) $, we have $ ((\Q * \P) * )^R \cong ( \P * )^R  \circ (\Q * )^R $.  This isomorphism is compatible with the isomorphism between $ (\Q * \P)_R $ and $ \P_R * \Q_R $.  Hence we conclude that there is a commutative diagram of natural transformations of functors:
\begin{equation*}
\begin{CD}
(\Q * ) \circ (\P * ) \circ (\P_R * ) \circ (\Q_R *) @>{ \Q* \gamma_\P(\Q_R * )}>> (\Q * ) \circ (\Q_R * ) \\
@VVV @V{\gamma_\Q}VV \\
((\Q * \P) *) \circ ((\Q * \P)_R * ) @>{\gamma_{\Q * \P}}>> id
\end{CD}
\end{equation*}

Now, apply these natural transformations to the object $ \O_\Delta $.  The top arrow becomes $ \Q *(\gamma_\P(\Q_R)) = \Q * \beta_\P * \Q_R $ by the above observation.   Hence the desired result follows.
\end{proof}

Let $ \sT_\P $ denote the Fourier-Mukai kernel in $ D(Y\times Y) $ defined (up to isomorphism) as the cone
\begin{equation*}
\sT_\P := \Cone( \P \ast \P_R \xrightarrow{\beta_\P} \O_\Delta).
\end{equation*}
So there is a distinguished triangle
\begin{equation} \label{eq:distsT}
\P \ast \P_R \xrightarrow{\beta_\P} \O_\Delta \rightarrow \sT_\P.
\end{equation}

Also we see that for any object $ \sF \in D(Y) $ we have a distinguished triangle
\begin{equation*}
\Phi_\P \Phi_\P^R(\sF) \rightarrow \sF \rightarrow \Phi_{\sT_\P} (\sF),
\end{equation*}
where the left map is the adjunction morphism.

Taking $ ()_L $ in (\ref{eq:distsT}), we see that there is a distinguished triangle
\begin{equation} \label{eq:distsTL}
(\sT_\P)_L \rightarrow \O_\Delta \rightarrow \P * \P_L \rightarrow (\sT_\P)_L[1]
\end{equation}

There are a number of results giving criteria for when the Fourier-Mukai transform with respect to the kernel $ \sT_\P $ is an equivalence.  The first is due to Seidel-Thomas \cite{ST} in the case that $ X $ is just a point and so $ \P $ is just an object of $ D(Y) $.  This was later generalized by Horja \cite{Ho} in a geometric context and Rouquier \cite{R} in a more abstract categorical context.  We also learnt about these generalizations from conversations with Anno and Bezrukavnikov.  

\begin{Theorem} \label{th:twistequiv}
Let $X, Y$ be two smooth projective varieties equipped with $\C^\times$ actions and $\P \in D(X \times Y)$. Assume that $\P_R \cong \P_L[k]\{l\}$ where $k = \dim(X) - \dim(Y)$ and $l$ is any integer.  Then there exists a sequence of adjoint maps
$$\O_\Delta \rightarrow \P_R \ast \P \rightarrow \O_\Delta[k]\{l\}$$
where the right map is the shift of the adjoint map $\P_L \ast \P \rightarrow \O_\Delta$. Suppose that this sequence forms a distinguished triangle and moreover suppose that 
\begin{eqnarray*}
\Hom(\P,\P[i]\{j\}) \cong \left\{
\begin{array}{ll}
\C & \mbox{if } i=0 \text{ and } j=0\\
0  & \mbox{if } i=k,k+1 \text{ and } j=l 
\end{array}
\right.
\end{eqnarray*}
Then $\Phi_{\sT_\P}$ is an equivalence.
\end{Theorem}

If the hypotheses of this theorem are satisfied, then we call $ \Phi_{\P} $ a \textbf{spherical functor} and we call $ \Phi_{\sT_\P} $ the \textbf{twist in the spherical functor}.  We also write $ \T_{\P} $ for $ \Phi_{\sT_\P} $.  This generalizes the notion of spherical object due to Seidel-Thomas \cite{ST} and also the notion of EZ spherical object due to Horja \cite{Ho}.

If $\P$ is a spherical functor then the inverse of $\T_{\P}$ is $(\T_{\P})^R \cong (\T_{\P})^L \cong \Phi_{(\sT_\P)_L} $.  

The proof that we give below is taken from Horja \cite{Ho}, although the hypotheses of the theorem are more general than his hypotheses.  

\begin{proof}
Applying $ * \P $ to (\ref{eq:distsT}), we obtain the distinguished triangle
\begin{equation*}
\sT_\P * \P[-1] \rightarrow \P * \P_R * \P \xrightarrow{\beta_\P * \P} \P
\end{equation*}
By hypothesis we have the distinguished triangle
\begin{equation*}
\O_\Delta \xrightarrow{\sigma_\P} \P_R * \P \rightarrow \O_\Delta[k]\{l\}
\end{equation*}
Applying $ \P * $ to this distinguished triangle gives
\begin{equation*}
\P \xrightarrow{\P* \sigma_\P} \P * \P_R * \P \rightarrow \P[k]\{l\}
\end{equation*}

By general theory of adjoints, we know that the composition of natural transformations 
\begin{equation*}
\P *  \xrightarrow{\P * \tau_\P} (\P * ) \circ (\P_R * ) \circ (\P * ) \xrightarrow{\gamma_\P(\P * )} \P * 
\end{equation*}
is the identity.  Applying this to $ \O_\Delta $, we see that the composition
\begin{equation*}
\P \xrightarrow{\P * \sigma_\P} \P * \P_R * \P \xrightarrow{\beta_\P * \P} \P
\end{equation*}
is the identity. Thus by Lemma \ref{th:twodist} we conclude that $\sT_\P * \P \cong \P [k+1]\{l\}.$  A similar argument shows that $ \P_R * \sT_P \cong \P_R[k+1]\{l\}$.

Applying $\sT_\P \ast$ to (\ref{eq:distsTL}), we obtain the distinguished triangle
$$\sT_\P \rightarrow \sT_\P \ast \P \ast \P_L \rightarrow \sT_\P \ast (\sT_P)_L[1]$$
which, using $\sT_\P \ast \P \cong \P[k+1]\{l\}$, simplifies to 
$$\sT_\P \rightarrow \P \ast \P_L[k+1]\{l\} \rightarrow \sT_\P \ast (\sT_\P)_L[1].$$
Using that $\P_R \cong \P_L[k]\{l\}$ and shifting by $-1$ we obtain the distinguished triangle
\begin{equation} \label{eq:firstdist}
\sT_\P[-1] \rightarrow \P \ast \P_R \rightarrow \sT_\P \ast (\sT_\P)_L.
\end{equation}
Notice that the map $\sT_\P[-1] \rightarrow \P \ast \P_R$ comes from the map $\sT_\P \rightarrow \sT_\P \ast \P \ast \P_L$ which is adjoint to the identity map $\sT_\P \ast \P \rightarrow \sT_\P \ast \P$ and hence non-zero. On the other hand we also have the standard distinguished triangle 
\begin{equation} \label{eq:seconddist}
\sT_\P[-1] \rightarrow \P \ast \P_R \rightarrow \O_\Delta.
\end{equation}
If the map $\sT_\P[-1] \rightarrow \P * \P_R $ in this triangle were zero then we would have $\O_\Delta \cong (\P \ast \P_R) \oplus \sT_\P$, which is impossible since $\Hom(\O_\Delta, \O_\Delta) \cong \C $.

We will now show that 
$$\Hom(\sT_\P[-1], \P \ast \P_R) \cong \C$$
and hence the two maps $\sT_\P[-1] \rightarrow \P * \P_R $ in (\ref{eq:firstdist}), (\ref{eq:seconddist}) must be equal up to a non-zero multiple. This implies that their cones are isomorphic and hence $\sT_\P \ast (\sT_\P)_L \cong \O_\Delta$.  A similar argument (starting from $ \P_R * \sT_P \cong \P_R[k+1] $) shows that $ (\sT_\P)_L * \sT_\P \cong \O_\Delta $.  Hence $\Phi_{\sT_\P}$ is an equivalence. 

To show $\Hom(\sT_\P[-1], \P \ast \P_R) \cong \C$ we apply $\Hom(-, \P \ast \P_R)$ to the standard distinguished triangle $\P \ast \P_R \rightarrow \O_\Delta \rightarrow \sT_\P$ to obtain the long exact sequence
\begin{equation*}
\rightarrow \Hom(\O_\Delta, \P * \P_R) \rightarrow \Hom(\P \ast \P_R, \P \ast \P_R) \rightarrow \Hom(\sT_\P[-1], \P \ast \P_R) \rightarrow \Hom(\O_\Delta[-1], \P \ast \P_R) \rightarrow
\end{equation*}
But,
$$\Hom(\O_\Delta, \P \ast \P_R) \cong \Hom (\P_L, \P_R) \cong \Hom (\P_L, \P_L[k]\{l\}) \cong \Hom (\P,\P[k]\{l\}) = 0$$ and
$$\Hom(\O_\Delta, \P * \P_R) = 0 \text{ and }\Hom(\O_\Delta[-1], \P \ast \P_R) \cong \Hom (\P, \P[k+1]\{l\}) = 0$$
so it suffices to show that $\Hom(\P \ast \P_R, \P \ast \P_R) \cong \C$. To do this consider the distinguished triangle $(\sT_\P)_L \rightarrow \O_\Delta \rightarrow \P \ast \P_L$ and apply $\P_R \ast$. Now
\begin{equation*}
\P_R \ast (\sT_\P)_L \cong \P_L \ast (\sT_\P)_L [k]\{l\} \cong (\sT_\P \ast \P)_L [k]\{l\} \cong (\P [k+1]\{l\})_L [k]\{l\} = \P_L[-1]
\end{equation*}
so we get the distinguished triangle
$$\P_L[-1] \rightarrow \P_R \rightarrow \P_R \ast \P \ast \P_L.$$
Applying the functor $\Hom(\P_R, -)$ we get the long exact sequence
\begin{equation*}
\begin{gathered}
\dots \rightarrow \Hom(\P_R, \P_R[k]\{l\}) \rightarrow \Hom(\P_R, \P_R \ast \P \ast \P_L[k]\{l\}) \rightarrow \\ \quad \Hom(\P_R, \P_L[k]\{l\}) \rightarrow \Hom(\P_R, \P_R[k+1]\{l\}) \rightarrow \dots
\end{gathered}
\end{equation*}
Since $\P_L[k]\{l\} \cong \P_R$ and $\Hom(\P_R,\P_R[i]\{l\}) \cong \Hom(\P,\P[i]\{l\}) = 0$ for $i=k,k+1$ we find that 
$$\Hom(\P \ast \P_R, \P \ast \P_R) \cong \Hom (\P_R, \P_R \ast \P \ast \P_R) \cong \Hom(\P_R, \P_R) \cong \Hom(\P, \P) \cong \C.$$
where the second isomorphism follows from the long exact sequence. 
\end{proof}

\begin{Lemma} \label{th:twodist}
Let $ A,A', B, C $ be objects of some triangulated category $ \mathcal{D} $.

Suppose that there are two distinguished triangles
\begin{equation*}
A \xrightarrow{\nu}  B \xrightarrow{\phi} C \xrightarrow{\rho} A[1] 
\qquad C \xrightarrow{\psi} B \rightarrow A' \rightarrow C[1]
\end{equation*}
in $ \mathcal{D} $ such that the composition $ \phi \circ \psi $ is a nonzero multiple of the identity.  Then $B \cong A \oplus C$, and $A \cong A'$.  
\end{Lemma}

\begin{proof}
Without loss of generality we may assume $\phi \circ \psi$ is the identity. The composition $C \xrightarrow{\psi} B \xrightarrow{\phi} C \xrightarrow{\rho} A[1] $ is the zero map since $\rho \circ \phi = 0$. But $\phi \circ \psi = id$ which means $\rho$ must be the zero map. Since $B[1] \cong \Cone(\rho)$ this means there is an isomorphism $\beta: B \cong A \oplus C$ making the following diagram commute
\begin{equation*}
\begin{CD}
A @>{\nu}>> B @>{\phi}>> C \\
@V{\mbox{id}}VV @V{\beta}VV @V{\mbox{id}}VV \\
A @>{i_1}>> A \oplus C @>{p_2}>> C
\end{CD}
\end{equation*}
where $i_1,i_2$ are the standard injections of $A$ and $C$ into $A \oplus C$ and $p_1,p_2$ the standard projections from $A \oplus C$ to $A$ and $C$ respectively. 

To show that $A \cong A'$ we will construct an isomorphism $\theta: B \rightarrow A \oplus C$ such that the following square commutes
\begin{equation*}
\begin{CD}
C @>{\psi}>> B \\
@V{\mbox{id}}VV @V{\theta}VV \\
C @>{i_2}>> A \oplus C
\end{CD}
\end{equation*}

Showing this implies there is an isomorphism $A' = \Cone(C \xrightarrow{\psi} B)  \rightarrow \Cone(C \xrightarrow{i_2} A \oplus C) = A$. We define
$$\theta = \beta - i_1 \circ p_1 \circ \beta \circ \psi \circ \phi.$$
Then
\begin{equation*}
\begin{aligned}
\theta \circ \psi &= \beta \circ \psi - i_1 \circ p_1 \circ \beta \circ \psi \circ \phi \circ \psi \\
&= (\mbox{id} - i_1 \circ p_1) \circ \beta \circ \psi \\
&= i_2 \circ p_2 \circ \beta \circ \psi \\
&= i_2 \circ \phi \circ \psi = i_2\
\end{aligned}
\end{equation*}
where we use $\phi \circ \psi = \mbox{id}$ twice. This shows the diagram commutes. 

The inverse of $\theta$ is $\theta^{-1} = \nu \circ p_1 + \psi \circ p_2$ since 
\begin{equation*}
\theta \circ (\nu \circ p_1 + \psi \circ p_2) = \theta \circ \nu \circ p_1 + \theta \circ \psi \circ p_2 
= i_1 \circ p_1 + i_2 \circ p_2 = \text{id}.
\end{equation*}
\end{proof}

\subsection{Braid relations among twists} \label{se:twistbraid}
The following results about braid relations among twists are due to Anno \cite{A}.  We have adapted them for use in our setting.

The following result generalizes Lemma 8.21 of \cite{H}, which is a version of Lemma 2.11 of \cite{ST}.
\begin{Lemma} \label{th:twists}
Let $ \T_{\P} : D(X) \rightarrow D(Y) $ be a spherical functor, and let $ \Phi_\Q $ denote any autoequivalence of $ D(Y) $.  Then there is an isomorphism of functors
\begin{equation*}
\Phi_\Q \circ \T_{\P} \cong \T_{\Q \ast \P} \circ \Phi_\Q.
\end{equation*}
\end{Lemma}

\begin{proof}
Let $ \mathcal{R} = \Q \ast \sT_{\P} \ast \Q_R $.  By the relationship between FM kernels and transforms, it suffices to show that $ \mathcal{R} \cong \sT_{\Q \ast \P}$.

By definition,
\begin{equation*}
\sT_\P := \Cone(\P \ast \P_R \xrightarrow{\beta_\P} \O_\Delta)
\end{equation*}
and hence
\begin{equation*}
\mathcal{R} = \Cone(\Q \ast \P \ast \P_R \ast \Q_R \xrightarrow{\Q * \beta_\P * \Q_R} \Q \ast \O_\Delta \ast \Q_R).
\end{equation*}

Since $ (\Q \ast \P)_R \cong \P_R \ast \Q_R $ and since $ \Q \ast \O_\Delta \ast \Q_R \cong \O_\Delta$, we have that
\begin{equation*}
\mathcal{R} \cong \Cone((\Q \ast \P) \ast (\Q \ast \P)_R \rightarrow \O_\Delta).
\end{equation*}

Moreover by Lemma \ref{th:convbeta}, this morphism is $ \beta_{\Q * \P} $.  Thus, $ \mathcal{R} \cong \sT_{\Q\ast \P} $ as desired. 
\end{proof}

\begin{Theorem} \label{th:braidrel}
Suppose that $ \Phi_\P $ and $ \Phi_{\P'} $ are both spherical functors from $ D(X) \rightarrow D(Y) $.  Suppose also that $\T_{\sT_{\P'} \ast \P} \cong \T_{{\sT_\P}_R \ast \P'} $.  Then $ \T_{\P} $ and $ \T_{\P'} $ satisfy the braid relation
\begin{equation*}
\T_\P \circ \T_{\P'} \circ \T_\P \cong \T_{\P'} \circ \T_\P \circ \T_{\P'}.
\end{equation*}
\end{Theorem}

\begin{proof}
Using Lemma \ref{th:twists} in the first step and the hypothesis in the second step and the lemma again in the third step, we have
\begin{equation*}
\begin{aligned}
 \T_\P \circ \T_{\P'} \circ \T_\P &\cong \T_\P \circ \T_{\sT_{\P'} \ast \P} \circ \T_{\P'} \\
 &\cong \T_\P \circ \T_{{\sT_{\P}}_R \ast \P'} \circ \T_{\P'} 
 \cong \T_{\sT_{\P} \ast {\sT_{\P}}_R \ast \P'} \circ \T_\P \circ \T_{\P'}
 \cong \T_{\P'} \circ \T_\P \circ \T_{\P'}
 \end{aligned}
 \end{equation*}
\end{proof}

\section{Functors from tangles} \label{se:funtangle}

An $(n,m) $ \textbf{tangle} is a proper, smooth embedding of $ (n+m)/2 $ arcs and a finite number of circles into $ \mathbb{R}^2 \times [0,1] $ such that the boundary points of the arcs map bijectively on the $n+m $ points $ (1,0,0), \dots, (n,0,0), (1,0,1), \dots, (m,0,1) $.  A $ (0,0) $ tangle is a \textbf{link}.

Given an $(n,m)$ tangle $ T $ and a $ (m,p) $ tangle $ U$, there is a composition tangle $ T \circ U $, which is the $ (n,p) $ tangle made by stacking $ U $ on top of $ T $ and the shrinking the $z$-direction.

\subsection{Generators and Relations}

\begin{figure}
\begin{center}
\psfrag{cup}{cup} 
\psfrag{cap}{cap} 
\psfrag{crossing 1}{crossing \#1} 
\psfrag{crossing 2}{crossing \#2} 
\psfrag{crossing 3}{crossing \#3} 
\psfrag{crossing 4}{crossing \#4} 
\puteps[0.35]{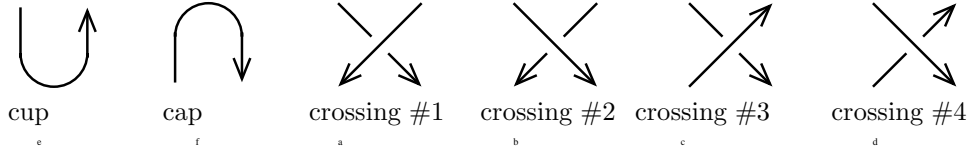}
\end{center}
\caption{(half) the generators for tangle diagrams.}\label{f1}
\end{figure}

By projecting to $\mathbb{R} \times [0,1]$ from a generic point we can represent any tangle as a planar diagram. By scanning the diagram of an $(n,m)$ tangle from top to bottom we can decompose it as the composition of cups, caps and crossings. The list of all such building blocks consists of those tangles from Figure \ref{f1} together with those obtained from them by switching all the orientations (arrows). 

We will give these generators names depending on their position.  A cap creating the $ i $ and $ i+1 $ strands in an $ n $ strand tangle will be denoted $ g_n^i $.  A cup connecting the $ i $ and $ i+1 $ strands in an $ n $ strand tangle will be denoted by $ f_n^i $.  A crossing of the $ i $ and $ i+1 $ strands in an $ n $ strand tangle will be denoted $ t_n^i(l) $, where $ l $ varies from 1 to 4, depending on the type of crossing, as shown in Figure \ref{f1}. Any generator which is obtained from these generators by reversing the direction of all arrows involved will be denoted by the same symbol.

The following theorem tells us when two tangle diagrams represent isotopic tangles.

\begin{Lemma}[{\cite[Lemma X.3.5]{K}}]\label{lem:relations} Two tangle diagrams represent isotopic tangles if and only if one can be obtained from the other by applying a finite number of the following operations:
\begin{itemize}
\item a Reidemeister move of type (0),(I),(II) or (III).
\item an isotopy exchanging the order with respect to height of two caps, cups, or crossings (e.g. the left figure in \ref{f3} shows such an isotopy involving a cup and a cap). 
\item the rightmost two isotopies in figure \ref{f3}, which we call the pitchfork move.
\end{itemize}
More concisely, we have the following relations, plus those obtained by changing the directions of the strands:
\begin{enumerate}
\item Reidemeister (0) : $ f_n^i \circ g_n^{i+1} = id = f_n^{i+1} \circ g_n^i $ 
\item Reidemeister (I) : $f_n^i \circ t_n^{i \pm 1}(2) \circ g_n^i = id = f_n^i \circ t_n^{i \pm 1}(1) \circ g_n^i$ 
\item Reidemeister (II) : $ t_n^i(2) \circ t_n^i(1) = id = t_n^i(1) \circ t_n^i(2) $ 
\item Reidemeister (III) : $ t_n^i(2) \circ t_n^{i+1}(2) \circ t_n^i(2) = t_n^{i+1}(2) \circ t_n^i(2) \circ t_n^{i+1}(2) $ 
\item cap-cap isotopy : $ g_{n+2}^{i+k} \circ g_n^i = g_{n+2}^i \circ g_n^{i+k-2} $ 
\item cup-cup isotopy : $ f_n^{i+k-2} \circ f_{n+2}^i = f_n^i \circ f_{n+2}^{i+k} $
\item cup-cap isotopy : $ g_n^{i+k-2} \circ f_n^i = f_{n+2}^i \circ g_{n+2}^{i+k},  \quad g_n^i \circ f_n^{i+k-2} = f_{n+2}^{i+k} \circ g_{n+2}^i $
\item cap-crossing isotopy: $ g_n^i \circ t_{n-2}^{i+k-2}(l) = t_n^{i+k}(l) \circ g_n^i, 
\quad g_n^{i+k} \circ t_{n-2}^i(l) = t_n^i(l) \circ g_n^{i+k} $
\item cup-crossing isotopy: $ f_n^i \circ t_n^{i+k}(l) = t_{n-2}^{i+k-2}(l) \circ f_n^i, \quad 
f_n^{i+k} \circ t_n^i(l) = t_{n-2}^i(l) \circ f_n^{i+k} $
\item crossing-crossing isotopy: $ t_n^i(l) \circ t_n^j(m) = t_n^j(m) \circ t_n^i(l) $
\item pitchfork move : $ t_{n}^i(1) \circ g_n^{i+1} = t_{n}^{i+1}(4) \circ g_n^i, \quad t_{n}^i(2) \circ g_n^{i+1} = t_{n}^{i+1}(3) \circ g_n^i $
\end{enumerate}
where in each case $ k \ge 2 $ and $ 1 \le l \le 4 $.
\end{Lemma}

\begin{figure}
\begin{center}
\psfrag{move (0)}{R-move (0)} 
\psfrag{move (I)}{R-move (I)} 
\psfrag{move (II)}{R-move (II)} 
\psfrag{move (III)}{R-move (III)} 
\psfrag{sim}{$\sim$} 
\puteps[0.20]{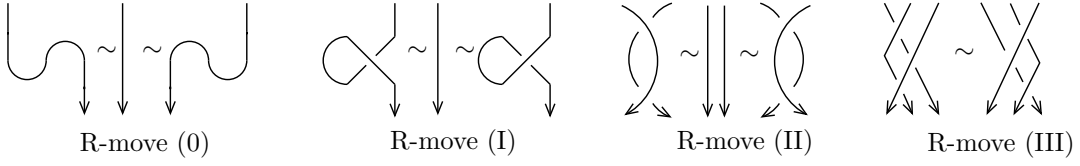}
\end{center}
\caption{Reidemeister relations for tangle diagrams.}\label{f2}
\end{figure}

\begin{figure}
\begin{center}
\psfrag{sim}{$\sim$} 
\puteps[0.25]{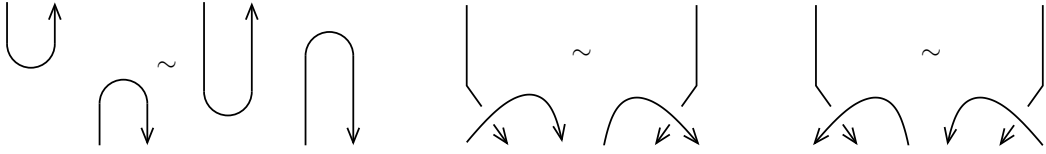}
\end{center}
\caption{Relations for tangle diagrams.}\label{f3}
\end{figure}

\subsection{The Functor $ \Psi(T): D(Y_n) \rightarrow D(Y_m) $}\label{sse:functor}
To each tangle $(n,m) $ tangle $ T $, we will construct an isomorphism class of functor $ \Psi(T) : D(Y_n) \rightarrow D(Y_m) $.  These functors will satisfy the property that $ \Psi(T) \circ \Psi(U) \cong \Psi(T \circ U) $. We begin by defining functors for each of the elementary tangles in figure \ref{f1}.  The tangle obtained from a tangle in figure \ref{f1} by reversing all the strand orientations will be assigned the same functor. 

\subsubsection{Cups and caps} \label{se:capcup}
Recall that we have the equivariant subvariety $ X_n^i $ of $ Y_n $ which projects to $ Y_{n-2}$.  Thus, we may regard $ X_n^i $ as a subvariety of the product $ Y_{n-2} \times Y_n $.  Let $ \sG_n^i $ denote the $\C^\times$-equivariant sheaf on $ Y_{n-2} \times Y_n $ defined by 
\begin{equation*} 
\sG_n^i := \O_{X_n^i} \otimes \E'_i \{-i+1\}.
\end{equation*}

Now, we define the functor $ \G_n^i : D(Y_{n-2}) \rightarrow D(Y_n) $ to be the Fourier-Mukai transform with respect to the kernel $ \sG_n^i $. We will use this functor for the cap, so we define 
\begin{equation} \label{eq:Gdef}
\Psi(g_n^i) := \G_n^i.
\end{equation}

We can give a different description of the functor $ \G_n^i $.  Namely, $\G_n^i(\cdot) = i_\ast(q^\ast(\cdot) \otimes \E_i\{-i+1\})$ where we make use of the diagram 
\begin{equation*}\begin{CD}X_n^i @>i>> Y_n \\@VqVV \\ Y_{n-2}\end{CD} \end{equation*}

Now, we define $\F_n^i : Y_n \rightarrow Y_{n-2} $ by $\F_n^i(\cdot) = q_\ast(i^\ast(\cdot) \otimes \E_{i+1}^\vee)\{i\} $ and we define 
\begin{equation} \label{eq:Fdef}
\Psi(f_n^i) := \F_n^i.
\end{equation} 
As with $ \G_n^i $, the functor $ \F_n^i $ can also be described as a Fourier-Mukai transform with respect to the kernel $\sF_n^i := \O_{X_n^i} \otimes \E_{i+1}^\vee \{i\}$.

\subsubsection{Crossings} \label{se:crossing}
Consider the $\C^\times$-equivariant subvariety
\begin{equation*}
 Z_n^i := \{ (L_\cdot, L'_\cdot) : L_j = L_j' \text{ for $ j \ne i $ } \} \subset Y_n \times Y_n.
\end{equation*}
$Z_n^i$ has two smooth irreducible components, each of dimension $ n $.   The first component corresponds to the locus of points where $L_i=L'_i$, and so is the diagonal $ \Delta \subset Y_n \times Y_n $.   The second component is the closure of the locus of points where $L_i \ne L'_i$.  Note that if $ L_i \ne L_i' $, then $ zL_{i+1} \subset L_i \cap L_i' =  L_{i-1}$.  Thus we see that on this closure, $ zL_{i+1} = L_{i-1} $.  Hence this second component is the subvariety $\othercomp := X_n^i \times_{Y_{n-2}} X_n^i$ where the fibre product is with respect to the map $q: X_n^i \rightarrow Y_{n-2}$. Thus $Z_n^i = \Delta \cup \othercomp \subset Y_n \times Y_n$. 

To a crossing connecting boundary points $i$ and $i+1$ we assign a Fourier-Mukai kernel $ \mathcal{T}_n^i(l) \in  D(Y_n \times Y_n)$ according to the type of crossing: 
\begin{itemize}
\item crossing \#1: $ \sT_n^i(1) := \O_{Z_n^i}[1]\{-1\}$ 
\item crossing \#2: $ \sT_n^i(2) := \O_{Z_n^i} \otimes \E_{i+1}^\vee \otimes \E'_i [-1] \{3\}$
\item crossing \#3: $ \sT_n^i(3) := \O_{Z_n^i}[-1]\{2\} $
\item crossing \#4: $ \sT_n^i(4) := \O_{Z_n^i} \otimes \E_{i+1}^\vee \otimes \E'_i [1]$
\end{itemize}
Now that we have these kernels, we associate to each crossing a functor 
\begin{equation}\label{eq:Cdef}
\Psi(t_n^i(l)) := \T_n^i(l) := \Phi_{\sT_n^i(l)}
\end{equation}
where as usual $ l $ runs from 1 to 4. 

\subsubsection{Functor for a tangle}
Let $ T $ be a tangle.  Scanning a projection of $ T $ from top to bottom and composing along the way gives us a functor $ \Psi(T) : D(Y_n) \rightarrow D(Y_m) $.  However, this functor may depend on the choice of tangle projection.  

\begin{Theorem}\label{thm:main}
The isomorphism class of the functor $\Psi(T): D(Y_n) \rightarrow D(Y_m)$ associated to the planar diagram of an $(n,m)$ tangle $T$ is a tangle invariant. 
\end{Theorem}

To prove this theorem, we must check that the functors assigned to the elementary tangles satisfy the relations from Lemma \ref{lem:relations}.  This will be done in section \ref{se:invariance}.

This construction associates to any link $L$ a functor $\Psi(L): D(Y_0) \rightarrow D(Y_0)$. Since $Y_0$ is a point with the trivial action of $\C^\times$, $\Psi_L$ is determined by $\Psi(L)(\C) \in D(Y_0)$ which is a complex of graded vector spaces. We denote by $\Kh^{i,j}(L)$ the $j$-graded piece of the $i$th cohomology group of $\Psi(L)(\C)$ (so $i$ marks the cohomological degree and $j$ marks the graded degree). Since $\Psi(L)$ is a tangle invariant $\Kh^{i,j}(L)$ is an invariant of the link $L$.

\subsection{Properties of the kernels}

Before continuing we mention a few technical results which we will need later. 
\begin{Lemma}\label{lem:facts} If $i$ is the equivariant inclusion $X_n^i \rightarrow Y_n$ and $q$ the equivariant projection $X_n^i \rightarrow Y_{n-2}$ then
\begin{enumerate}
\item $\O_{Y_n}(X_n^i) \cong \E_{i+1}^\vee \otimes \E_i\{2\}$,
\item $\omega_{X_n^i} \otimes i^\ast \omega^\vee_{Y_n} \cong \O_{X_n^i}(X_n^i) \cong i^\ast(\E_{i+1}^\vee \otimes \E_i)\{2\} \cong \omega_{X_n^i} \otimes q^\ast \omega^\vee_{Y_{n-2}}\{2\}$ so that, in particular, $i^\ast \omega_{Y_n} \cong q^\ast \omega_{Y_{n-2}}\{2\}$,
\item $\omega_{\othercomp} \cong \O_{\othercomp} \otimes \E_i \otimes \E_{i+1}^\vee \otimes \E'_i \otimes \E_{i+1}^{'\vee} \otimes \pi_1^\ast \omega_{Y_n} \{2\} \cong \O_{\othercomp} \otimes \E_i \otimes \E_{i+1}^\vee \otimes \E'_i \otimes \E_{i+1}^{'\vee} \otimes \pi_2^\ast \omega_{Y_n} \{2\}$ as $\C^\times$-equivariant sheaves on $Y_n \times Y_n$, 
\item $\O_{\othercomp}(D) \cong \O_{\othercomp} \otimes \E_i^\vee \otimes \E'_{i+1}$ where $D$ is the divisor $\othercomp \cap \Delta$ inside $\othercomp$ .
\end{enumerate}
\end{Lemma}
\begin{proof}
The map $z: L_{i+1} \rightarrow L_i \{2\}$ induces a morphism $L_{i+1}/L_i \rightarrow L_i/L_{i-1}\{2\}$ and hence a section of $\Hom_{Y_n}(\E_{i+1},\E_i\{2\})$. This section is zero precisely over the locus where $z$ maps $L_{i+1}$ to $L_{i-1}$, namely $X_n^i$. Thus $\O_{Y_n}(X_n^i) \cong \E_{i+1}^\vee \otimes \E_i\{2\}$ as long as the section is transverse to the zero section. Since $X_n^i$ is smooth it suffices to show that $\E_{i+1}^\vee \otimes \E_i$ is not a multiple of $\O(X_n^i)$. This would follow if we can show they both restrict to the same (non-trivial) line bundle on a subvariety. We take this subvariety to be a $\p$ fibre of $q: X_n^{i+1} \rightarrow Y_{n-2}$. Since $X_n^i$ and $X_n^{i+1}$ intersect transversely in a section of $q$ (see the remark following lemma \ref{lem:transverse}) the restriction of $X_n^i$ to the fibre is $\O_\p(1)$. Meanwhile, $\E_{i+1}^\vee$ restricts to the dual of the tautological line bundle, namely $\O_\p(1)$, whereas $\E_i$ restricts to the trivial line bundle $\O_\p$ so that $\E_{i+1}^\vee \otimes \E_i$ also restricts to $\O_\p(1)$. This concludes (i).

For (ii), by the equivariant adjunction formula we get $\omega_{X_n^i} \otimes i^\ast \omega^\vee_{Y_n} \cong \O_{X_n^i}(X_n^i)$ which is isomorphic to $i^\ast(\E_{i+1}^\vee \otimes \E_i) \{2\}$ by (i). To show the last equality note that $q: X_n^i \rightarrow Y_{n-2}$ is the $\C^\times$-equivariant $\p$ bundle $\mathbb{P}(z^{-1}L_{i-1}/L_{i-1}) \rightarrow Y_{n-2}$. So the relative dualizing sheaf $\omega_{X_n^i} \otimes q^\ast \omega^\vee_{Y_{n-2}}$ of $q$ is 
$$\Hom(z^{-1}L_{i-1}/L_i, L_i/L_{i-1}) = i^\ast(\E_{i+1}^\vee \otimes \E_i).$$

For (iii), consider the projection $\pi_1: \othercomp \rightarrow X_n^i$ which is the $\p$ bundle $\mathbb{P}(z^{-1}L_{i-1}/L_{i-1}) \rightarrow X_n^i$. The relative dualizing sheaf is then $\Hom(z^{-1}L_{i-1}/L'_i,L'_i/L_{i-1}) = \E_{i+1}^{'\vee} \otimes \E'_i$. At the same time this is isomorphic to $\omega_{\othercomp} \otimes \pi_1^\ast \omega^\vee_{X_n^i}$ which from the calculations above is isomorphic to $\omega_{\othercomp} \otimes \E_{i+1} \otimes \E_i^\vee \otimes \pi_1^\ast \omega^\vee_{Y_n} \{-2\}$. It follows that 
$$\omega_{\othercomp} \cong \E_i \otimes \E_{i+1}^\vee \otimes \E_i' \otimes \E_{i+1}^{\vee'} \otimes \pi_1^\ast \omega_{Y_n} \{2\}$$
while the second isomorphism involving $\omega_{\othercomp}$ follows similarly.

On $\othercomp$ consider the equivariant inclusion map $L_i/L_{i-1} \rightarrow L_{i+1}/L_{i-1} = L_{i+1}/L'_{i-1}$ composed with the equivariant projection $L_{i+1}/L'_{i-1} \rightarrow L_{i+1}/L'_i$. This gives a $\C^\times$-equivariant map $L_i/L_{i-1} \rightarrow L_{i+1}/L'_i$ which vanishes precisely along the locus $D$ where $L_i = L'_i$. If this section is transverse to the zero section this implies $\O_{\othercomp}(D) \cong \O_{\othercomp} \otimes \E_i^\vee \otimes \E'_{i+1}$. Since $D$ is smooth, to check the section is transverse it suffices to show that $\E_i^\vee \otimes \E'_{i+1}$ is not a multiple of $D$. This would follow if we can show they both restrict to the same (non-trivial) line bundle on a subvariety. We take this subvariety to be a $\p$ fibre of $\pi_1: \othercomp \rightarrow X_n^i$. Then $D$ restricts to $\O_\p(1)$ on such a fibre while $\E_i^\vee$ and $\E'_{i+1}$ restrict to $\O_\p(1)$ and $\O_\p$ respectively.  This concludes (iv).
\end{proof}

The functors $ \F_n^i $ and $ \G_n^i $ for caps and cups are related by the following Lemma.
\begin{Lemma} \label{th:Gadj} 
We have that $ \sF_n^i \cong {\sG_n^i}_R[1]\{-1\} \cong {\sG_n^i}_L[-1]\{1\}.$
In particular, $$\F_n^i(\cdot) = q_\ast(i^\ast(\cdot) \otimes \E_{i+1}^\vee)\{i\} \cong {\G_n^i}^R(\cdot)[1]\{-1\} \cong {\G_n^i}^L(\cdot)[-1]\{1\}.$$  
\end{Lemma}
\begin{proof}
The right adjoints of $q^\ast$, $i_\ast$ and $(\cdot) \otimes \E_i$ are $q_\ast$, $i^!$ and $(\cdot) \otimes \E_i^\vee$ respectively. Thus ${\G_n^i}^R(\cdot) = q_\ast(i^!(\cdot) \otimes \E_i^\vee)\{i-1\}$. Since $i: X_n^i \rightarrow Y_n$ is codimension one we find $i^!(\cdot) = i^\ast(\cdot) \otimes \omega_{X_n^i} \otimes i^\ast \omega^\vee_{Y_n}[-1]$ so that ${\G_n^i}^R(\cdot) = q_\ast(i^\ast(\cdot) \otimes \omega_{X_n^i} \otimes i^\ast \omega^\vee_{Y_n} \otimes \E_i^\vee)[-1]\{i-1\}$.  Finally by lemma \ref{lem:facts} below we see that $\omega_{X_n^i} \otimes i^\ast \omega^\vee_{Y_n} \otimes \E_i^\vee \cong \E_{i+1}^\vee \{2\}$ so 
$${\G_n^i}^R(\cdot)[1]\{-1\} \cong q_\ast(i^\ast(\cdot) \otimes \E_{i+1}^\vee)\{i\}.$$  

The calculation of the left adjoint ${\G_n^i}^L$ follows similarly from the facts that the left adjoints of $q^\ast$ and $i_\ast$ are $q_!$ and $i^\ast$. This time, using that $q: X_n^i \rightarrow Y_{n-2}$ a $\p$ bundle, we have $q_!(\cdot) = q_\ast((\cdot) \otimes \omega_{X_n^i} \otimes q^\ast \omega^\vee_{Y_{n-2}})[1]$ so that 
${\G_n^i}^L(\cdot) = q_\ast(i^\ast(\cdot) \otimes \omega_{X_n^i} \otimes q^\ast \omega^\vee_{Y_{n-2}} \otimes \E_i^\vee)[1]\{i-1\}.$ 
By lemma \ref{lem:facts} we see that $\omega_{X_n^i} \otimes q^\ast \omega^\vee_{Y_{n-2}} \otimes \E_i^\vee \cong \E_{i+1}^\vee$ so 
$${\G_n^i}^L(\cdot)[-1]\{1\} \cong q_\ast(i^\ast(\cdot) \otimes \E_{i+1}^\vee)\{i\}.$$

The statement on the level of kernels follows similarly.
\end{proof}

The functors for under and over crossings are related by the following Lemma.
\begin{Lemma} \label{th:twistadjoint}
$ {\sT_n^i(2)}_L \cong \sT_n^i(1) $
\end{Lemma}
\begin{proof}
Consider the standard $\C^\times$-equivariant short exact sequence $0 \rightarrow \O_\Delta(-D) \rightarrow \O_{Z_n^i} \rightarrow \O_{\othercomp} \rightarrow 0$ where $D = \othercomp \cap \Delta$ is a divisor both in $\othercomp$ and $\Delta$. Since $ D \subset \Delta $ is just $ X_n^i \subset Y_n $, we see that $\O_\Delta(D) $ is the restriction of the globally defined line bundle $\E_{i+1}^{\vee} \otimes \E'_i \{2\}$ on $Y_n \times Y_n$. Dualizing the above short exact sequence we get the distinguished triangle
$$\O^\vee_{\othercomp} \rightarrow \O^\vee_{Z_n^i} \rightarrow \O^\vee_\Delta(D) $$
where the connecting morphism $\O^\vee_\Delta(D)[-1] \rightarrow \O^\vee_{\othercomp}$ is non-zero. Now 
$$\O^\vee_\Delta \cong \omega_\Delta \otimes \omega_{Y_n \times Y_n}^\vee[-\mbox{codim}(\Delta)] \cong \O_\Delta \otimes \pi_2^\ast \omega^\vee_{Y_n}[-n]$$
so that $\O^\vee_\Delta(D) \cong \O_\Delta \otimes \E_{i+1}^{\vee} \otimes \E'_i \otimes \pi_2^\ast \omega^\vee_{Y_n}[-n] \{2\}$. Similarly, using lemma \ref{lem:facts}
\begin{eqnarray*}
\O^\vee_{\othercomp} &\cong& \omega_{\othercomp} \otimes \omega_{Y_n \times Y_n}^\vee[-n] \\
&\cong& \O_{\othercomp} \otimes \E_i \otimes \E^\vee_{i+1} \otimes \E'_i \otimes \E^{'\vee}_{i+1} \otimes \pi_1^\ast \omega_{Y_n} \otimes \omega_{Y_n \times Y_n}^\vee [-n] \{2\} \\
&\cong& \O_{\othercomp} \otimes \E_i \otimes \E^\vee_{i+1} \otimes \E'_i \otimes \E^{'\vee}_{i+1} \otimes \pi_2^\ast \omega^\vee_{Y_n}[-n] \{2\} \\
\end{eqnarray*}
Tensoring the above distinguished triangle by $\E_{i+1} \otimes \E_i^{'\vee} \otimes \pi_2^\ast \omega_{Y_n}[n] \{-2\}$ and simplifying we get 
$$\O_{\othercomp} \otimes \E_i \otimes \E_{i+1}^{'\vee} \rightarrow \O^\vee_{Z_n^i} \otimes \E_{i+1} \otimes \E_i^{'\vee} \otimes \pi_2^\ast \omega_{Y_n}[n] \{-2\} \rightarrow \O_\Delta$$
Using lemma \ref{lem:facts} we can simplify further to obtain 
$$\O_{\othercomp}(-D) \rightarrow \O^\vee_{Z_n^i} \otimes \E_{i+1} \otimes \E_i^{'\vee} \otimes \pi_2^\ast \omega_{Y_n}[n] \{-2\} \rightarrow \O_\Delta$$
where the connecting morphism $\O_\Delta[-1] \rightarrow \O_{\othercomp}(-D)$ is non-zero. By corollary \ref{cor:abstractnonsense3} there exists, up to a non-zero multiple, a unique such map, namely the one coming from the standard sequence $0 \rightarrow \O_{\othercomp}(-D) \rightarrow \O_{Z_n^i} \rightarrow \O_\Delta \rightarrow 0$. Hence $\O^\vee_{Z_n^i} \otimes \E_{i+1} \otimes \E_i^{'\vee} \otimes \pi_2^\ast \omega_{Y_n}[n]\{-2\} \cong \O_{Z_n^i}$ which means 
$$ \O^\vee_{Z_n^i} \cong \O_{Z_n^i} \otimes \E_{i+1}^\vee \otimes \E'_i \otimes \pi_2^\ast \omega^\vee_{Y_n}[-n]\{2\}.$$
Consequently
\begin{eqnarray*}
\sT_n^i(2)_L &\cong& (\O_{Z_n^i} \otimes \E_{i+1}^\vee \otimes \E'_i [-1] \{3\})^\vee \otimes \pi_2^\ast \omega_{Y_n}[n] \\
&\cong& \O_{Z_n^i} \otimes \E_{i+1}^\vee \otimes \E'_i \otimes \pi_2^\ast \omega^\vee_{Y_n}[-n] \{2\} \otimes \E_{i+1} \otimes \E_i^{'\vee} [1] \{-3\} \otimes \pi_2^\ast \omega_{Y_n}[n] \\
&\cong& \O_{Z_n^i}[1]\{-1\} = \sT_n^i(1)
\end{eqnarray*}
\end{proof}

It turns out that the functors associated to crossings can also be described as spherical twists. The essential reason for this is that $ Z_n^i $ has two components, one of which is the diagonal $ \Delta $ and the other is the fibre product $ \othercomp = X_n^i \times_{Y_{n-2}} X_n^i $. 

\begin{Theorem} \label{th:kerneltwist}
The functor $ \T_n^i(2) $ is the twist in the functor $ \G_n^i $ shifted by $[-1]\{1\}$. In particular
\begin{equation*}
\sT_n^i(2) \cong \sT_{\sG_n^i}[-1]\{1\}.
\end{equation*}
Similarly, the kernel of the functors $\T_n^i(1)$, $\T_n^i(3)$ and $\T_n^i(4)$ are given by 
$$\sT_n^i(1) \cong (\sT_{\sG_n^i})_L[1]\{-1\} \text{ and }
\sT_n^i(3) \cong (\sT_{\sG_n^i})_L[-1]\{2\} \text{ and }
\sT_n^i(4) \cong \sT_{\sG_n^i}[1]\{-2\}.$$
\end{Theorem}
\begin{proof} 
We need to show that 
$$\O_{Z_n^i} \otimes \E_{i+1}^\vee \otimes \E'_i \{2\} \cong \Cone(\P \ast \P_R \rightarrow \O_\Delta) \in D(Y_{n-2} \times Y_n)$$
where $\P = \O_{X_n^i} \otimes \E_i' \{-i+1\}$.  By Lemma \ref{th:Gadj}, we have $ \P_R \cong \O_{X_n^i} \otimes \E_{i+1}^{\vee}[-1] \{i+1\}$.

Now we compute $\P \ast \P_R \in D(Y_n \times Y_{n-2} \times Y_n)$. We have
\begin{equation*}
\P \ast \P_R = \pi_{13 \ast} (\pi_{12}^\ast \O_{X_n^i} \otimes \E_{i+1}^\vee[-1]\{i+1\} \otimes \pi_{23}^\ast \O_{X_n^i} \otimes \E_i^{''} \{-i+1\})
\end{equation*}
By corollary \ref{cor:transverse} we know $\pi_{12}^{-1}(X_n^i)$ and $\pi_{23}^{-1}(X_n^i)$ intersect transversely. So by lemma \ref{th:tensor} we have $\pi_{12}^\ast \O_{X_n^i} \otimes \pi_{23}^\ast \O_{X_n^i} \cong \O_W$ where $W = \pi_{12}^{-1}(X_n^i)$ and $\pi_{23}^{-1}(X_n^i)$ inside $Y_n \times Y_{n-2} \times Y_n$. The projection $\pi_{13}$ maps $W$ isomorphically onto $\othercomp \subset Y_n \times Y_n$. Using the projection formula we find
$$\P \ast \P_R \cong \O_{\othercomp} \otimes \E_{i+1}^\vee \otimes \E'_i [-1]\{2\}.$$
Since the map $\P \ast \P_R \rightarrow \O_\Delta$ is obtained through a series of adjunctions it is not identically zero. From the distinguished triangle 
$$\O_\Delta \rightarrow \Cone(\P \ast \P_R \rightarrow \O_\Delta) \rightarrow \P \ast \P_R[1]$$
we find that $\Cone(\P \ast \P_R \rightarrow \O_\Delta)$ is a non-trivial extension of $\O_{\othercomp} \otimes \E_{i+1}^\vee \otimes \E'_i \{2\}$ by $\O_\Delta$ which is supported in degree zero. On the other hand, $\othercomp \cap \Delta$ is the divisor $X_n^i \hookrightarrow \Delta$ so we have the standard exact sequence
$$0 \rightarrow \O_\Delta(-X_n^i) \rightarrow \O_{Z_n^i} \rightarrow \O_{\othercomp} \rightarrow 0.$$
By lemma \ref{lem:facts} we know $\O_{Y_n}(-X_n^i) \cong \E_{i+1} \otimes \E_i^\vee \{-2\}$ which means $\O_\Delta(-X_n^i) \cong \O_\Delta \otimes \E_{i+1} \otimes \E_i^{'\vee} \{-2\}$. So
$$0 \rightarrow \O_\Delta \rightarrow \O_{Z_n^i} \otimes \E_{i+1}^\vee \otimes \E'_i \{2\} \rightarrow \O_{\othercomp} \otimes \E_{i+1}^\vee \otimes \E_i' \{2\} \rightarrow 0.$$
Thus $\O_{Z_n^i} \otimes \E_{i+1}^\vee \otimes \E'_i \{2\}$ is a non-trivial extension of $\O_{\othercomp} \otimes \E_{i+1}^\vee \otimes \E'_i \{2\}$ by $\O_{\Delta}$. Since the same is true of $\Cone(\P \ast \P_R \rightarrow \O_\Delta)$ it is enough to show that there exists a unique non-trivial such extension since then 
$$\Cone(\P \ast \P_R \rightarrow \O_\Delta) \cong \O_{Z_n^i} \otimes \E_{i+1}^\vee \otimes \E'_i \{2\}.$$ 
We do this by showing that $\Ext^1(\O_{\othercomp} \otimes \E_{i+1}^\vee \otimes \E'_i \{2\}, \O_\Delta)$ is one dimensional, or equivalently that $\Ext^1(\O_{\othercomp}, \O_\Delta(-X_n^i))$ is one dimensional. Since $\Delta$ and $\othercomp$ are smooth with $\Delta \cap \othercomp = X_n^i$ a smooth divisor in $\Delta$, we know from corollary \ref{cor:abstractnonsense3} that $\Hom^1_{Y_n \times Y_n}(\O_{\othercomp},\O_\Delta(-X_n^i)) \cong \Hom_{X_n^i}(\O_{X_n^i}, \O_{X_n^i})$ is one dimensional. 

To see that $\sT_n^i(1) \cong (\sT_{\sG_n^i})_L[1]\{-1\}$ we use Lemma \ref{th:twistadjoint} above together with the fact that $(\sT_{\P}[-1]\{1\})_L = (\sT_\P)_L[1]\{-1\}$ which we apply to $\P = \sG_n^i$. The expressions for $\sT_n^i(3)$ and $\sT_n^i(4)$ follow similarly. 
\end{proof}

\begin{Lemma}\label{lem:abstractnonsense2} Let $S,S' \subset T$ be smooth $\C^\times$-equivariant subvarieties of a smooth quasi-projective variety $T$ equipped with a $\C^\times$ action. Suppose the scheme theoretic intersection $S \cap S' \subset T$ is smooth. If $\sF$ is a $\C^\times$-equivariant locally free sheaf on $S$ then 
\begin{eqnarray*}
\H^l(\O_{S'}^\vee \otimes i_\ast \sF) =
\left\{
\begin{array}{ll}
0 & \mbox{if $l<c$}  \\
\O_{S \cap S'}(j^\ast \sF \otimes \det N_{S \cap S'/S}) & \mbox{if $l=c$} 
\end{array}
\right.
\end{eqnarray*}
where $i: S \hookrightarrow T$, $c$ is the codimension of the embedding $j: S \cap S' \hookrightarrow S$ and $N_{S \cap S'/S}$ is the normal bundle of $S \cap S' \subset S$. 
\end{Lemma}
\begin{proof}
We use the sequence of inclusions $S \cap S' \xrightarrow{j} S \xrightarrow{i} T$ where $k = i \circ j$. From corollary 11.2 of \cite{H} we know $\H^{-l}(i^\ast i_\ast \sF) \cong \wedge^l N_{S/T}^\vee \otimes \sF$ where $N_{S/T}$ is the normal bundle of $S$ in $T$. Since $N_{S/T}$ as well as $\sF$ are locally free on $S$ we find $\H^{-l}(k^\ast i_\ast \sF) \cong j^\ast(\wedge^l N_{S/T}^\vee \otimes \sF) = \wedge^l N_{S/T}^\vee|_{S \cap S'} \otimes j^\ast \sF$. In particular, taking $\sF = \O_S$, this means $\H^{-l}(k^\ast \O_S) \cong \wedge^l N_{S/T}^\vee|_{S \cap S'}$. Similarly, $\H^{-l}(k^\ast \O_{S'}) \cong \wedge^l N_{S'/T}^\vee|_{S \cap S'}$. 

Note that (the cohomology of) $\O_{S'}^\vee \otimes i_\ast \sF$ is supported on $S \cap S'$. So to prove the lemma it is enough to show that for any sheaf $\sG$ on $S \cap S'$ the sheaf ${\H}om^l(k_\ast(\sG), \O_{S'}^\vee \otimes i_\ast \sF)$ is equal to 0 if $l < c$ and equal to ${\H}om(k_\ast(\sG), \O_{S \cap S'}(j^\ast \sF \otimes \det N_{S \cap S'/S}))$ if $l = c$. By adjunction, we only need to show that $\H^l(k^!(\O_{S'}^\vee \otimes i_\ast \sF))$ is zero if $l < c$ and $\O_{S \cap S'}(j^\ast \sF \otimes \det N_{S \cap S'/S})$ if $l=c$. 

Now $k^!(\cdot) = k^\ast(\cdot) \otimes \omega_{S \cap S'} \otimes \omega_T^\vee|_{S \cap S'}[-s-c] \cong k^\ast(\cdot) \otimes \det N_{S \cap S'/T}[-s-c]$ where $s$ is the codimension of $S \hookrightarrow T$. Thus 
$$k^!(\O_{S'}^\vee \otimes i_\ast \sF) \cong k^\ast(\O_{S'}^\vee) \otimes k^\ast i_\ast \sF \otimes \det N_{S \cap S'/T}[-s-c].$$
Now $k^\ast(\O_{S'}^\vee) = k^\ast(\O_{S'})^\vee$ while $\H^l(k^\ast \O_{S'}) \cong \wedge^{-l} N_{S'/T}^\vee|_{S \cap S'}$ which are all locally free sheaves. Thus $\H^l(k^\ast \O_{S'}^\vee) \cong \wedge^l N_{S'/T}|_{S \cap S'}$. In particular, $k^\ast \O_{S'}^\vee$ is supported in degrees $0, \dots, s'$ where $s'$ is the codimension of $S' \hookrightarrow T$ and $\H^0(k^\ast \O_{S'}^\vee) \cong \O_{S \cap S'}$. On the other hand, $k^\ast i_\ast \sF$ is supported in degrees $-s, \dots, 0$ with $\H^{-s}(k^\ast i_\ast(\sF)) \cong \det N_{S/T}^\vee|_{S \cap S'} \otimes j^\ast \sF$. From this it follows that $k^\ast(\O_{S'}^\vee \otimes i_\ast \sF)$ is supported in degrees $\ge -s$ so that
$$\H^l(k^\ast \O_{S'}^\vee \otimes k^\ast i_\ast \sF \otimes \det N_{S \cap S'/T}[-s-c]) \cong
\left\{\begin{array}{ll}
0 & \mbox{if $l<c$}  \\
\det N_{S/T}^\vee|_{S \cap S'} \otimes j^\ast(\sF) \otimes \det N_{S \cap S'/T} &  \mbox{if $l=c$.} 
\end{array} \right.
$$
But from the standard short exact sequence $0 \rightarrow N_{S \cap S'/S} \rightarrow N_{S \cap S'/T} \rightarrow N_{S/T}|_{S \cap S'} \rightarrow 0$ we see that $\det N_{S/T}^\vee|_{S \cap S'} \otimes \det N_{S \cap S'/T} \cong \det N_{S \cap S'/S}$. The result follows. 
\end{proof}

\begin{Corollary}\label{cor:abstractnonsense3} Let $S,S' \subset T$ be smooth $\C^\times$-equivariant subvarieties of a smooth quasi-projective variety $T$ equipped with a $\C^\times$ action. Suppose the scheme theoretic intersection $D = S \cap S' \subset S$ is a smooth divisor. Then $\Hom_T(\O_{S'}[-1], \O_S(-D)) \cong \Hom_D(\O_D,\O_D)$ so, in particular, if $T$ is projective there exists a non-zero (equivariant) map $\O_{S'}[-1] \rightarrow \O_S(-D)$ which is unique up to a non-zero multiple. 
\end{Corollary}
\begin{proof}
By lemma \ref{lem:abstractnonsense2} we know 
\begin{eqnarray*}
\Ext^l_T(\O_{S'},\O_S(-D)) = \H^l(\O_{S'}^\vee \otimes \O_S(-D)) \cong
\left\{\begin{array}{ll}
0 & \mbox{if $l<1$}  \\
\O_{S \cap S'}(-D + N_{S \cap S'/S}) &  \mbox{if $l=1$.}
\end{array}
\right.
\end{eqnarray*}
This means 
$$\Hom_T(\O_{S'}[-1],\O_S(-D)) = H^1(T, \O_{S'}^\vee \otimes \O_S(-D)) \cong H^0(T, \O_{S \cap S'}(-D + N_{S \cap S'/S})).$$
On the other hand, since $D = S \cap S' \subset S$ is a divisor, $N_{S \cap S'/S} \cong \O_{S \cap S'}(D)$. Hence 
$$\Hom_T(\O_{S'}[-1],\O_S(-D)) \cong H^0(T, \O_{S \cap S'}) \cong \Hom_D(\O_D,\O_D).$$ 
\end{proof}

\section{Invariance of $\Psi$} \label{se:invariance}
In the previous section we defined the functor $\Psi(T)$ by first choosing a planar diagram representation of $T$. In this section we prove that $\Psi$ does not depend on this choice. To do this we will show that $\Psi$ is invariant under the operations described in lemma \ref{lem:relations}. 

\subsection{Intersection of subvarieties}
We begin with some comments concerning intersection of subvarieties which we will use for proving many of the invariance results.
\begin{Lemma}\label{lem:transverse} For $1 \le i < j \le n-1$ subvarieties $X_n^i \hookrightarrow Y_n$ and $X_n^j \hookrightarrow Y_n$ intersect transversely in a smooth subvariety of codimension two. 
\end{Lemma}
\begin{proof}
Since $X_n^i$ and $X_n^j$ are smooth submanifolds of $ Y_n $, it is enough to work with the topological $ Y_n $.  By Theorem \ref{th:isoman}, $ Y_n \cong \p^{n} $ with $ X_n^i $ taken to $ A_n^i $.  So it suffices to show that $ A_n^i $ and $ A_n^j $ intersect transversely.  Recall that 
\begin{equation*}
A_n^i := \{ (l_1, \dots, l_i, a(l_i), l_{i+2}, \dots , l_{n}) \}.
\end{equation*}

If $ | i - j | > 1$ then $ A_n^i $ and $ A_n^j $ are clearly transverse. 

If $ j = i+1 $ then $ A_n^i \cap A_n^j = \{ (l_1, \dots, l_i, a(l_i), l_i, l_{i+3}, \dots, l_{n}) \} $, since $ a^2 = 1$.  So, either lemma \ref{lem:transversepullback} or a direct examination shows that the intersection is transverse.
\end{proof}

\begin{Remark} The intersection $X_n^i \cap X_n^{i+1} \subset Y_n$ consists of points $(L_1, \dots, L_{n}) \in Y_n$ where $L_{i+1} = z^{-1}(L_{i-1})$ and $L_{i+2} = z^{-1}(L_i)$. This means that setwise, as a subvariety of $X_n^i$, it is the image of the section of $q: X_n^i \rightarrow Y_{n-2}$ given by 
$$s_{i+1}^i: (L_1, \dots, L_{n-2}) \mapsto (L_1, \dots, L_{i-1}, L_i, z^{-1}L_{i-1}, z^{-1}L_i, \dots, z^{-1}L_{n-2}).$$
Since by lemma \ref{lem:transverse} the intersection of $X_n^i$ and $X_n^{i+1}$ is transverse we get that, as a subscheme of $X_n^i$, the scheme theoretic intersection $X_n^i \cap X_n^{i+1}$ is $s_{i+1}^i$. Similarly, the scheme theoretic intersection of $X_n^i$ and $X_n^{i-1}$ inside $X_n^i$ is also given by a section of $q: X_n^i \rightarrow Y_{n-2}$, namely
$$s_{i-1}^i: (L_1, \dots, L_{n-2}) \mapsto (L_1, \dots, L_{i-1}, z^{-1}L_{i-2}, z^{-1}L_{i-1},z^{-1}L_i, \dots, z^{-1}L_{n-2}).$$
Note that if $|i-j|>1$, the intersection of $X_n^i$ and $X_n^j$ as a subscheme of $X_n^i$ is not a section of $q: X_n^i \rightarrow Y_{n-2}$. 
\end{Remark}

\begin{Lemma}\label{lem:transversepullback} Consider the sequence of smooth varieties and morphisms $Y \xleftarrow{f} X \xrightarrow{g} Y'$ and $Y' \xleftarrow{f'} X' \xrightarrow{g'} Y''$ and suppose the maps $(f,g): X \rightarrow Y \times Y'$ and $(f',g'): X' \rightarrow Y' \times Y''$ are embeddings and that $g$ and $f'$ are smooth. Let $\tilde{X} = \pi_{12}^{-1} \circ (f,g)(X)$ and $\tilde{X'} = \pi_{23}^{-1} \circ (f',g')(X')$ where $\pi_{12}$ and $\pi_{23}$ are the projection maps from $Y \times Y' \times Y''$ to $Y \times Y'$ and $Y' \times Y''$. Then $\tilde{X}$ and $\tilde{X'}$ intersect transversely in $Y \times Y' \times Y''$ if and only if $g(X)$ and $f'(X')$ intersect transversely in $Y'$. In particular, if $g$ or $f'$ is surjective then $\tilde{X}$ and $\tilde{X'}$ intersect transversely. 
\end{Lemma}
\begin{proof}
To prove that $\tilde{X}$ and $\tilde{X'}$ intersect transversely we need to show that at every point of $\tilde{X} \cap \tilde{X'}$ their tangent spaces intersect transversely. Suppose $g(X)$ and $f'(X')$ intersect transversely. Fix a point $p = (p_1,p_2,p_3) \in \tilde{X} \cap \tilde{X'} \subset Y \times Y' \times Y''$.  It suffices to show that $ T_p \tilde{X} + T_p \tilde{X'} = Y_{p_1} \oplus Y'_{p_2} \oplus Y''_{p_3} $. 

Fix $ (a,b,c) \in T_{p_1} Y \oplus T_{p_2} Y' \oplus T_{p_3} Y'' $. Since $g(X)$ and $f'(X')$ intersect transversely and $g,f'$ are smooth we know $g_\ast(T_{(p_1, p_2)} X)$ and $f'_\ast(T_{(p_2, p_3)} X')$ intersect transversely and so $g_\ast(T_{(p_1, p_2)} X) + f'_\ast(T_{(p_2, p_3)} X') = T_{p_2} Y'$.  Then there exists $ x \in T_{(p_1,p_2)} X, x' \in T_{(p_2, p_3)} X' $ such that $ g_\ast(x) + f'_\ast(x') = b $. Hence $(f_\ast(x), g_\ast(x), c - g'_\ast(x')) \in T_p \tilde{X} $ and $ (a - f_*(x), f'_*(x'), g'_*(x')) \in T_p \tilde{X'} $ and these two vectors add up to $(a,b,c)$.  Hence $T_p \tilde{X} $ and $T_p \tilde{X'} $ intersect transversely.  The converse is easier.
\end{proof}

\begin{Corollary} \label{cor:transverse}
The following intersections are all transverse:
\begin{enumerate}
\item $\pi_{12}^{-1}(X_n^i) \cap \pi_{23}^{-1}(X_{n+2}^j) $ in $ Y_{n-2} \times Y_n \times Y_{n+2} $,
\item $\pi_{12}^{-1}(X_n^i) \cap \pi_{23}^{-1}(X_n^j) $ in $ Y_n \times Y_{n-2} \times Y_n $, if $ i \ne j$.
\item $\pi_{12}^{-1}(X_{n+2}^i) \cap \pi_{23}^{-1}(X_{n+2}^j) $ in $ Y_n \times Y_{n+2} \times Y_n $, if $ i\ne j$,
\item $\pi_{12}^{-1}(X_n^i \times_{Y_{n-2}} X_n^i) \cap \pi_{23}^{-1} X_n^j $ in $ Y_n \times Y_n \times Y_{n-2}$, if $ i \ne j $.
\item $\pi_{12}^{-1}(X_n^i \times_{Y_{n-2}} X_n^i) \cap \pi_{23}^{-1} X_{n+2}^j $ in $ Y_n \times Y_n \times Y_{n+2} $, if $ i \ne j $.
\item $\pi_{12}^{-1}(X_n^i \times_{Y_{n-2}} X_n^i) \cap \pi_{23}^{-1} (X_n^j \times_{Y_{n-2}} X_n^j) $ in $ Y_n \times Y_n \times Y_n $ if $ i \ne j $.
\end{enumerate}
\end{Corollary}

\begin{proof}
These all follow from applications from Lemma \ref{lem:transversepullback}.

For (i), (ii), (v), we use that $ X_{n+2}^i \rightarrow Y_n $ is surjective and so by the lemma the intersection is transverse.

For (iii), (iv), (vi), by Lemma \ref{lem:transverse}, the image of $ X_n^i $ in $ Y_n $ is transverse to the image of $ X_n^j $ in $ Y_n $ for $ i \ne j$.  Hence by the lemma, the intersections (iii), (iv), (vii) are transverse.
\end{proof}

The main reason for proving that certain intersections are transverse is that transverse intersections are useful in  computing tensor products of structure sheaves, as shown by the following lemma.

\begin{Lemma} \label{th:tensor}
Let $S_1, \dots, S_m \subset T$ and $S_1', \dots, S_n' \subset T$ be irreducible subvarieties such that all intersections $S_i \cap S_j \subset S_i$ and $S_i' \cap S_j' \subset S_i'$ are divisors. If all pairs $S_i$ and $S_j'$ intersect transversely then $\O_S \otimes \O_{S'} \cong \O_{S \cap S'}$ where $S = \cup_i S_i$ and $S' = \cup_i S_i'$.
\end{Lemma}
\begin{proof}
We need to show $\H^l(\O_S \otimes \O_{S'}) = 0$ if $l < 0$. We proceed by induction on the total number of components. The base case is clear so we just need to prove the inductive step. Suppose we have shown $\H^l(\O_S \otimes \O_{S'})=0$ for $l<0$ and $S = S_1 \cup \dots \cup S_{k-1}$. The inductive step is to show $\H^l(\O_{S \cup S_k} \otimes \O_{S'})=0$ if $l<0$. Consider the short exact sequence 
$$0 \rightarrow \O_{S_k}(\mathcal{I}_{S}) \rightarrow \O_{S \cup S_k} \rightarrow \O_{S} \rightarrow 0$$
where $\mathcal{I}_S$ is the ideal sheaf of $S$. Since $S \cap S_k \subset S_k$ is a divisor we get $\O_{S_k}(\mathcal{I}_S) = \O_{S_k}(-D)$. Now $\H^l(\O_S \otimes \O_{S'})=0$ if $l<0$. So it suffices to show $\H^l(\O_{S_k}(-D) \otimes \O_{S'})=0$ if $l<0$. 

Notice $\O_{S_k}(-D) \otimes \O_{S'} = i_\ast(\O_{S_k}(-D) \otimes i^\ast \O_{S'})$ where $i$ is the inclusion $S_k \hookrightarrow T$. Since $i_\ast$ is exact it is enough to show $\H^l(\O_{S_k}(-D) \otimes i^\ast \O_{S'})=0$ for $l < 0$. Since tensoring in $S_k$ by $\O_{S_k}(-D)$ is exact we just need to show $\H^l(i^\ast \O_{S'})=0$ for $l<0$. This is true since by induction $\H^l(i_\ast i^\ast \O_{S'}) = \H^l(\O_{S_k} \otimes \O_{S'}) = 0$ for $l<0$. 
\end{proof}

\subsection{Invariance Under Reidemeister Move (0)}

The first relation we deal with is Reidemeister move (0) from figure \ref{f2} which follows from the following identity.

\begin{Proposition}\label{prop:R0} $\F_n^{i+1} \circ \G_n^i \cong id \cong \F_n^i \circ \G_n^{i+1}$.
\end{Proposition}

\begin{proof}
To prove the first isomorphism recall that the functors $ \F_n^{i+1}, \G_n^i $ are the Fourier-Mukai transforms with respect to the kernels $ \sF_n^{i+1} = \O_{X_n^{i+1}} \otimes \E_{i+2}^\vee \{i+1\}, \sG_n^i = \O_{X_n^i} \otimes \E'_i \{-i+1\}$.

By Corollary \ref{cor:transverse}, the intersection $ W := \pi_{12}^{-1}(X_n^i) \cap \pi_{23}^{-1}(X_n^{i+1}) $ inside $ Y_{n-2} \times Y_n \times Y_{n-2} $ is transverse.  Moreover 
\begin{equation*}
\begin{aligned}
W = \{ (L_\cdot, L'_\cdot, L''_\cdot) : &L_j = L'_j \text{ for } j \le i,\ L_j = zL'_{j+2} \text{ for } j \ge i, \\ &L''_j = L'_j \text{ for } j \le i+1,\ L''_j = zL'_{j+2} \text{ for } j \ge i+1 \} .
\end{aligned}
\end{equation*}

Examining this variety closely, we see that both projections, $ \pi_1, \pi_3 : W \rightarrow Y_{n-2} $, are isomorphisms and in fact $ \pi_{13} $ maps $ W $ isomorphically to the diagonal in $Y_{n-2} \times Y_{n-2}$.  Moreover, as $ zL'_{i+2} = L'_i\{2\} $ and $ zL'_{i+1} = L'_{i-1}\{2\}$, we see that the operator $ z $ induces an isomorphism of line bundles $ \E'_{i+2} \cong \E'_i\{2\}$ on $W$.

Hence the composition $ \F_n^{i+1} \circ \G_n^i $ is a Fourier-Mukai transform with respect to the kernel
\begin{equation*}
\begin{aligned}
{\pi_{13}}_*( \pi_{12}^*(\O_{X_n^i} &\otimes \E'_i\{-i+1\}) \otimes \pi_{23}^*(\O_{X_n^{i+1}} \otimes \E_{i+2}^\vee \{i+1\})) \\
&= {\pi_{13}}_*( \O_{\pi_{12}^{-1}(X_n^i)} \otimes \O_{\pi_{23}^{-1}(X_n^{i+1})} \otimes \E'_i \otimes {\E'_{i+2}}^\vee \{2\}) \\
&\cong {\pi_{13}}_*( \O_W \otimes \E'_i \otimes {\E'_{i+2}}^\vee \{2\}) \\
&\cong {\pi_{13}}_*(\O_W) = \O_{\Delta}
\end{aligned}
\end{equation*}

Thus, we conclude that $ \F_n^{i+1} \circ \G_n^i \cong id $. 

The second isomorphism is actually easier to prove. $\F_n^i \circ \G_n^{i+1}$ is the Fourier-Mukai transform with respect to the kernel
\begin{equation*}
\begin{aligned}
{\pi_{13}}_*( \pi_{12}^*(\O_{X_n^{i+1}} \otimes \E'_{i+1}\{-i\}) \otimes \pi_{23}^*(\O_{X_n^{i}} \otimes {\E'_{i+1}}^{\vee} \{i\}))
&= {\pi_{13}}_*( \O_{\pi_{12}^{-1}(X_n^{i+1})} \otimes \O_{\pi_{23}^{-1}(X_n^{i})} ) \\
&\cong \O_{\pi_{13}(\pi_{12}^{-1}(X_n^{i+1}) \cap \pi_{23}^{-1}(X_n^{i}))} = \O_{\Delta} \\
\end{aligned}
\end{equation*}

\end{proof}

\subsection{Invariance Under Reidemeister Move (I)}
We begin with a lemma on non-transverse intersections.
\begin{Lemma} \label{th:nontransint}
Let $ X $ be a smooth projective variety with a $\C^\times$ action and let $ D $ be a smooth $\C^\times$-equivariant divisor in $ X $.  Let $ Y, Z $ be smooth equivariant subvarieties of $ D $ and assume that they meet transversely in $ D$.  Then we have the following distinguished triangle in $ D(X) $,
\begin{equation*}
\O_{Y \cap Z} \otimes \O_X(-D) [1] \rightarrow \O_Y \otimes \O_Z \rightarrow \O_{Y \cap Z} 
\end{equation*}
\end{Lemma}
\begin{proof}
Denote by $i$ the equivariant inclusion $Y \rightarrow X$ and by $i_2$ and $i_1$ the sequence of equivariant inclusions $Y \rightarrow D \rightarrow X$ so that $i = i_2 \circ i_1$. To compute $i_\ast i^\ast \O_Z = \O_Y \otimes \O_Z$ we first calculate $i_1^\ast \O_Z$. Since $Z \subset D$ we find that $\H^0(i_1^\ast \O_Z) \cong \O_Z$ and $\H^{-1}(i_1^\ast \O_Z) \cong \O_Z(-D)$ while the rest of the cohomology groups vanish (see for instance proposition 11.7 of \cite{H}). This means there is a distinguished triangle
$$\O_Z(-D)[1] \rightarrow i_1^\ast \O_Z \rightarrow \O_Z.$$
Since $Z$ and $Y$ intersect transversely inside $D$ we get $i_{2\ast}i_2^\ast \O_Z \cong \O_{Y \cap Z}$. Thus we get a distinguished triangle
$$\O_{Y \cap Z} \otimes i_1^\ast \O_X(-D)[1] \rightarrow i_{2\ast} i_2^\ast i_1^\ast \O_Z \rightarrow \O_{Y \cap Z}.$$
Applying $i_{1\ast}$ and remembering that $i_{1\ast}i_{2\ast}i_2^\ast i_1^\ast \O_Z = i_\ast i^\ast \O_Z = \O_Y \otimes \O_Z$ we get the distinguished triangle
$$\O_{Y \cap Z} \otimes \O_X(-D)[1] \rightarrow \O_Y \otimes \O_Z \rightarrow \O_{Y \cap Z}.$$
\end{proof}

To prove invariance under Reidemeister move (I) we need to understand the functor $\F_n^i \circ \G_n^i(\cdot)$. We will eventually see that $\F_n^i \circ \G_n^i(\cdot) \cong (\cdot) \otimes_\C V$ where $V = \C[-1]\{1\} \oplus \C[1]\{-1\}$, but first we need the following result. 

\begin{Proposition} \label{th:circ}
There is a distinguished triangle
\begin{equation*}
\O_\Delta[1]\{-1\} \rightarrow \sF_n^i * \sG_n^i \rightarrow \O_\Delta[-1]\{1\} 
\end{equation*}
where the first map is the adjoint map $\O_\Delta \rightarrow \sF_n^i \ast {\sF_n^i}_L$ and the second map is the adjoint map $\sF_n^i \ast {\sF_n^i}_R \rightarrow \O_\Delta$ via the isomorphisms $ \sG_n^i \cong {\sF_n^i}_R[-1]\{1\} \cong {\sF_n^i}_L[1]\{-1\}$.
\end{Proposition}
\begin{proof}
We have 
\begin{equation*}
\begin{aligned}
\sF_n^i * \sG_n^i &= {\pi_{13}}_* ( \pi_{12}^*(\O_{X_n^i} \otimes \E'_i \{-i+1\}) \otimes \pi_{23}^*(\O_{X_n^i} \otimes \E_{i+1}^\vee \{i\})) \\
&= {\pi_{13}}_*( \O_{\pi_{12}^{-1}(X_n^i)} \otimes \O_{\pi_{23}^{-1}(X_n^i)} \otimes \E'_i \otimes {\E'_{i+1}}^\vee) \{1\}
\end{aligned}
\end{equation*}

Now the intersection $ \pi_{12}^{-1} (X_n^i) = X_n^i \times Y_{n-2} \cap \pi_{23}^{-1}(X_n^i) = Y_{n-2} \times X_n^i $ inside $ Y_{n-2} \times Y_n \times Y_{n-2} $ is not transverse.  In particular, both $ X_n^i \times Y_{n-2} $ and $ Y_{n-2} \times X_n^i $ are contained inside the divisor $ D := Y_{n-2} \times i(X_n^i) \times Y_{n-2} $.  (To avoid confusion, in this section we are writing $ i(X_n^i) $ when we view $ X_n^i $ as a subvariety of $ Y_n $ and $ X_n^i $ when we view it as a subvariety of $ Y_{n-2} \times Y_n $ or $ Y_n \times Y_{n-2}$.)

So Lemma \ref{th:nontransint} applies in our situation.  Let $ W = X_n^i \times Y_{n-2} \cap Y_{n-2} \times X_n^i$.  Also, we have by Lemma \ref{lem:facts} that $ \O_{Y_{n-2} \times Y_n \times Y_{n-2}}(-D) \cong {\E'_i}^\vee \otimes \E'_{i+1} \{-2\}$. We conclude that there is a distinguished triangle
\begin{equation*}
\O_W \otimes {\E'_i}^\vee \otimes \E'_{i+1} [1] \{-2\} \rightarrow \O_{\pi_{12}^{-1}(X_n^i)} \otimes \O_{\pi_{23}^{-1}(X_n^i)} \rightarrow \O_W .
\end{equation*}

Tensoring by $\E'_i \otimes \E_{i+1}^{'\vee}\{1\}$ we get a distinguished triangle
\begin{equation*}
{\pi_{13}}_* \O_W [1]\{-1\} \rightarrow \sF_n^i * \sG_n^i \rightarrow {\pi_{13}}_*(\O_W \otimes \E'_i \otimes {\E'_{i+1}}^\vee) \{1\}
\end{equation*}

Now $ \pi_{13}(W) = \Delta \subset Y_{n-2} \times Y_{n-2} $ where $W \rightarrow \Delta$ is the projection $q: X_n^i \rightarrow Y_{n-2}$. In particular we get 
\begin{equation*} 
{\pi_{13}}_*(\O_W )[1] \cong \O_\Delta [1]
\end{equation*}
Meanwhile, by corollary \ref{lem:facts} we know that $\E'_i \otimes {\E'_{i+1}}^\vee$ restricted to $W = X_n^i$ is precisely the relative dualizing sheaf of $q$. Consequently, 
\begin{equation*}
{\pi_{13}}_* (\O_W \otimes \E'_i \otimes {\E'_{i+1}}^\vee) \cong \O_\Delta[-1].
\end{equation*}

Thus we conclude that there is a distinguished triangle
\begin{equation*}
\O_\Delta[1]\{-1\} \rightarrow \sF_n^i * \sG_n^i \rightarrow \O_\Delta [-1]\{1\}
\end{equation*}
as desired. To show the maps are as claimed we prove that the Hom spaces $\O_\Delta[1]\{-1\} \rightarrow \sF_n^i \ast \sG_n^i$ and $\sF_n^i \ast \sG_n^i \rightarrow \O_\Delta[-1]\{1\}$ are one dimensional. Since the maps in our distinguished triangle are non-zero and adjoint maps are non-zero this implies they must agree up to a non-zero multiple. 

Applying the functor $\Hom(-,\O_\Delta[-1]\{1\})$ to our distinguished triangle we obtain
$$\Hom(\O_\Delta[-1]\{1\},\O_\Delta[-1]\{1\}) \rightarrow \Hom(\sF_n^i \ast \sG_n^i,\O_\Delta[-1]\{1\}) \rightarrow \Hom(\O_\Delta[1]\{-1\},\O_\Delta[-1]\{1\}).$$
But $\Hom(\O_\Delta[-1]\{1\},\O_\Delta[-1]\{1\}) \cong \C$ since $\Delta = Y_n$ is projective while 
\begin{equation*}
\Hom(\O_\Delta[1]\{-1\},\O_\Delta[-1]\{1\}) = \Ext^{-2}(\O_\Delta,\O_\Delta\{2\})=0.
\end{equation*}
Thus $\Hom(\sF_n^i \ast \sG_n^i,\O_\Delta[-1]\{1\})$ is at most one dimensional which means that up to a non-zero multiple there is at most one non-zero map $\sF_n^i \ast \sG_n^i \rightarrow \O_\Delta[-1]\{1\}$. Similarly, applying the functor $\Hom(\O_\Delta[1]\{-1\},-)$ to our distinguished triangle one finds that up to a non-zero multiple there is at most one non-zero map $\O_\Delta[1]\{-1\} \rightarrow \sF_n^i \ast \sG_n^i$. 
\end{proof}

Now we will show invariance under Reidemeister move (I).
\begin{Theorem}
$\F_n^i \circ \T_n^{i \pm 1}(2) \circ \G_n^i \cong id \cong \F_n^i \circ \T_n^{i \pm 1}(1) \circ \G_n^i $
\end{Theorem}

\begin{proof}
We will just prove the left statement, since the right statement follows from taking left adjoints.

To prove the left equality it suffices to prove the stronger statement
\begin{equation*}
\sF_n^i * \sT_n^{i \pm 1}(2) * \sG_n^i \cong \O_\Delta.
\end{equation*}

Let $\R = \sF_n^i * \sT_n^{i \pm 1}(2) * \sG_n^i$.

By theorem \ref{th:kerneltwist}, $\sT_n^{i \pm 1}(2) \cong \Cone(\sG_n^{i \pm 1} * \sF_n^{i \pm 1}[-1]\{1\} \xrightarrow{\beta_{\sG_n^{i \pm 1}}} \O_\Delta)[-1]\{1\}$.  Hence the left hand side is
\begin{equation*}
\R \cong \Cone(\sF_n^i * \sG_n^{i \pm 1} * \sF_n^{i \pm 1} * \sG_n^i [-1] \xrightarrow{\sF_n^i * \beta_{\sF_n^i * \sG_n^{i \pm 1}}* \sG_n^i} \sF_n^i * \sG_n^i) [-1] \{1\}.
\end{equation*}

Now, we know from the proof of \ref{prop:R0} that $ \sF_n^i * \sG_n^{i \pm 1} \cong \O_\Delta \cong \sF_n^{i \pm 1} * \sG_n^i$.  Hence we see that there is a distinguished triangle
\begin{equation*}
\O_\Delta[-1]\{1\} \xrightarrow{\psi} \sF_n^i * \sG_n^i \rightarrow \R [1]\{-1\}.
\end{equation*}
where $ \psi $ denotes the result of transferring $ \sF_n^i * \beta_{\sG_n^{i \pm 1}} * \sG_n^i $ via the isomorphisms $ \sF_n^i * \sG_n^{i \pm 1} * \sF_n^{i \pm 1} * \sG_n^i [-1]\{1\} \cong \O_\Delta * \O_\Delta[-1]\{1\} \cong \O_\Delta[-1]\{1\} $.

Now, by Proposition \ref{th:circ}, there is a distinguished triangle
\begin{equation*}
\O_\Delta[1]\{-1\} \rightarrow \sF_n^i * \sG_n^i \xrightarrow{\phi} \O_\Delta[-1]\{1\}
\end{equation*}

We claim that $ \phi \circ \psi $ is a non-zero multiple of the identity.  Once we know this we will apply Lemma \ref{th:twodist} to conclude that $ \O_\Delta \cong \R $.  Also from this lemma, we see that 
\begin{equation} \label{eq:split}
\sF_n^i * \sG_n^i \cong \O_\Delta[1]\{-1\} \oplus \O_\Delta[-1]\{1\}.
\end{equation}

Now we show that $ \phi \circ \psi $ is a non-zero multiple of the identity.  We have $ {\sG_n^{i \pm 1}}_R \cong \sF_n^{i \pm 1}[-1]\{1\} $ and $ {\sF_n^i}_R \cong \sG_n^i[1]\{-1\} $.  Hence by Lemma \ref{th:convbeta} we deduce the commutativity of the diagram
\begin{equation*}
\begin{CD}
\sF_n^i * \sG_n^{i \pm 1} * \sF_n^{i \pm 1} * \sG_n^i [-1]\{1\} @>{\sF_n^i * \beta_{\sG_n^{i \pm 1}} * \sG_n^i}>> \sF_n^i * \sG_n^i \\
@VVV @V{\beta_{\sF_n^i}[-1]\{1\}}VV \\
(\sF_n^i * \sG_n^{i+1}) * (\sF_n^i * \sG_n^{i+1})_R[-1]\{1\} @>{\beta_{\sF_n^i * \sG_n^{i \pm 1}}[-1]\{1\}}>> \O_\Delta[-1]\{1\}
\end{CD}
\end{equation*}

Now, by Proposition \ref{th:circ}, $ c \phi = \beta_{\sF_n^i}[-1]\{1\} $ (where $ c \in \C^\times$) and by definition $ \psi $ is the composition of $ \beta_{\sG_n^{i \pm 1}}$ with the inverses of the isomorphisms $\sF_n^i * \sG_n^{i \pm 1} * \sF_n^{i \pm 1} * \sG_n^i [-1]\{1\} \cong \O_\Delta * \O_\Delta[-1]\{1\} \cong \O_\Delta[-1]\{1\}$.  So up to the identifications that we have made, $ c \phi $ is the right arrow and $ \psi $ is the top arrow, and so $ c \phi \circ \psi $ is the result of transferring $ \beta_{\sF_n^i * \sG_n^{i \pm 1}}[-1]\{1\} $ to a map $ \O_\Delta[-1]\{1\} \rightarrow \O_\Delta[-1]\{1\} $.

By the naturality of $ \beta $ the diagram
\begin{equation*}
\begin{CD}
(\sF_n^i * \sG_n^{i \pm 1}) * (\sF_n^i * \sG_n^{i \pm 1})_R @>{\beta_{\sF_n^i * \sG_n^{i \pm 1}}}>> \O_\Delta \\
@VVV @| \\
\O_\Delta * \O_\Delta @>{\beta_{\O_\Delta}}>> \O_\Delta \\
\end{CD}
\end{equation*}
commutes. Now $ \beta_{\O_\Delta} $ becomes the identity when we identify $ \O_\Delta * \O_\Delta $ with $ \O_\Delta $.  Thus we conclude that $ c \phi \circ \psi = id $. This proves the left equality. 

\end{proof}

From (\ref{eq:split}) we immediately deduce the following corollary of our proof.
\begin{Corollary}\label{cor:circle} $\F_n^i \circ \G_n^i (\cdot) \cong (\cdot)[-1]\{1\} \oplus (\cdot)[1]\{-1\}$ so if $T$ is a tangle then $\Psi(T \cup O)(\cdot) \cong \Psi(T)(\cdot) \otimes_\C V$ where $V = \C[-1]\{1\} \oplus \C[1]\{-1\} $. 
\end{Corollary}

\begin{Remark} \label{re:Dolb}
In general, let $X$ be a smooth projective variety and consider the embedding as the zero section $X \rightarrow T^* X \{l\}$ where $\C^\times$ acts trivially on the base of $\pi: T^* X \{l\} \rightarrow X$ and by $t \cdot v \mapsto t^{-l} v$ on the fibres. Consider the natural Koszul resolution of $\O_X$ in $ D(T^* X) $:
$$0 \leftarrow \O_X \leftarrow \pi^\ast \O_X \leftarrow \pi^\ast T_X \{-l\} \leftarrow \pi^\ast \wedge^2 T_X \{-2l\} \leftarrow \dots.$$
If we take its dual and tensor with $\O_X$ we get 
$$\O_X \rightarrow \Omega^1_X \{l\} \rightarrow \Omega^2_X \{2l\} \rightarrow \dots$$
where all the maps are zero. Thus we find that
$$\Ext^k_{D(T^* X \{l\})}(\O_X,\O_X) = \oplus_{i+j=k} H^j(X, \Omega_X^i)\{il\} = \oplus_{i+j=k} H^{i,j}(X)\{il\}$$
which is a graded version of Dolbeault cohomology. 

In our case $Y_2$ is nothing but the Hirzebruch surface ${\mathbb{P}}(\O_\p \oplus \O_\p(-2))$ with $\p \cong  X_2^1 \subset Y_2$ embedded as the $-2$ section. The $\C^\times$ action fixes $\p \subset Y_2$ and acts on the fibres by $t \cdot [x,y] = [t^{-2}x, y]$. Thus 
$$\Ext_{D(Y_2)}(\O_\p,\O_\p) \cong \Ext_{D(T^*\p \{2\})}(\O_\p, \O_\p) \cong \oplus_k \oplus_{i+j=k} H^{i,j}(\p)[-k]\{2i\}$$ 
which is isomorphic to $\C \oplus \C[-2]\{2\}$. Since the bigraded vector space $ V $ appearing in the corollary is $ \F_2^1 \circ \G_2^1(\C) = \Ext_{D(Y_2)}(\O_\p,\O_\p) [1]\{-1\} $ we can identify $V$ with the $[1]\{-1\}$ shift of the graded Dolbeault cohomology of $\p$. In the non-equivariant situation this reduces to the $[1]$ shift of the cohomology of $\p$. 
\end{Remark}

\subsection{Invariance Under Reidemeister Moves (II) and (III)}

\begin{Proposition}\label{prop:R2} 
$\T_n^i(2) \circ \T_n^i(1) \cong id \cong \T_n^i(1) \circ \T_n^i(2) $
\end{Proposition}
\begin{proof}
Since $\sF_n^i \cong (\sG_n^i)_R [1]\{-1\}$, Proposition (\ref{th:circ}) yields the distinguished triangle
\begin{equation*}
\O_\Delta \rightarrow {\sG_n^i}_R * \sG_n^i \rightarrow \O_\Delta[-2]\{2\} 
\end{equation*}
By Lemma \ref{th:Gadj}, $ {\sG_n^i}_R \cong {\sG_n^i}_L[-2]\{2\}$. Also, since $\sG_n^i = \O_{X_n^i} \otimes \E'_i \{-i+1\}$ we have $\Hom(\sG_n^i, \sG_n^i) = \Hom(\O_{X_n^i}, \O_{X_n^i})$ where we view $X_n^i$ as a subvariety of $Y_{n-2} \times Y_n$. Thus, for any $k \in {\mathbb{Z}}$,
\begin{eqnarray*}
\Hom(\sG_n^i, \sG_n^i[j]\{k\}) \cong \left\{
\begin{array}{ll}
\C & \mbox{if } j=0=k \\
0  & \mbox{if } j<0 
\end{array}
\right.
\end{eqnarray*} 
Since $\dim(Y_{n-2}) - \dim(Y_n)=-2$ the conditions of Theorem \ref{th:twistequiv} are satisfied. So the twist in $ \sG_n^i $ is an equivalence. By Theorem \ref{th:kerneltwist}, the functor $ \T_n^i(2) $ is isomorphic to the twist in $ \sG_n^i$ shifted by $[-1]\{1\}$. Hence the functor $ \T_n^i(2) $ is also an equivalence.

Since $\T_n^i(2)$ is an equivalence its inverse is its left adjoint. By Lemma \ref{th:twistadjoint} its left adjoint is $ \T_n^i(1) $. 
\end{proof}

So $ \G_n^i $ is a spherical functor. Now, we verify the hypothesis which ensures the braid relation for the twists.
\begin{Lemma} \label{th:forpitch}
\begin{equation*}
\sT_{\sG_n^{i+1}} \ast \sG_n^i \cong {\sT_{\sG_n^i}}_R \ast \sG_n^{i+1} \{1\} \text{ and } 
\sT_{\sG_n^{i}} \ast \sG_n^{i+1} \cong {\sT_{\sG_n^{i+1}}}_R \ast \sG_n^{i} \{1\}
\end{equation*}
\end{Lemma}

\begin{proof}
We first prove the left equality. By definition and Theorem \ref{th:kerneltwist} we have that 
\begin{equation*} 
\sG_n^i = \O_{X_n^i} \otimes \E'_i \{-i+1\} \text{ and }
\sT_{\sG_n^i} = \O_{Z_n^i} \otimes \E_{i+1}^\vee \otimes \E'_i \{2\} \text{ and }
{\sT_{{\sG_n^i}}}_R = {\sT_{{\sG_n^i}}}_L = \O_{Z_n^i}
\end{equation*}

Now by Corollary \ref{cor:transverse}, $ \pi_{12}^{-1}(X_n^i) $ is transverse to each of the two components of $ \pi_{23}^{-1}(Z_n^{i+1}) $.  Let $ W = \pi_{12}^{-1}(X_n^i) \cap \pi_{23}^{-1}(Z_n^{i+1}) $.  Note that $ W $ is the locus
\begin{equation*}
\begin{aligned}
W = \{ (L_\cdot, L'_\cdot, L''_\cdot) \in Y_{n-2} \times Y_n \times Y_n :\ &L_j = L'_j \text{ for } j \le i-1,\ L_j = zL'_{j+2}\text{ for } j \ge i-1,\\
& L'_j = L''_j \text{ for } j \ne i+1 \}.
\end{aligned}
\end{equation*}
In particular, considering the definitions of the line bundles, we see that on $ W $ we have isomorphisms of line bundles 
\begin{equation} \label{eq:isoline1}
 \E_i \{2\} \cong \E'_{i+2} \text{ and } \quad \E'_i \cong \E''_i,
\end{equation}
two facts which we will use later.

Note that the fibres of $ \pi_{13} $ restricted to $ W $ are all points and that  
\begin{equation*}
\pi_{13}(W) = \{ (L_\cdot, L'_\cdot) \in Y_{n-2} \times Y_n : L_j = L'_j \text { for all } j \le i-1, L_j = z L'_{j+2} \text{ for all } j \ge i \}.
\end{equation*}

By Lemma \ref{th:tensor}, we see that 
\begin{equation*}
\pi_{12}^\ast(\O_{X_n^i}) \otimes \pi_{23}^\ast(\O_{Z_n^{i+1}}) = \O_{\pi_{12}^{-1}(X_n^i)} \otimes \O_{\pi_{23}^{-1}(Z_n^{i+1})} \cong \O_W
\end{equation*}

Hence, we have that
\begin{equation} \label{eq:lhs}
\begin{aligned}
\sT_{\sG_n^{i+1}} \ast \sG_n^i &= {\pi_{13}}_\ast(\pi_{12}^*(\O_{X_n^i} \otimes \E'_i \{-i+1\}) \otimes \pi_{23}^*(\O_{Z_n^{i+1}} \otimes \E_{i+2}^\vee \otimes \E'_{i+1} \{2\})) \\
&= {\pi_{13}}_\ast( \pi_{12}^*(\O_{X_n^i}) \otimes \pi_{23}^*(\O_{Z_n^{i+1}}) \otimes \E'_i \otimes {\E'_{i+2}}^\vee \otimes \E''_{i+1}) \{-i+3\} \\
& \cong {\pi_{13}}_\ast( \O_W \otimes \E'_i \otimes {\E'_{i+2}}^\vee \otimes \E''_{i+1}) \{-i+3\} \\
& \cong {\pi_{13}}_\ast( \O_W \otimes {\E_i}^\vee \otimes \E''_i \otimes \E''_{i+1}) \{-i+1\} \\ 
&= \O_{\pi_{13}(W)} \otimes {\E_i}^\vee \otimes \E'_i \otimes \E'_{i+1} \{-i+1\}
\end{aligned}
\end{equation}
where in the second last step we use the isomorphisms of line bundles (\ref{eq:isoline1}).

Now, we compute the right hand side $ {\sT_{\sG_n^i}}_R \ast \sG_n^{i+1} \{1\} $.
Again the intersection $ W' := \pi_{12}^{-1} X_n^{i+1} \cap \pi_{23}^{-1} Z_n^i $ is transverse. We have 
\begin{equation*}
W' = \{ (L_\cdot, L'_\cdot, L''_\cdot) \in Y_{n-2} \times Y_n \times Y_n : L_j = L_j' \text{ for } j \le i,\ L_j = z L_{j+2}'  \text{ for } j \ge i,\ L'_j = L''_j \text{ for } j \ne i \}.
\end{equation*}

The fibres of $ \pi_{13} $ restricted to $ W' $ are again points and we see that
\begin{equation*}
\pi_{13}(W') = \pi_{13}(W).
\end{equation*}

Hence we have
\begin{equation} \label{eq:rhs}
\begin{aligned}
{\sT_{\sG_n^i}}_R \ast \sG_n^{i+1} \{1\} &= {\pi_{13}}_\ast(\pi_{12}^*(\O_{X_n^{i+1}} \otimes \E'_{i+1} \{-i\}) \otimes \pi_{23}^*(\O_{Z_n^i})) \{1\} \\
&= {\pi_{13}}_\ast( \pi_{12}^*(\O_{X_n^{i+1}}) \otimes \pi_{23}^*(\O_{Z_n^i}) \otimes \E'_{i+1} ) \{-i+1\} \\
& \cong {\pi_{13}}_\ast( \O_{W'} \otimes  \E'_{i+1} ) \{-i+1\} \\
\end{aligned}
\end{equation}

Suppose that we have a diamond of vector bundles $ V_1, V_2, V_3, V_4 $ on any variety,
\begin{equation*}
\xymatrix{
&V_1 \\
V_2 \ar@{-}[ur] & & V_3 \ar@{-}[ul] \\
&V_4 \ar@{-}[ul] \ar@{-}[ur] \\
}
\end{equation*}
where each line indicates containment of dimension 1.  In this case, we can factor $ \det V_1 / V_4 $ in two different ways.  Namely we have isomorphisms of line bundles 
\begin{equation*}
 V_1/V_2 \otimes V_2/V_4 \cong \det V_1/ V_4 \cong  V_1/V_3 \otimes V_3/V_4.
\end{equation*}  
On $ W' $ we have a diamond 
\begin{equation*}
\xymatrix{
&L'_{i+1} = L''_{i+1} \\
L_i = L'_i \ar@{-}[ur] & & L''_i \ar@{-}[ul] \\
&L_{i-1} = L'_{i-1} = L''_{i-1} \ar@{-}[ul] \ar@{-}[ur] \\
}
\end{equation*}
Hence on $ W' $, we have
\begin{equation*}
\E'_{i+1} \cong {\E_i}^\vee \otimes \E''_{i+1} \otimes \E''_i.
\end{equation*}
Combining with (\ref{eq:rhs}) shows that
\begin{equation*}
{\sT_{\sG_n^i}}_R \ast \sG_n^{i+1} \{1\} \cong {\pi_{13}}_*( \O_{W'} \otimes \E_i^\vee \otimes \E''_i \otimes \E''_{i+1}) \{-i+1\} = \O_{\pi_{13}(W')} \otimes \E_i^\vee \otimes \E'_i \otimes \E'_{i+1} \{-i+1\}
\end{equation*}

Comparing this with (\ref{eq:lhs}) completes the proof. 

The proof of the second equality is a little easier so we just sketch it. We have
\begin{equation*} 
\begin{aligned}
\sT_{\sG_n^{i}} \ast \sG_n^{i+1} &= {\pi_{13}}_\ast(\pi_{12}^*(\O_{X_n^{i+1}} \otimes \E'_{i+1} \{-i\}) \otimes \pi_{23}^*(\O_{Z_n^{i}} \otimes \E_{i+1}^\vee \otimes \E'_{i} \{2\})) \\
&= {\pi_{13}}_\ast( \pi_{12}^*(\O_{X_n^{i+1}}) \otimes \pi_{23}^*(\O_{Z_n^{i}}) \otimes \E'_{i+1} \otimes {\E'_{i+1}}^\vee \otimes \E''_{i}) \{-i+2\} \\
& \cong {\pi_{13}}_\ast( \O_{\pi_{12}^{-1}(X_n^{i+1}) \cap \pi_{23}^{-1}(Z_n^i)} \otimes \E''_{i}) \{-i+2\} \\
&\cong \O_{\pi_{13}(\pi_{12}^{-1}(X_n^{i+1}) \cap \pi_{23}^{-1}(Z_n^i))} \otimes \E'_{i} \{-i+2\}
\end{aligned}
\end{equation*}
whereas 
\begin{equation*}
\begin{aligned}
{\sT_{\sG_n^{i+1}}}_R \ast \sG_n^{i} \{1\} &= {\pi_{13}}_\ast(\pi_{12}^*(\O_{X_n^{i}} \otimes \E'_{i} \{-i+1\}) \otimes \pi_{23}^*(\O_{Z_n^{i+1}})) \{1\} \\
&= {\pi_{13}}_\ast( \pi_{12}^*(\O_{X_n^{i}}) \otimes \pi_{23}^*(\O_{Z_n^{i+1}}) \otimes \E'_{i} ) \{-i+2\} \\
&\cong {\pi_{13}}_\ast( \O_{\pi_{12}^{-1}(X_n^i) \cap \pi_{23}^{-1}(Z_n^{i+1})} \otimes  \E''_{i} ) \{-i+2\} \\
&= \O_{\pi_{13}(\pi_{12}^{-1}(X_n^i) \cap \pi_{23}^{-1}(Z_n^{i+1}))} \otimes \E'_i \{-i+2\}
\end{aligned}
\end{equation*}
where the third equality is a consequence of the fact that $\E'_i \cong \E''_i$ on $\pi_{12}^{-1}(X_n^i) \cap \pi_{23}^{-1}(Z_n^{i+1})$. We have already checked that  
$$\pi_{13}(\pi_{12}^{-1}(X_n^{i+1}) \cap \pi_{23}^{-1}(Z_n^i)) = \pi_{13}(\pi_{12}^{-1}(X_n^i) \cap \pi_{23}^{-1}(Z_n^{i+1}))$$
so $\sT_{\sG_n^{i}} \ast \sG_n^{i+1} \cong {\sT_{\sG_n^{i+1}}}_R \ast \sG_n^{i} \{1\}$. 
\end{proof}

\begin{Proposition}\label{prop:R3} 
The map $ \Psi $ is invariant under Reidemeister III moves, ie.
$$\T_n^i(l) \circ \T_n^{i+1}(l') \circ \T_n^i(l'') \cong \T_n^{i+1}(l'') \circ \T_n^i(l') \circ \T_n^{i+1}(l)$$
\end{Proposition}
\begin{proof}
We show invariance under the R-move (III) shown in figure \ref{f2} (ie $ l = l' = l'' = 2$). Reversing the orientations of some of the strands only changes the twist functors by a shift so it is then easy to check that invariance under any R-move of type (III) follows from invariance under the R-move (III) in figure \ref{f2}. 

By Theorem \ref{th:kerneltwist}, we know that $ \T_n^i(2) $ is the twist in the spherical functor $ \G_n^i $ (shifted by $[-1]\{1\}$).  Since for any $k \in {\mathbb{Z}}$ and FM kernel $\P$ we have $\T_{\P} \cong \T_{\P\{k\}}$, Lemma \ref{th:forpitch} implies
$$ \T_{{\sT_{\sG_n^{i+1}} \ast \sG_n^i}} \cong \T_{{\sT_{\sG_n^i}}_R \ast \sG_n^{i+1}} \text{ and } 
\T_{\sT_{\sG_n^{i}} \ast \sG_n^{i+1}} \cong \T_{{\sT_{\sG_n^{i+1}}}_R \ast \sG_n^{i}}. $$
Thus by Theorem \ref{th:braidrel}, we see that the twists in the spherical functors $\G_n^i$ and $\G_n^{i+1}$ satisfy the braid relations. Hence $ \T_n^i(2), \T_n^{i+1}(2) $ satisfy the braid relations.
\end{proof}

\begin{Remark}
Propositions \ref{prop:R2}, \ref{prop:R3}, and \ref{prop:tangleiso} show that the functors $ \T_n^i(2), \T_n^i(1) $ give a weak action of the braid group $ B_{n} $ on the derived category of coherent sheaves on $ Y_n $.  It can be shown that similar functors (given by restricting the kernels from $ Y_{2n} \times Y_{2n} $ to $ U_n \times U_n $) give a weak action of $ B_{2n} $ on the derived category of coherent sheaves on the open subset $ U_n $.  Recall that in Proposition \ref{prop:Unslice} we showed that there was an isomorphism $ U_n $ with the resolution of a slice to a nilpotent orbit in $ \mathfrak{gl}_{2n} $.  In \cite{B}, Bezrukavnikov describes an action of the braid group $ B_{\mathfrak{g}} $ on the resolution of a slice to any nilpotent orbit in any semisimple Lie algebra $ \mathfrak{g} $.  This action was one of the original motivations for our construction and it can be shown that the two constructions coincide.

While our paper was in preparation, the paper \cite{KT} by Khovanov-Thomas appeared.  In that paper (among other things), they construct an action of $ B_n$ on the derived category of coherent sheaves on the cotangent bundle to the flag variety.  The generators of the braid group act as Fourier-Mukai transforms with respect to kernels similar to our kernels $ \sT_n^i(l)$.  They also use the fact that these functors are spherical twists.  Their proof of the braid relation is more complicated than ours because the geometry of their situation is more complicated.
\end{Remark}

\subsection{Invariance Under Vertical Isotopies From Figure \ref{f3}}\label{se:invdistant}

\begin{Proposition}\label{prop:capscups} The functor $ \Psi $ is invariant with respect to isotopies exchanging the order with respect to height of caps and cups.  We have
\begin{itemize}
\item cap-cap: $ \G_{n+2}^{i+k} \circ \G_n^i \cong \G_{n+2}^i \circ \G_n^{i+k-2}$ for $k \ge 2$. 
\item cup-cup: $\F_n^{i+k-2} \circ \F_{n+2}^i \cong \F_n^i \circ \F_{n+2}^{i+k}$ for $k \ge 2$.
\item cup-cap: $\G_n^{i+k-2} \circ \F_n^i \cong \F_{n+2}^i \circ \G_{n+2}^{i+k}$ and $\G_n^i \circ \F_n^{i+k-2} \cong \F_{n+2}^{i+k} \circ \G_{n+2}^i$ for $k \ge 2$. 
\end{itemize}
\end{Proposition}
\begin{proof}
The kernel of $\G_{n+2}^{i+k} \circ \G_n^{i}: D(Y_{n-2}) \rightarrow D(Y_{n+2}) $ is given by 
$$\sG_{n+2}^{i+k} \ast \sG_n^{i} = \pi_{13\ast} \left( \pi_{12}^\ast(\sG_n^{i}) \otimes \pi_{23}^\ast(\sG_{n+2}^{i+k}) \right) \in D(Y_{n-2} \times Y_{n+2}).$$

We have $ \sG_n^i = \O_{X_n^i} \otimes \E'_i \{-i+1\} $ and $\sG_{n+2}^{i+k} = \O_{X_{n+2}^{i+k}} \otimes \E'_{i+k} \{-i-k+1\}$. Let $ W = \pi_{12}^{-1}(X_n^i) \cap \pi_{23}^{-1}(X_{n+2}^{i+k}) $.  By Corollary \ref{cor:transverse}, this intersection is transverse and so by Lemma \ref{th:tensor}, we have 
\begin{equation*}
\pi_{12}^*(\O_{X_n^i}) \otimes \pi_{23}^*(\O_{X_{n+2}^{i+k}}) \cong \O_W.
\end{equation*}

A routine computation shows that
\begin{equation*}
\begin{aligned}
W = \{ (L_\cdot, L'_\cdot, L''_\cdot) \in Y_{n-2} \times Y_n \times Y_{n+2} :\ &L_j = L'_j \text{ for } j \le i-1,\ L_j = zL'_{j+2} \text{ for } j \ge i-1, \\
&L'_j = L''_j \text{ for } j \le i+k-1,\ L'_j = zL''_{j+2} \text{ for } j \ge i+k-1 \}.
\end{aligned}
\end{equation*}
Hence on $ W $ there is an isomorphism of line bundles 
\begin{equation} \label{eq:lineiso2}
 \E'_i \cong \E''_i 
 \end{equation}

The fibres of $ \pi_{13} $ restricted to $ W $ are points and 
\begin{equation*}
\begin{aligned}
\pi_{13}(W) = \{ (L_\cdot, L'_\cdot) \in Y_{n-2} \times Y_{n+2} :\ &L_j = L'_j \text{ for } j \le i-1,\ L_j = zL'_{j+2} \text{ for } i-1 \le j \le i + k -1,\\
&L_j = z^2L'_{j+4} \text{ for } j \ge i+k -1 \}.
\end{aligned}
\end{equation*}

Hence
\begin{equation} \label{eq:lhs2}
\begin{aligned}
{\pi_{13}}_\ast (\pi_{12}^\ast(\sG_n^i) \otimes \pi_{23}^\ast(\sG_n^{i+k})) &= {\pi_{13}}_\ast( \O_{\pi_{12}^{-1}(X_n^i)} \otimes \O_{\pi_{23}^{-1}(X_{n+2}^{i+k})} \otimes \E'_i \otimes \E''_{i+k} ) \{-2i-k+2\} \\
&\cong {\pi_{13}}_\ast( \O_W \otimes  \E'_i \otimes \E''_{i+k}) \{-2i-k+2\} \\
&\cong {\pi_{13}}_\ast(\O_W \otimes \E''_i \otimes \E''_{i+k}) \{-2i-k+2\} \\
&= \O_{\pi_{13}(W)} \otimes  \E'_i \otimes \E'_{i+k} \{-2i-k+2\}
\end{aligned}
\end{equation}
where in the second last step we used the above isomorphism (\ref{eq:lineiso2}).

Now, we want to compute $ \sG_{n+2}^i \ast \sG_n^{i+k-2} $. Again the intersection $ W' := \pi_{12}^{-1}(X_n^{i+k-2}) \cap \pi_{23}^{-1}(X_{n+2}^i) $ is transverse.  We have
\begin{equation*}
\begin{aligned}
W' = \{ (L_\cdot, L'_\cdot, L''_\cdot) :\ &L_j = L'_j \text{ for } j \le i+k-3,\ L_j = zL'_{j+2} \text{ for } j \ge i+k-3, \\
&L'_j = L''_j \text{ for } j \le i-1,\ L'_j = zL''_{j+2} \text{ for } j \ge i+1 \}.
\end{aligned}
\end{equation*}
Hence on $ W' $ the map $ z $ induces an isomorphism of line bundles
\begin{equation} \label{eq:lineiso3}
\E'_{i+k-2}\{2\} \cong \E''_{i+k}.
\end{equation}

Also we see that the fibres of $ \pi_{13} $ restricted to $ W' $ are points and $ \pi_{13}(W') = \pi_{13}(W)$. Hence
\begin{equation} \label{eq:rhs2}
\begin{aligned}
{\pi_{13}}_\ast (\pi_{12}^\ast(\sG_n^{i+k-2}) \otimes \pi_{23}^\ast(\sG_{n+2}^i)) &= {\pi_{13}}_\ast( \O_{\pi_{12}^{-1}(X_n^{i+k-2})} \otimes \O_{\pi_{23}^{-1}(X_{n+2}^i)} \otimes \E'_{i+k-2} \otimes \E''_i ) \{-2i-k+4\} \\
&\cong {\pi_{13}}_\ast( \O_{W'} \otimes \E'_{i+k-2} \otimes \E''_i) \{-2i-k+4\} \\
&\cong {\pi_{13}}_\ast( \O_{W'} \otimes \E''_{i+k} \{-2\} \otimes \E''_i) \{-2i-k+4\}\\
&= \O_{\pi_{13}(W)} \otimes \E'_i \otimes \E'_{i+k} \{-2i-k+2\}
\end{aligned}
\end{equation}
where in the second last step we have used the isomorphism (\ref{eq:lineiso3}).

Since (\ref{eq:lhs2}) and (\ref{eq:rhs2}) agree, we conclude that $ \sG_{n+2}^{i+k} \ast \sG_n^i \cong \sG_{n+2}^i \ast \sG_n^{i+k-2} $ and hence that $ \G_{n+2}^{i+k} \circ \G_n^i \cong \G_{n+2}^i \circ \G_n^{i+k-2} $.

The other equations all follow similarly.
\end{proof}

\begin{Proposition} \label{prop:capcupcross}
The functor $ \Psi $ is invariant with respect to isotopies exchanging the order with respect to height of cap/cups with crossings.
\begin{itemize}
\item cap-crossing: $ \G_n^i \circ \T_{n-2}^{i+k-2}(l) \cong \T_n^{i+k}(l) \circ \G_n^i $ for $ k \ge 2, 1 \le l \le 4 $.
\item crossing-cap:  $ \G_n^{i+k} \circ \T_{n-2}^i(l) \cong \T_n^i(l) \circ \G_n^{i+k} $ for $ k \ge 0, 1 \le l \le 4 $.
\item cup-crossing: $ \F_n^i \circ \T_n^{i+k}(l) \cong \T_{n-2}^{i+k-2}(l) \circ \F_n^i $ for $ k \ge 2, 1 \le l \le 4  $.
\item crossing-cup:  $ \F_n^{i+k} \circ \T_n^i(l) \cong \T_{n-2}^i(l) \circ \F_n^{i+k} $ for $ k \ge 0, 1 \le l \le 4  $.
\end{itemize}
\end{Proposition}

\begin{proof}
We will prove the first equation with $ l=1$.

First, note that $ \G_n^i \circ \T_{n-2}^{i+k-2}(1) $ is the Fourier-Mukai transform with respect to the kernel
\begin{eqnarray*}
\sG_n^i \ast \sT_{n-2}^{i+k-2}(1) &=& {\pi_{13}}_*(\pi_{12}^*(\O_{Z_{n-2}^{i+k-2}}[1]\{-1\}) \otimes \pi_{23}^* (\O_{X_n^i} \otimes \E'_i \{-i+1\})) \\
&=& \pi_{13\ast}(\O_{\pi_{12}^{-1}(Z_{n-2}^{i+k-2})} \otimes \O_{\pi_{23}^{-1}(X_n^i)}) \otimes \E'_i[1] \{-i\}.
\end{eqnarray*}

By Corollary \ref{cor:transverse}, each component of $\pi_{12}^{-1}(Z_{n-2}^{i+k-2}) $ intersects $ \pi_{23}^{-1}(X_n^i) $ transversely.  Hence by Lemma \ref{th:tensor}, we have that 
\begin{equation*}
\O_{\pi_{12}^{-1}(Z_{n-2}^{i+k-2})} \otimes \O_{\pi_{23}^{-1}(X_n^i)} \cong \O_W
\end{equation*}
where $ W = \pi_{12}^{-1}(Z_{n-2}^{i+k-2}) \cap \pi_{23}^{-1}(X_n^i) $.  We have that
\begin{equation*}
W = \{ (L_\cdot, L'_\cdot, L''_\cdot) : L_j = L'_j \text{ for }  j \ne i+k-2, L'_j = L''_j \text{ for } j \le i, L'_j = zL''_{j+2} \text{ for } j \ge i \}.
\end{equation*}
The fibres of $ \pi_{13} $ restricted to $ W $ are just points and so we conclude that
\begin{equation*}
\sG_n^i \ast \sT_{n-2}^{i+k-2}(1) \cong \O_{\pi_{13}(W)} \otimes \E'_i[1]\{-i\}.
\end{equation*}

On the other hand, using transversality again, we see that
\begin{equation*}
\sT_n^{i+k}(1) \ast \sG_n^i \cong \O_{\pi_{13}(W')} \otimes \E'_i[1]\{-i\},
\end{equation*}
where $ W' = \pi_{12}^{-1}(X_n^i) \cap \pi_{23}^{-1}(Z_n^{i+k})$.  We see that 
\begin{equation*}
W' = \{ (L_\cdot, L'_\cdot, L''_\cdot) : L_j = L'_j \text{ for } j \le i,\ L_j = zL'_{j+2} \text{ for } j \ge i,\ L'_j = L''_j \text{ for } j \ne i+k \}.
\end{equation*}

A quick inspection of the definitions, shows that 
\begin{equation*}
\pi_{13}(W) = \pi_{13}(W') = \{ (L_\cdot, L'_\cdot) : L_j = L'_j \text{ for } j \le i, L_j = zL'_{j+2} \text{ for } j \ge i, j \ne i+k-2 \}
\end{equation*}
whence we conclude that $ \sG_n^i \ast \sT_{n-2}^{i+k-2}(1) \cong \sT_n^{i+k}(1) * \sG_n^i $. The other equations follow similarly. 
\end{proof}

\begin{Proposition} \label{prop:tangleiso}
The functor $ \Psi $ is invariant with respect to isotopies exchanging the heights of crossings.  We have 
$$ \T_n^i(l) \circ \T_n^j(m) \cong \T_n^j(m) \circ \T_n^i(l) \text{ for } 1 \le l,m \le 4 \text{ and } | i -j | \ge 2. $$
\end{Proposition}

\begin{proof}
We will prove the case $ l=m=1 $ (the other cases being similar).

As before, we compute the convolution of the Fourier-Mukai kernels.  Recall that $ Z_n^i = \Delta \cup \othercomp$ where $\othercomp = X_n^i \times_{Y_{n-2}} X_n^i $.  By Corollary \ref{cor:transverse}, we see that each component of $ Z_n^i $ is transverse to each component of $ Z_n^j $.  Hence by Lemma \ref{th:tensor}, we have 
\begin{equation*}
\begin{aligned}
\sT_n^i(1) * \sT_n^j(1) &= {\pi_{13}}_* (\pi_{12}^*(\O_{Z_n^j}[1]\{-1\}) \otimes \pi_{23}^*(\O_{Z_n^i}[1]\{-1\})) \\
&\cong {\pi_{13}}_\ast(\O_W)[2]\{-2\} = \O_{\pi_{13}(W)}[2]\{-2\}
\end{aligned}
\end{equation*}
where $ W = \pi_{12}^{-1}(Z_n^j) \cap \pi_{23}^{-1}(Z_n^i) = \{ (L, L', L'') : L_k = L'_k \text{ for } k \ne j, L'_k = L''_k \text{ for } k \ne i  \} $.

A similar computation shows that 
\begin{equation*}
\sT_n^j(1) * \sT_n^i(1) \cong \O_{\pi_{13}(W')}[2]\{-2\}
\end{equation*}
where $ W' =  \{ (L_\cdot, L_\cdot', L_\cdot'') : L_k = L'_k \text{ for } k \ne i, L'_k = L''_k \text{ for } k \ne j  \} $.

Since $ |i-j| \ge 2$, we see that
\begin{equation*}
\pi_{13}(W) = \pi_{13}(W') = \{ (L_\cdot, L_\cdot') : L_k = L'_k \text{ for } k \ne i, j \},
\end{equation*}
and so the result follows.
\end{proof}

\subsection{Invariance under pitchfork move}


\begin{Proposition}
$\T_{n}^i(1) \circ \G_n^{i+1} \cong \T_{n}^{i+1}(4) \circ \G_n^{i}, \quad \T_{n}^i(2) \circ \G_n^{i+1} \cong \T_{n}^{i+1}(3) \circ \G_n^i$
\end{Proposition}

\begin{proof}
We must show that 
\begin{equation*}
\sT_{n}^{i}(1) * \sG_n^{i+1} \cong \sT_{n}^{i+1}(4) * \sG_n^i \text{ and } 
\sT_{n}^i(2) * \sG_n^{i+1} \cong \sT_{n}^{i+1}(3) * \sG_n^i
\end{equation*}
But by Theorem \ref{th:kerneltwist}, 
$$ \sT_{n}^i(1) \cong ({\sT_{\sG_{n}^i}})_R[1]\{-1\}, \ \sT_{n}^{i+1}(4) \cong \sT_{\sG_n^{i+1}}[1]\{-2\}, \quad 
\sT_{n}^i(2) \cong \sT_{\sG_{n}^i}[-1]\{1\}, \ \sT_{n}^{i+1}(3) \cong ({\sT_{\sG_{n}^{i+1}}})_R[-1]\{2\}. $$  
So we need to show that
$$({\sT_{\sG_{n}^i}})_R[1]\{-1\} * \sG_n^{i+1} \cong \sT_{\sG_n^{i+1}}[1]\{-2\} * \sG_n^i \ \text{ and } \ 
\sT_{\sG_{n}^i}[-1]\{1\} * \sG_n^{i+1} \cong ({\sT_{\sG_{n}^{i+1}}})_R[-1]\{2\} * \sG_n^i$$
which follows from Lemma \ref{th:forpitch}.
\end{proof}

\section{The Induced Action on Grothendieck groups} \label{se:Kgroups}

The goal of this section is to prove that on the level of the Grothendieck groups, the geometrically defined $ \Psi(T) $ induces the representation theoretic map $ \psi(T) $.

\subsection{Representation theoretic map}
We begin by examining more closely the Reshetikhin-Turaev invariant of tangles, as discussed in the introduction.

Let $ V $ denote the standard representation of $\uq $.  It is a free $ \base $-module with  basis $ v_0, v_1 $.  The action of the generators of $ \uq $ on this basis are
\begin{gather*}
E v_0 = 0, \ F v_0 = v_1, \ K v_0 = q v_0, \\
E v_1 = v_0, \ F v_1 = v_0, \ K v_1 = q^{-1} v_1.
\end{gather*}

We will now write down explicitly a version of the map
\begin{equation*}
\psi : \biggl\{ (n,m) \text{ tangles } \biggr\} \rightarrow \Hom_{\uq}(V^{\otimes n}, V^{\otimes m}).
\end{equation*}
discussed in the introduction.  We must define $ \psi(f_n^i), \psi(g_n^i), \psi(t_n^i(l)) $.  However, since $$ \psi(g_n^i) = id_{V^{\otimes i-1}} \otimes \psi(g_2^1) \otimes id_{V^{\otimes n-2-i+1}} $$ and similarly for the caps and crossings, it suffices to define $ \psi(g_2^1), \psi(f_2^1), \psi(t_2^1(l)) $.  For the caps and cups, we define
\begin{gather*}
\psi(g_2^1) : \base \rightarrow V \otimes V \\
1 \mapsto q^{-1} v_1 \otimes v_0 - v_0 \otimes v_1 
\end{gather*} 
\begin{gather*}
\psi(f_2^1) : V \otimes V \rightarrow \base \\
v_0 \otimes v_0 \mapsto 0,\ v_0 \otimes v_1 \mapsto q,\ v_1 \otimes v_0 \mapsto -1, \ v_1 \otimes v_1 \mapsto 0 
\end{gather*}
For the crossings, we begin by defining the braiding
\begin{gather*}
\beta : V\otimes V \rightarrow V\otimes V \\
v_0 \otimes v_0 \mapsto q^{1/2} v_0 \otimes v_0,\ v_0 \otimes v_1 \mapsto q^{-1/2} v_1 \otimes v_0, \\
 v_1 \otimes v_0 \mapsto q^{-1/2} (v_0 \otimes v_1 + (q - q^{-1})v_1 \otimes v_0, \ v_1 \otimes v_1 \mapsto q^{1/2} v_1 \otimes v_1 
\end{gather*}
and then define the crossings in terms of the braiding as
\begin{equation*}
\psi(t_2^1(1)) := -q^{3/2} \beta^{-1}, \quad \psi(t_2^1(2)) := -q^{-3/2} \beta, \quad \psi(t_2^1(3)) := - q^{-3/2} \beta^{-1}, \quad \psi(t_2^1(4)) := q^{3/2} \beta
 \end{equation*}
Note in particular we have the ``Kauffman'' relation 
\begin{equation} \label{eq:Kauff}
 \psi(t_2^1(2)) = -q^{-1}\beta =  -q^{-1}(q^{-1} \psi(g_1^1) \circ \psi(f_1^1) +  id) 
 \end{equation}

\begin{Remark} \label{re:Resh}
These formulas can be obtained from the Reshtikhin-Turaev invariants \cite{RT}, with the following caveats.  First, we are identifying $ V^\vee $ with $ V $ via the $ \uq$-equivariant isomorphism $ v^0 \mapsto -v_1 $ and $ v^1 \mapsto q^{-1} v_0 $ (here $v^0, v^1 $ denote the dual basis).  Second, as our ``ribbon'' element $ v $, we are using the negative of the usual element.  This has the effect of negating the maps assigned to right moving caps and cups, and the maps assigned to all of the crossings.
\end{Remark}

\subsection{Equivariant K-theory of $ Y_n $}
First, note that the Grothendieck group $ K(D(Y_n))$ is a module over $ \base $, where $ q [\sF] = [\sF\{-1\}] $.

To calculate $ K(D(Y_n)) $ our main tool will be a result from Chriss-Ginzburg \cite{CG}, called the cellular fibration lemma.  Using this result we will establish a basis for $ K(D(Y_n)) $.
\begin{Theorem}
$ K(D(Y_n)) $ is a free $ \base $ module with basis 
\begin{equation*}
\{ \prod_j [\E_j]^{\delta_j} \}_{\delta \in \{0,1\}^{n}}
\end{equation*}
\end{Theorem}

\begin{proof}
We proceed by induction on $ n $.  As the base case, $ Y_1 \cong \p $ with the trivial $ \C^\times $ action.  So  $ K(D(Y_1)) $ is a free $ \base $ module with basis $ [\O], [\E_i]$.

We have a $ \C^\times $ equivariant fibration $ \pi: Y_n \rightarrow Y_{n-1} $.  This map admits an equivariant section 
\begin{equation*}
\begin{aligned}
s: Y_{n-1} &\rightarrow Y_n \\
(L_1, \dots, L_{n-1}) &\mapsto (L_1, \dots, L_{n-1}, z^{-1}L_{n-2}).
\end{aligned}
\end{equation*}
The image of this section is $ X_n^{n-1} $.  The restriction of $ \pi $ to the compliment $ \pi : Y_n \smallsetminus X_n^{n-1} \rightarrow Y_{n-1} $ is an affine fibration, that is, it is a locally trivial bundle with fibres isomorphic to $ \C $ and transition functions affine linear.

Hence the fibration $ \pi $ along with the filtration $ Y_n \supset X_n^{n-1} \supset \emptyset $ is a cellular fibration in the sense of \cite{CG}, section 5.5.  By induction $K(D(Y_{n-1})) $ is a free $\base $ module with basis
$ \{ \prod [\E_j]^{\delta_j} \}_{\delta \in \{0,1\}^{n-1}} $.  Thus by Lemma 5.5.1 of \cite{CG}, $K(D(Y_n))$ is a free $\base $ module with basis 
\begin{equation*}
\{ \pi^*(\prod_{j < n} [\E_j]^{\delta_j}) \}_{\delta \in \{0,1\}^{n-1}} \cup \{ s_*(\prod_{j < n} [\E_j]^{\delta_j}) \}_{\delta \in \{0,1\}^{n-1}}.
\end{equation*}
Note that, $ \pi^* \E_i = \E_i $ by our usual notation, while $ s_* \E_i = \O_{X_n^{n-1}} \otimes \E_i $.

Hence $ K(D(Y_n)) $ is a free $ \base $ module with basis
\begin{equation} \label{eq:basis1}
\{ \prod_{j < n} [\E_j]^{\delta_j} \}_{\delta \in \{0,1\}^{n-1}} \cup \{ \prod_{j < n} [\E_j]^{\delta_j} [\O_{X_n^{n-1}}] \}_{\delta \in \{0,1\}^{n-1}}.
\end{equation}

Now, we have the short exact sequence
\begin{equation*}
0 \rightarrow \O_{Y_n}(-X_n^{n-1}) \rightarrow \O_{Y_n} \rightarrow \O_{X_n^{n-1}} \rightarrow 0
\end{equation*}
and also we know by Lemma \ref{lem:facts} that $ \O_{Y_n}(-X_n^{n-1}) = \E_{n-1}^\vee \otimes \E_n \{-2\} $.  Hence we conclude that 
\begin{equation*}
[\O_{X_n^{n-1}}] = [\O] - q^2[\E_{n-1}^\vee][\E_n].
\end{equation*}

Thus by subtracting the first kind of basis vector from the second, then scaling the second kind of basis vector and then splitting into two cases depending on $ \delta_{n-1} $, we see that 
\begin{equation} \label{eq:basis2}
\{ \prod_{j < n} [\E_j]^{\delta_j} \}_{\delta \in \{0,1\}^{n-1}} \cup \{ \prod_{j < n-1} [\E_j]^{\delta_j} [\E_n] \}_{\delta \in \{0,1\}^{n-2}} \cup  \{ \prod_{j < n-1} [\E_j]^{\delta_j} [\E_{n-1}^\vee] [\E_n] \}_{\delta \in \{0,1\}^{n-2}}
\end{equation}
is a basis.

We now claim it suffices to show that 
\begin{equation} \label{eq:toshow}
[\E_{n-1}^\vee] = c[\E_{n-1}] + [V]
\end{equation}
where $ c \in \base $ is invertible and where $ [V] $ is pulled back from $ Y_{n-2} $.  Indeed if this is the case, then by induction for any $ \delta \in \{0,1\}^{n-2} $, $ \prod_{j < n-1} [\E_j]^{\delta_j} [V] $ is a linear combination of products of distinct $ [E_j] $ for $ j < n-1 $.  Thus, $ \prod_{j < n-1} [\E_j]^{\delta_j} [V][\E_n] $ is a linear combination of the second kind of basis vector in (\ref{eq:basis2}).  Hence by adding multiples of the second kind of basis vector to the third kind of basis vector, we see that we can replace the third kind by 
\begin{equation*}
\prod_{j < n-1} [\E_j]^{\delta_j} c[\E_{n-1}] [\E_n].
\end{equation*}
Since $ c $ is invertible we can clear the $ c $ to obtain our desired basis.

Thus, it remains to establish (\ref{eq:toshow}).  Consider the short exact sequence
\begin{equation} \label{eq:shortexact}
0 \rightarrow \E_{n-1} \rightarrow z^{-1} L_{n-2} / L_{n-2} \rightarrow z^{-1} L_{n-2} / L_{n-1} \rightarrow 0
\end{equation}
This gives us 
\begin{equation} \label{eq:tryinginverse}
\E_{n-1} \otimes z^{-1} L_{n-2} / L_{n-1} \cong \det(z^{-1} L_{n-2} / L_{n-2}) 
\cong \det (z^{-1} L_{n-2}) \otimes (\det L_{n-2})^\vee .
\end{equation}

On the other hand we also have the short exact sequence
\begin{equation*}
0 \rightarrow \ker z \rightarrow z^{-1} L_{n-2} \rightarrow L_{n-2}\{2\} \rightarrow 0
\end{equation*}
where the third map is the action of $ z $.  This gives us
\begin{equation*}
\det(\ker z) \otimes \det (L_{n-2} \{2 \}) \cong \det (z^{-1} L_{n-2}).
\end{equation*}
Since $\ker z$ is spanned by $e_1$ and $f_1$ we find $\ker z = \O^2 \{2\}$ and hence $ \det (\ker z) \cong \O \{4\}$. We conclude that 
\begin{equation*}
\det (z^{-1} L_{n-2}) \otimes (\det L_{n-2})^\vee \cong \O \{2n\}.
\end{equation*} 

Substituting this back into (\ref{eq:tryinginverse}) we see that
\begin{equation*}
 z^{-1} L_{n-2} / L_{n-1}\{-2n\} \cong \E_{n-1}^\vee .
\end{equation*}

On the other hand (\ref{eq:shortexact}) also gives us that
\begin{equation*}
[\E_{n-1}] + [z^{-1} L_{n-2} / L_{n-1}] = [z^{-1} L_{n-2} / L_{n-2}]
\end{equation*}
and so we see that
\begin{equation*}
[\E_{n-1}]^\vee = q^{2n}(-[\E_{n-1}] + [z^{-1} L_{n-2} / L_{n-2}]).
\end{equation*}
Since $ z^{-1} L_{n-2} / L_{n-2} $ is pulled back from $ Y_{n-2} $ we have shown (\ref{eq:toshow}) as desired.
\end{proof}

\begin{Remark}
In the non-equivariant case, a much easier proof is available.  First we see that $ K(D(Y_n)) \cong K(Y_n) $, the topological K-homology of $ Y_n $.  Then we use the (non-equivariant) isomorphism of manifolds $ Y_n \cong (\p)^{n} $ (Theorem \ref{th:isoman}), to see that 
\begin{equation*}
K(D(Y_n)) \cong K((\p)^n) \cong (K(\p))^{\otimes n}.
\end{equation*}  The basis we have given comes from the natural basis of the right hand side.
\end{Remark}

\subsection{Action of basic functors on the Grothendieck groups}
We now examine actions of the caps and cups on our basis for $ K(Y_n) $.
\begin{Proposition} \label{th:capcupK}
We have 
\begin{gather*}
[\G_n^i](\prod_j [\E_j]^{\delta_j}) = q^{i-1 + 2 \sum_{j \ge i} \delta_j} \prod_{j < i} [\E_j]^{\delta_j} ([\E_i] - q^2 [\E_{i+1}]) \prod_{j\ge i} [\E_{j+2}]^{\delta_j}, \\
[\F_n^i](\prod_j [\E_j]^{\delta_j}) = 
\begin{cases}
0& \text{ if $ \delta_i = \delta_{i+1} $}, \\
q^{-i - 2 \sum_{j \ge i+2} \delta_j} \prod_{j < i} [\E_j]^{\delta_j} \prod_{j \ge i+2} \E_{j-2}^{\delta_j}& \text{ if $\delta_i = 0 $ and $ \delta_{i+1} = 1 $}, \\
- q^{-i - 2 \sum_{j \ge i+2} \delta_j} \prod_{j < i} [\E_j]^{\delta_j} \prod_{j \ge i+2} \E_{j-2}^{\delta_j}& \text{ if $ \delta_i = 1 $ and $ \delta_{i+1} = 0 $}. \\
\end{cases}
\end{gather*}
\end{Proposition}

\begin{proof}
First, we consider $ [\G_n^i] $.  By definition $ \G_n^i(\sF) = i_* q^* \sF \otimes \E_i \{-i+1\} $.  Note that if $ j < i $, then $ \E_j \cong \E'_j $ on $ X_n^i $, while if $ j \ge i $, then $ \E_j \{2\} \cong \E'_{j+2} $ on $ X_n^i $.  Hence we see that 
\begin{equation} \label{eq:Gonbasis}
\G_n^i(\bigotimes_j \E_j^{\delta_j}) = \bigotimes_{j < i} \E_j^{\delta_j} \otimes \bigotimes_{j \ge i} \E_{j+2}^{\delta_j} \otimes \O_{X_n^i} \otimes \E_i \{-i+1 - 2\sum_{j \ge i} \delta_j\} .
\end{equation}
Hence it will be necessary to examine $ [\O_{X_n^i}] $ in terms of our basis.  From the standard short exact sequence \begin{equation*}
0 \rightarrow \O_{Y_n}(-X_n^i) \rightarrow \O_{Y_n} \rightarrow \O_{X_n^i} \rightarrow 0
\end{equation*}
we have that
\begin{equation*} 
[\O_{X_n^i}] = [\O_{Y_n}] - [\O_{Y_n}(-X_n^i)] = [\O_{Y_n}] - [\E_{i+1} \otimes \E^\vee_i \{-2\}].
\end{equation*}
Using this relation, we can see that (\ref{eq:Gonbasis}) gives us
\begin{equation*}
[\G_n^i(\bigotimes_j \E_j^{\delta_j})] = q^{i-1 + 2 \sum_{j \ge i} \delta_j} \prod_{j < i} [\E_j]^{\delta_j} \prod_{j \ge i} [\E_{j+2}]^{\delta_j} ([\O_{Y_n}] - [\E_{i+1} \otimes \E^\vee_i\{-2\}])[\E_i]
\end{equation*}
and so first result follows.

For the calculation of $ [\F_n^i] $, we recall that $ \F_n^i(\sF) = q_* (i^* \sF \otimes \E_{i+1}^\vee)\{i\} $. Now $ \E_i $ and $ \E_{i+1}^\vee $ restrict to $ O_{\p}(-1) $ on the fibres of $ q $, while $ \E_i \otimes \E_{i+1}^\vee $ is isomorphic to the relative dualizing bundle by Lemma \ref{lem:facts}.  Hence we have
\begin{gather*}
\F_n^i(\O) = q_* \E_{i+1}^\vee \{i\} = 0, \quad \F_n^i(\E_i) = q_* \E_i \otimes \E_{i+1}^\vee \{i\} = \O[-1]\{i\}, \\
\F_n^i(\E_{i+1}) = q_* \O \{i\} = \O\{i\}, \quad \F_n^i(\E_i \otimes \E_{i+1}) = q_* \E_i \{i\} = 0.
\end{gather*}

From this knowledge, it is easy to deduce the desired description of $ [\F_n^i] $.
\end{proof}

\subsection{Comparison of Grothendieck group with representation}
Define a map of $ \base $ modules
\begin{equation*}
\begin{aligned}
\alpha : K(D(Y_n)) &\rightarrow V^{\otimes n} \\
\prod_i [\E_i]^{\delta_i} &\mapsto q^{-\sum_i i \delta_i} v_{\delta_1} \otimes \dots \otimes v_{\delta_{n}}.
\end{aligned}
\end{equation*}

\begin{Theorem} \label{th:Kgroup}
For any $(n,m) $ tangle $ T $, the following diagram commutes:
\begin{equation} \label{eq:commute?}
\begin{CD}
K(D(Y_n)) @>[\Psi(T)]>> K(D(Y_m)) \\
@V\alpha VV @VV\alpha V \\
V^{\otimes n} @>>\psi(T)> V^{\otimes m}
\end{CD}
\end{equation}
\end{Theorem}

\begin{proof}
It suffices to prove the statement when $ T $ is a basic tangle, ie. when $T = f_n^i, g_n^i, t_n^i(l) $.

However, note that by (\ref{eq:Kauff}), we have $\psi(t_n^i(2)) =  -q^{-1} ( q^{-1}\psi(g_n^i) \circ \psi(f_n^i) + id)$.  Also since $ \Psi(t_n^i(2)) $ is the twist in $ \G_n^i $ shifted by $ [-1]\{1\} $ (by Theorem \ref{th:kerneltwist}),  we have that 
\begin{eqnarray*}
[\Psi(t_n^i(2))] &=& -q^{-1}(-[\G_n^i \circ (\G_n^i)^R] + [id]) \\
&=& -q^{-1}(-[\G_n^i \circ \F_n^i[-1]\{1\}] + id) = -q^{-1}(q^{-1}[\Psi(g_n^i)] \circ [\Psi(f_n^i)] + id).
\end{eqnarray*}
A similar argument for the other $t_n^i(l)$ shows that it actually suffices to prove the statement when $ T $ is $g_n^i $ or $ f_n^i $.

When $ T = g_n^i $, using Proposition \ref{th:capcupK}, following (\ref{eq:commute?}) right and down we find
\begin{equation*}
\begin{aligned}
\prod_j [\E_j]^{\delta_j} &\mapsto q^{i-1 + 2 \sum_{j \ge i} \delta_j} \prod_{j < i} [\E_j]^{\delta_j} ([\E_i] - q^2 [\E_{i+1}]) \prod_{j\ge i} [\E_{j+2}]^{\delta_j} \\
&\mapsto q^{i-1 + 2 \sum_{j \ge i} \delta_j} q^{-\sum_{j < i} j\delta_j - \sum_{j \ge i} (j+2)\delta_j} v_{\delta_1} \otimes \dots \otimes  v_{\delta_{i-1}} \otimes \\
&\phantom{\mapsto q^{i-1 + 2 \sum_{j \ge i} \delta_j}}(q^{-i} v_1 \otimes v_0 - q^{-i-1} q^2 v_0 \otimes v_1) \otimes v_{\delta_i} \otimes \dots \otimes v_{\delta_{n-2}}\\
&= q^{- \sum_j j \delta_j} v_{\delta_1} \otimes \dots \otimes  v_{\delta_{i-1}} \otimes (q^{-1} v_1 \otimes v_0 - v_0 \otimes v_1) \otimes v_{\delta_i} \otimes \dots \otimes v_{\delta_{n-2}}.
\end{aligned}
\end{equation*}

On the other hand, following (\ref{eq:commute?}) down and right we find
\begin{equation*}
\begin{aligned}
\prod_j [\E_j]^{\delta_j} &\mapsto q^{-\sum_j j \delta_j} v_{\delta_1} \otimes \dots \otimes v_{\delta_{n-2}} \\
&\mapsto q^{-\sum_j j \delta_j} v_{\delta_1} \otimes \dots \otimes  v_{\delta_{i-1}} \otimes (q^{-1} v_1 \otimes v_0 - v_0 \otimes v_1) \otimes v_{\delta_i} \otimes \dots \otimes v_{\delta_{n-2}}.
\end{aligned}
\end{equation*}
Hence the diagram (\ref{eq:commute?}) commutes when $ T = g_n^i $.

Now, we consider the case $ T = f_n^i $.  There are four cases depending on the values of $ \delta_{i}, \delta_{i+1} $.  We will do the case $ \delta_i = 1, \delta_{i+1} = 0 $.  Using proposition \ref{th:capcupK}, following (\ref{eq:commute?}) right and down we find 
\begin{equation*}
\begin{aligned}
\prod_{j < i} [\E_j]^{\delta_j} [\E_i] \prod_{j \ge i+2} [\E_j]^{\delta_j}
&\mapsto - q^{-i - 2 \sum_{j > i+1} \delta_j} \prod_{j < i} [\E_j]^{\delta_j} \prod_{j \ge i+2} [\E_{j-2}]^{\delta_j} \\
&\mapsto -q^{-i-2 \sum_{j \ge i+2} \delta_j} q^{-\sum_{j <i} j \delta_j - \sum_{j \ge i+2} (j-2)\delta_j} v_{\delta_1} \otimes \dots v_{\delta_{i-1}} \otimes v_{\delta_{i+2}} \otimes \dots \otimes v_{\delta_{n}} \\
&= -q^{-i-\sum_j j \delta_j} v_{\delta_1} \otimes \dots \otimes v_{\delta_{i-1}} \otimes v_{\delta_{i+2}} \otimes \dots \otimes v_{\delta_{n}}.
\end{aligned}
\end{equation*}
On the other hand, following (\ref{eq:commute?}) down and right we find 
\begin{equation*}
\begin{aligned}
\prod_{j < i} [\E_j]^{\delta_j} [\E_i] \prod_{j \ge i+2} [\E_j]^{\delta_j} 
&\mapsto q^{-\sum_j j \delta_j - i} v_{\delta_1} \otimes \dots \otimes v_{\delta_{i-1}} \otimes v_1 \otimes v_0 \otimes v_{\delta_{i+2}} \otimes \dots \otimes v_{\delta_{n}} \\
&\mapsto -q^{-i-\sum_j j \delta_j} v_{\delta_1} \otimes \dots \otimes v_{\delta_{i-1}} \otimes v_{\delta_{i+2}} \otimes \dots \otimes v_{\delta_{n}},
\end{aligned}
\end{equation*}
and so the diagram commutes.  The case where $ \delta_i = 0, \delta_{i+1} = 1 $ is similar and the two cases where $ \delta_i = \delta_{i+1} $ result in both paths giving 0.
\end{proof}

\begin{Corollary}
If $ K $ is a link, then $ (-1)^{\text{\# of comp of $ K $}} \sum_{i,j} (-1)^i q^{j} \dim \Kh^{ij}(K) $ is the Jones polynomial.
\end{Corollary}

\begin{proof}
First note that by definition 
\begin{equation*}
\sum_{i,j} (-1)^i q^{j} \dim \Kh^{ij}(K) = [\Psi(K)(\mathbb{C})].
\end{equation*}
Also by Theorem \ref{th:Kgroup}, $[\Psi(K)(\C)] = \psi(K)(1) $.  

Now let $ \widehat{\psi} $ denote the ``official'' Reshetikhin-Turaev invariant, so that $ \widehat{\psi}(K)(1) $ is the Jones polynomial.  By Remark \ref{re:Resh} we see that $ \psi(K) = (-1)^{r(K)} \widehat{\psi}(K) $ where $ r(K) $ denotes the number of right moving caps, right moving cups and crossings in a projection of $ K $.  Note that $r(K) $ is a link invariant.

We claim that $ r(K) $ is the number of components of $ K $.  To see this, fix a projection of $ K $ and let $ K' $ be an unlink projection which is obtained from the projection of $ K $ by flipping crossings.  Hence $ r(K) = r(K') $.  On the other hand $ K' $ is the unlink, so it has a projection consisting of disjoint circles.  Thus, $ r(K') $ is the number of components of $ K' $.  But $ K' $ and $ K $ have the same number of components, so $ r(K) $ is the number of components of $ K $.

Thus combining everything we have that
\begin{equation*}
\psi'(K)(1) = (-1)^{\# \text{ of comp of } K} \psi(K)(1) = (-1)^{\# \text{ of comp of } K} \sum_{i,j} (-1)^i q^{j} \dim \Kh^{ij}(K).
\end{equation*}
\end{proof}

\section{An Invariant Of Tangle Cobordisms} \label{se:cobordism}

In this section we show that the functors $\Psi(T): D(Y_n) \rightarrow D(Y_m)$ associated to tangles $T$ admit natural transformations (up to scalar) $\Psi(T) \rightarrow \Psi(T')$ corresponding to tangle cobordisms $T \rightarrow T'$. To do this we use the combinatorial realization of tangle cobordisms as developed by Carter-Rieger-Saito \cite{CRS} and Baez-Langford \cite{BL}. We will follow the general argument and use the notation as well as diagrams from \cite{Kcobordism}.

\subsection{Combinatorial Model}

\begin{figure}
\begin{center}
\psfrag{birth}{birth}
\psfrag{death}{death}
\psfrag{saddle 1}{saddle \#1}
\psfrag{saddle 2}{saddle \#2}
\puteps[0.35]{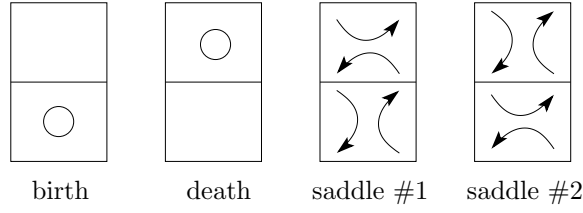}
\end{center}
\caption{Elementary movies (read from top to bottom).}\label{f4}
\end{figure}

To explain more precisely the combinatorial model for tangle cobordisms we introduce three 2-categories $\mathcal{C}, \mathcal{C'}$ and $\mathcal{T}$. In $\mathcal{C}$ the objects are non-negative integers, the 1-morphisms $[n] \rightarrow [m]$ are generic projections of $(n,m)$ tangles and the 2-morphisms are (combinatorial depictions of) cobordisms of tangle projections. The cobordisms are illustrated by movies showing the cobordism at various stages. Just as tangles are generated by elementary cups, caps and crossings the tangle cobordisms are generated by certain elementary movies: namely,
\begin{itemize}
\item movies depicting a Reidemeister move of type (0),(I),(II),(III) or a pitchfork move
\item movies depicting the shifting of distant crossings, local minima or local maxima with respect to height
\item birth, death and saddle cobordisms illustrated in figure \ref{f4}
\end{itemize}

Just as tangle projections represent the same tangle if and only if they differ by a series of elementary moves (such as Reidemeister moves) tangle movies represent the same cobordism if and only if they differ by a series of elementary movie moves. For instance, figure \ref{f5} illustrates three movie moves. The left movie move contains two movies each of which depicts a cobordism of the death of a circle. These two cobordisms live inside $S^3 \times [0,1]$ and are in fact isotopic through cobordisms. Thus the movie move reflects the fact that the cobordisms represented by the movies are isotopic to each other. All the movie moves are nicely illustrated in \cite{Kcobordism}. Since there are a total of 31 movie moves we only reproduced the three from figure \ref{f5}. 

We define the 2-category $\mathcal{C}'$ as $\mathcal{C}$ modulo the relation on 2-morphisms which identifies two movies if they differ by a series of movie moves. The tangle 2-category $\mathcal{T}$ is defined to have objects the non-negative integers, 1-morphisms $(n,m)$ tangles up to isotopy and 2-morphisms cobordisms up to isotopy. Theorem 17 in \cite{BL} shows that $\mathcal{C'}$ is a combinatorial model for $\mathcal{T}$, namely

\begin{Theorem} $\mathcal{C}'$ and $\mathcal{T}$ are isomorphic 2-categories. 
\end{Theorem}

\subsection{Cobordism Invariance of $\Psi$}

Using the spaces $Y_n$ we define the 2-category $\mathcal{D}$ as follows. The objects in $\mathcal{D}$ are $D(Y_n)$ for $n \ge 0$, the 1-morphisms are FM kernels $\P \in D(Y_n \times Y_m)$  and the 2-morphisms are (weighted) projective morphisms of FM kernels $\P \rightarrow \P'$. Here projective means that a 2-morphism is defined only up to a non-zero multiple so that a 2-morphism between $\P$ and $\P'$ is an element of $\{0\} \cup_{i,j} \mathbb{P}\left(\Hom(\P,\P'[i]\{j\})\right)$. 

Note that the 1-morphisms $ \P \in D(Y_n \times Y_m) $ induce functors $\Phi_\P: D(Y_n) \rightarrow D(Y_m)$.  Also, the 2-morphisms induce natural transformations of functors $\Phi_\P \rightarrow \Phi_{\P'}$ (up to scalar).  Hence there is a natural 2-functor from $\mathcal{D} $ to the 2-category of triangulated categories (where we quotient out the natural transformations by scalars).

We now describe a 2-functor $\Psi: \mathcal{C} \rightarrow \mathcal{D}$. To an object $[n] \in \mathcal{C}$ we associate $D(Y_n) \in \mathcal{D}$. To a 1-morphism $T$ in $\mathcal{C}$ we associate the FM kernel corresponding to the functor $\Psi(T)$ described in section \ref{sse:functor}. To define the 2-functor on 2-morphisms it is enough to associate a (weighted) morphism of FM kernels for each elementary cobordism (if $\P \rightarrow \P'[i]\{j\}$ is an equivariant map then the associated map $\P \rightarrow \P'$ has weight $[i]\{j\}$):
\begin{itemize}
\item to a movie depicting a Reidemeister move of type (0),(I),(II),(III) or a pitchfork move we choose an isomorphism between the two kernels (we get such isomorphisms from section \ref{se:invariance} and any such isomorphism is unique up to scalar multiple)
\item to a movie depicting the shifting of distant crossings or local minima or maxima we also assign an isomorphism of kernels (we get such isomorphisms from section \ref{se:invdistant} which again are unique)
\item to the birth cobordism we assign the morphism $\O_\Delta \rightarrow \O_\Delta[1]\{-1\} \oplus \O_\Delta[-1]\{1\}$ of weight $[-1]\{1\}$ given by $s \mapsto (s,0)$ 
\item to the death cobordism we assign the morphism $\O_\Delta[1]\{-1\} \oplus \O_\Delta[-1]\{1\} \rightarrow \O_\Delta$ of weight $[-1]\{1\} $ given by $(s,t) \mapsto t$ 
\item to saddle cobordism \#1 we assign the natural adjoint map $\sG_n^i \ast \sF_n^i \rightarrow \O_\Delta$ of weight $[1]\{-1\}$ coming from ${\sF_n^i}_L \ast \sF_n^i \cong \sG_n^i \ast {\sG_n^i}_R \rightarrow \O_\Delta$
\item to saddle cobordism \#2 we assign the natural adjoint map $\O_\Delta \rightarrow \sG_n^i \ast \sF_n^i$ of weight $[1]\{-1\}$ coming from $\O_\Delta \rightarrow {\sF_n^i}_R \ast \sF_n^i \cong \sG_n^i \ast {\sG_n^i}_L$
\end{itemize}

The fact that $\Psi: \mathcal{C} \rightarrow \mathcal{D}$ is a cobordism invariant of tangles is a consequence of the following result. 

\begin{Theorem}\label{thm:cobordism} The morphism of 2-categories $\Psi: \mathcal{C} \rightarrow \mathcal{D}$ factors as $\mathcal{C} \rightarrow \mathcal{C}' \rightarrow \mathcal{D}$. 
\end{Theorem} 
\begin{proof}

We first need the following lemma.

\begin{Lemma}\label{lem:homs} Let $\Q, \Q'$ be invertible kernels. If $\P$ a kernel corresponding to a planar tangle with no closed components or crossings then 
\begin{eqnarray*}
\Hom(\Q \ast \P \ast \Q',\Q \ast \P \ast \Q'[i]\{j\}) = \left\{
\begin{array}{ll}
0 & \mbox{if $i<0$}  \\
\C & \mbox{if $i=j=0$}
\end{array}
\right.
\end{eqnarray*}
Consequently, $\Hom(\Q \ast \sG_n^i \ast \sF_n^i \ast \Q', \Q \ast \O_\Delta \ast \Q' [1]\{-1\}) \cong \C \cong \Hom(\Q \ast \O_\Delta \ast \Q', \Q \ast \sG_n^i \ast \sF_n^i \ast \Q' [1]\{-1\})$. 
\end{Lemma}
\begin{proof}
Since $\Q$ is invertible $\Q_R \ast \Q \cong \O_\Delta$ and so $\Hom(\Q \ast \P,\Q \ast \P) = \Hom(\P,\Q_R \ast \Q \ast \P) \cong \Hom(\P,\P)$. A similar argument shows that it is enough to consider the case when $\Q \cong \O_\Delta \cong \Q'$. 

We can write the kernel $\P$ corresponding to a planar tangle as the convolution $\P_1 \ast \P_2$ where $\P_1$ is the convolution of many $\sG_n^i$ and $\P_2$ the convolution of many $\sF_n^i$. Now $\sF_n^i \ast {\sF_n^i}_L = \O_\Delta \oplus \O_\Delta[-2]\{2\} = V_\Delta$ so $\P_2 \ast {\P_2}_L = V_\Delta^{\ast m}$ for some non-negative integer $m$. Consequently
$$\Hom(\P,\P) = \Hom(\P_1, \P_1 \ast \P_2 \ast {\P_2}_L) = \Hom(\P_1, \P_1 \ast V_\Delta^{\ast m}).$$
On the other hand, $\P_1$ is the convolution of many $\sG_n^i$ and hence is isomorphic to a line bundle supported on a subvariety. This is clear if you recall that the associated functor $\G_n^i$ is $i_\ast(q^\ast(\cdot) \otimes \E_i\{-i+1\})$ where $q$ is flat and $i$ is an inclusion. Thus $\P_1 \cong \O_W(\mathcal{L})$ where $W$ is some variety and $\mathcal{L}$ is some line bundle. Hence $\P_1 \ast V_\Delta^{\ast m} = \oplus_{i \ge 0} \left(\O_W(\mathcal{L})[-i]\{i\}\right)^{n_i}$ where $n_0=1$. The result follows since $\Hom(\O_W,\O_W)=\C$ and $\Hom(\O_W,\O_W[i]\{j\})=0$ if $i<0$. 

Finally, $\Hom(\sG_n^i \ast \sF_n^i, \O_\Delta[1]\{-1\}) \cong \Hom(\sG_n^i, {\sF_n^i}_L[1]\{-1\}) \cong \Hom(\sG_n^i,\sG_n^i) \cong \C$ and similarly $\Hom(\O_\Delta, \sG_n^i \ast \sF_n^i[1]\{-1\}) \cong \Hom({\sG_n^i}_L, \sF_n^i[1]\{-1\}) \cong \Hom(\sF_n^i[1]\{-1\}, \sF_n^i[1]\{-1\}) \cong \C$. 
\end{proof}

The idea is to show that two movies related by a movie move induce non-zero morphisms between two kernels. If these two kernels are a pair from lemma \ref{lem:homs} then the two morphisms induced by the two movies must be equal (up to a non-zero multiple). 

A quick survey reveals that for movie moves 1--21, 23a, 25--26 the induced morphisms have weight $[0]\{0\}$ and are between pairs of isomorphic kernels of the form $\Q \ast \P \ast \Q'$ as in lemma \ref{lem:homs}. For example, movie move 14 induces an endomorphism of $\sT_n^{i+1}(l) \ast \sT_n^i(l') \ast \sG_n^{i-1}$ while movie move 19 induces an endomorphism of $\P \ast \sF_n^i$ where $\P$ is a cup, cap or twist distant from the cup $\sF_n^i$. Moreover, all these morphisms are the composition of isomorphisms and thus are non-zero. 

The movies in moves 24,29--30 induce two morphisms of weight $[1]\{-1\}$ between pairs of kernels of the form $\Q \ast \sG_n^i \ast \sF_n^i \ast \Q'$ and $\Q \ast \O_\Delta \ast \Q'$ as in lemma \ref{lem:homs}. Each of these two morphisms is the composition of several invertible maps and one non-zero (saddle) map and hence is non-zero. As usual, since by lemma \ref{lem:homs} there is a one dimensional space of weight $[1]\{-1\}$ maps between the two kernels, the two morphisms must be the same up to a non-zero multiple.

The two movies in move 28 yield morphisms of weight $[-1]\{1\}$ from $\O_\Delta$ to $\sF_n^i \ast \sG_n^i = \O_\Delta[1]\{-1\} \oplus \O_\Delta[-1]\{1\}$ of which, by lemma \ref{lem:homs}, there is a unique non-zero map. Since the movie morphisms are obtained by composing a birth morphism (which is non-zero) with two isomorphisms they are also non-zero and hence equal (up to non-zero multiple). 

Movie 31 says that if you have two tangles $T$ and $U$ and cobordisms $T \sim T'$ and $U \sim U'$ then to get a cobordism $T \circ U \sim T' \circ U'$ it does not matter in which order you apply the cobordisms. This translates to the following obvious fact about morphisms of kernels. Given kernels $\P$ and $\Q$ with morphisms $f: \P \rightarrow \P'$ and $g: \Q \rightarrow \Q'$ then $(id_{\P'} \ast g) \circ (f \ast id_\Q) = (f \ast id_{\Q'}) \circ (id_\P \ast g)$. 

Probably the most interesting movie moves are 22, 23b and 27 as shown in figure \ref{f5}. The box in movie move 22 represents either a cup, cap or crossing. Both movies yield maps as a composition 
$$\P[1]\{-1\} \oplus \P[-1]\{1\} \rightarrow \P[1]\{-1\} \oplus \P[-1]\{1\} \rightarrow \P$$
where $\P$ is $\sF_n^i, \sG_n^i$ or $\sT_n^i(l)$, the first map is an isomorphism while the second is the death morphism of weight $[-1]\{1\}$. The first map is not necessarily unique since there is potentially a map $\P[-1]\{1\} \rightarrow \P[1]\{-1\}$. Fortunately, we chose the death morphism so that the second map is zero on the $\P[1]\{-1\}$ component. Thus the composition is indeed unique (up to multiple). This explains why we chose the death map to have weight $[-1]\{1\}$. If you read the movie backwards it also explains why we chose the birth map to have weight $[-1]\{1\}$. 

Similarly, in movie move 27 the two movies yield two isomorphisms of $\sF_n^i \ast \sG_n^i = \O_\Delta[1]\{-1\} \oplus \O_\Delta[-1]\{1\}$. These isomorphisms might be different since we have the potential for maps $\O_\Delta \rightarrow \O_\Delta[2]\{-2\}$ but if we precompose with the birth morphism the composition becomes unique (up to non-zero multiple). 

Movie move 23b yields a morphism $\sG_n^i \rightarrow \sG_n^i \ast \sF_n^i \ast \sG_n^i \rightarrow \sG_n^i$ of weight $[0]\{0\}$ where the first map is a birth and the second is a saddle. It is easy to check the composition is the identity (up to non-zero multiple). 

\begin{figure}
\begin{center}
\psfrag{Move 22}{move 22}
\psfrag{Move 23b}{move 23b}
\psfrag{Move 27}{move 27}
\puteps[0.40]{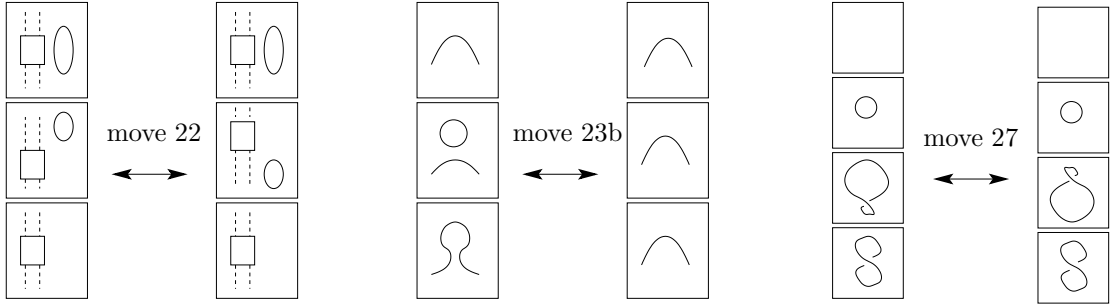}
\end{center}
\caption{Movie moves 22, 23b and 27.}\label{f5}
\end{figure}
\end{proof}

\begin{Remark}
A quick check reveals that a cobordism between two $(n,m)$ tangles has weight $[w]\{-w\}$ where $w=(n+m)/2 - \chi(S)$, $S$ is the surface forming the cobordism and $\chi(S)$ is its Euler characteristic. Moreover, a birth followed immediately by a death gives the zero morphism.
\end{Remark}

\section{Relation to Khovanov Homology} \label{se:unorient}
\subsection{Unoriented theory} \label{se:uorientdef}
In order to relate our theory to Khovanov homology, it will be convenient for us to define an unoriented version of our theory.  Let $ T $ be a $(n,m) $ unoriented tangle diagram.  We define a kernel $ \Q(T) \in D(Y_n \times Y_m)  $, by defining it on basic unoriented tangles.  There are only two unoriented crossings, denoted $ t_n^i(1), t_n^i(2) $.  We define 
\begin{gather*}
\Q(f_n^i) = \sF_n^i, \ \Q(g_n^i) = \sG_n^i, \\ 
\Q(t_n^i(1)) = (\sT_{\sG_n^i})_L = \sT_n^i(1)[-1]\{1\} = \sT_n^i(3)[1]\{-2\} = \O_{Z_n^i}, \\ 
\Q(t_n^i(2)) = \sT_{\sG_n^i}\{-1\} = \T_n^i(2)[1]\{-2\} = \T_n^i(4)[-1]\{1\} = \O_{Z_n^i} \otimes \E^\vee_{i+1} \otimes \E'_i \{1\} .
\end{gather*}

If $ T $ is an oriented tangle, then $ \Psi(T) = \Phi_{\Q(T)}[r]\{s\} $ where $ r = k_1 -k_2 - k_3 + k_4 $, $s = -k_1 + 2k_2 + 2k_3 - k_4 $ are determined by the numbers $ k_1, k_2, k_3, k_4 $ of each type of crossing in the oriented tangle $ T $.  Note a slight abuse of notation here, as $ \Psi(T) $ is defined for $ T $ an oriented tangle, while $ \Q(T) $ is defined for $ T $ an unoriented tangle diagram.  If $K$ is a link diagram then we denote by $\oH^{i,j}(K) = \mathrm{H}^{i,j}(\Q(K))$ the unoriented homology of $K$. 

\subsection{A long exact sequence}
The main advantage of the unoriented theory is that we can easily relate the kernels obtained by modifying a tangle at a crossing.  Let $ T $ be an unoriented tangle diagram and let $ s $ be a crossing in $ T $.  Let $ T_{\ver}, T_{\hor} $ denote the results of resolving the crossing $ s $ as either two vertical lines or two horizonal lines.  By the definitions above, the kernels $ \uoker(T), \uoker(T_{\ver}), \uoker(T_{\hor}) $ are of the form 
\begin{equation*}
\P_1 * (\sT_{\sG_n^i})_L * \P_2 \ \text{ or }\  \P_1 * \sT_{\sG_n^i} * \P_2\{-1\},\quad \P_1 * \P_2 ,\quad \P_1 * \sG_n^i * \sF_n^i * \P_2  
\end{equation*}
respectively.  Hence we immediately see that.
\begin{enumerate}
\item
Suppose that $ s $ is a crossing of type 1.  The natural map $ \O_\Delta \rightarrow \sG_n^i * (\sG_n^i)_L $ induces a map $ \P_1 * \P_2[-1] \rightarrow \P_1 * \sG_n^i * \sF_n^i\{-1\} * \P_2 $ and we have a distinguished triangle
\begin{equation*}
\Q(T_{\ver})[-1] \rightarrow \Q(T_{\hor})\{-1\} \rightarrow \Q(T).
\end{equation*} 
\item
Suppose that $ s $ is a crossing of type 2.  The natural map $ \sG_n^i * (\sG_n^i)_R [1]\{-1\} \rightarrow \O_\Delta [1]\{-1\} $ induces a map $ \P_1 * \sG_n^i * \sF_n^i [-1] * \P_2 \rightarrow \P_1 * \P_2 \{-1\} $ and we have a distinguished triangle
\begin{equation*}
\Q(T_{\hor})[-1] \rightarrow \Q(T_{\ver})\{-1\} \rightarrow \Q(T).
\end{equation*}
\end{enumerate}

We can view the two cases above as just a single case if we change our notation.  If $ s $ is a crossing of type 1, we refer to $T_{\ver} $ as the 0 smoothing of $ T $ at $ s $ and denote it $ T_0^s$, and we refer to $ T_{\hor} $ as the 1 smoothing of $ T $ at $ s $ and denote it $ T_1^s $.  On the other hand, if $ s $ is a crossing of type 2, we refer to $ T_{\hor} $ as the 0 smoothing of $ T $ at $ s $ and denote it $ T_0^s $ and we refer to $ T_{\ver} $ as the 1 smoothing of $ T $ at $ s$ and denote it $ T_1^s $.  This terminology matches the usual notions of 0 and 1 smoothing of a tangle at a crossing used in knot theory (see for example \cite{BN}).  With this notation, we have a distinguished triangle 
\begin{equation} \label{eq:tricross}
\Q(T_0^s)[-1] \xrightarrow{\alpha_T^s} \Q(T_1^s)\{-1\} \rightarrow \Q(T)
\end{equation}
for any tangle $ T $ and any crossing $ s $.

We can use this observation to get a long exact sequence for the cohomology of link diagrams.  

\begin{Corollary} \label{th:longexact}
Let $ K $ be a link diagram and let $ s $ be a crossing in $ K $.  There is a long exact sequence
\begin{equation*}
\dots \rightarrow \oH^{i-1,j}(K_0^s) \rightarrow \oH^{i,j-1}(K_1^s) \rightarrow \oH^{i,j}(K) \rightarrow \oH^{i,j}(K_0^s) \rightarrow \dots
\end{equation*}
\end{Corollary}

Recall that there is another knot homology theory called Khovanov homology which can be defined completely combinatorially (see \cite{Kknot}).  Let $ \Khh^{i,j}(K) $ denote Khovanov homology of a link $ K $.  Khovanov homology can be defined by first starting with an unoriented theory (which is denoted [[ ]] in \cite{BN}) and then correcting it by degree and grading shifts.  Unoriented Khovanov homology obeys the same long exact sequence as in Corollary \ref{th:longexact}.  The Khovanov homology of a loop is $ \C\{-1\} \oplus \C\{1\} $, whereas our homology of a loop is $ \C[1]\{-1\} \oplus \C[-1]\{1\}$.  

From this it is reasonable to expect a relation between Khovanov homology and our cohomology theory.  The remainder of this section is devoted to the proof of the following result.
\begin{Theorem} \label{th:Khov}
For any link $ K $,
\begin{equation*}
\Kh^{i,j}(K) = \Khh^{i+j, j}(K)
\end{equation*}
\end{Theorem}

\subsection{Convolutions of Complexes} \label{se:cones}
To prove the equivalence with Khovanov homology we will need a notion of cones for complexes of objects in triangulated categories. Thus, we recall the notion of convolution introduced by Orlov \cite{O}. 

Let $ \mathcal{D} $ be a triangulated category.  Let $ (A_\bullet, f_\bullet) = A_0 \xrightarrow{f_1} A_1 \rightarrow \cdots \xrightarrow{f_n} A_n $ be a sequence of objects and morphisms in $ \mathcal{D}$ such that $ f_{i+1} \circ f_i = 0 $. Such a sequence is called a \textbf{complex}.

A \textbf{(right) convolution} of a complex $ (A_\bullet, f_\bullet) $ is any object $ B $ such that there exist
\begin{enumerate}
\item objects $ A_0 = B_0, B_1, \dots, B_{n-1}, B_n = B $ and
\item morphisms $ g_i : B_{i-1} \rightarrow A_i $, $ h_i : A_i \rightarrow B_i $ (with $ h_0 = id $)
\end{enumerate}
such that 
\begin{equation} \label{eq:distB}
B_{i-1} \xrightarrow{g_i} A_i \xrightarrow{h_i} B_i
\end{equation}
is a distinguished triangle for each $ i $ and $ g_i \circ h_{i-1} = f_i $. Such a collection of data is called a \textbf{Postnikov system}. Notice that in a Postnikov system we also have $f_{i+1} \circ g_i = (g_{i+1} \circ h_i) \circ g_i = 0 $ since $h_i \circ g_i = 0$.

The following result is a sharper version of Lemma 1.5 from \cite{O}.
\begin{Proposition} \label{th:uniquecone}
Let $ (A_\bullet, f_\bullet) $ be a complex.  The following existence and uniqueness results hold.
\begin{enumerate}
\item
If $\Hom(A_i[k], A_{i+k+1}) = 0 $ for all $ i \ge 0, k \ge 1 $, then any two convolutions of $ (A_\bullet, f_\bullet) $ are isomorphic.

\item
If $ \Hom(A_i[k], A_{i+k+2}) = 0 $ for all $ i \ge 0, k \ge 1 $, then $ (A_\bullet, f_\bullet) $ has a convolution.
\end{enumerate}
\end{Proposition}

\begin{proof}
We start with (i).
Let $ B, B' $ be convolutions and let $ B_i, g_i, h_i, B'_i, g'_i, h'_i $ be the corresponding Postnikov systems.  We will prove that for all $ i$, $ B_i \cong B'_i $ and, using this isomorphism, that $ g_i = g'_i, h_i = h'_i $.

First, we show by induction on $ i $ that $\Hom(B_i[k], A_{i+k+1}) = 0 $ for all $ i \ge 0, k \ge 1 $.  The base case of $ i = 0 $ is covered by the hypothesis.  Now assume $ i > 0 $ and consider the long exact sequence coming from applying $ \Hom(,A_{i+k+1}) $ to the distinguished triangle (\ref{eq:distB})
\begin{equation*}
\Hom(B_{i-1}[k+1], A_{i+k+1}) \rightarrow \Hom(B_i[k], A_{i+k+1}) \rightarrow \Hom(A_i[k], A_{i+k+1})
\end{equation*}
By hypothesis $\Hom(A_i[k], A_{i+k+1}) = 0 $ and by induction $\Hom(B_{i-1}[k+1], A_{i+k+1}) = 0 $.  Thus $ \Hom(B_i[k], A_{i+k+1}) = 0 $.  Hence for all $ i \ge 0, k \ge 1 $, $ \Hom(B_i[k], A_{i+k+1}) = 0$.  In particular $ \Hom(B_{i-1}[1], A_{i+1}) = 0 $.

Now, we are in a position to prove the uniqueness statement.  By induction, we assume that $ i\ge 1 $, $ B_{i-1} = B'_{i-1} $ and $ g_i = g'_i, h_{i-1} = h'_{i-1} $.  Then both $ B_i, h_i $ and $ B'_i, h'_i $ give distinguished triangles (\ref{eq:distB}).  From the axioms of triangulated categories we see that there is an isomorphism $ B_i \cong B'_i $ intertwining $ h_i $ and $ h'_i $.  So it remains to show that this isomorphism intertwines $ g_{i+1} $ and $g'_{i+1} $.  For simplicity assume $ B_i = B'_i $.  Now consider the long exact sequence which comes from applying $ \Hom(\cdot,A_{i+1}) $ to the distinguished triangle (\ref{eq:distB}).  We see that
\begin{equation*}
\Hom(B_{i-1}[1], A_{i+1}) \rightarrow \Hom(B_i, A_{i+1}) \rightarrow \Hom(A_i, A_{i+1}).
\end{equation*}
Now, $ g_{i+1} $ and $g'_{i+1} $ lie in the middle group here and both are sent to $ f_{i+1} = g_{i+1} \circ h_i = g'_{i+1} \circ h_i  $ under the right arrow.  Since $ \Hom(B_{i-1}[1], A_{i+1}) = 0 $, this forces $ g_{i+1} = g'_{i+1} $ as desired.  By induction we deduce that $ B \cong B' $ as desired.

For (ii), we must show that there exists a chain $ B_\bullet, g_\bullet, h_\bullet $ with the desired properties.  Suppose that we have constructed $ B_0, \dots, B_i $, $ g_1, \dots, g_i$, $ h_0, \dots, h_i $.  We need to construct $ g_{i+1} $ such that $g_{i+1} \circ h_i = f_{i+1}$ and $f_{i+2} \circ g_{i+1}=0$. This is the only obstruction since we then choose $B_{i+1} = \Cone(B_i \xrightarrow{g_{i+1}} A_{i+1})$. 

As above, starting from the hypothesis $\Hom(A_i[k], A_{i+k+2}) = 0 $ for all $ i \ge 0, k \ge 1 $, we can prove by induction that $\Hom(B_i[k], A_{i+k+2}) = 0 $ for all $ i \ge 0, k \ge 1 $ and in particular $ \Hom(B_{i-1}[1], A_{i+2}) = 0$.

Now, we have the two distinguished triangles
\begin{equation*}
\begin{CD}
B_{i-1} @>g_i>> A_i @>h_i>> B_i \\
@VVV @VVf_{i+1}V \\
0 @>>> A_{i+1} @>>id> A_{i+1}
\end{CD}
\end{equation*}
The left square commutes since by induction $ f_{i+1} \circ g_i = 0 $.  Thus by the axioms of triangulated categories, we may choose a morphism $g_{i+1} $ making this into a commutative diagram of distinguished triangles.  This will mean that $ f_{i+1} = g_{i+1} \circ h_i $.  However, we might not have $ f_{i+2} \circ g_{i+1} = 0 $. 

To prove this last relation, consider the long exact sequence which comes from applying $ \Hom(\cdot, A_{i+2}) $ to the distinguished triangle (\ref{eq:distB}).  We see that
\begin{equation*}
\Hom(B_{i-1}[1], A_{i+2}) \rightarrow \Hom(B_i, A_{i+2}) \rightarrow \Hom(A_i, A_{i+2}).
\end{equation*}
Now, $ f_{i+2} \circ g_{i+1} $ is an element in the middle group and is sent to $ f_{i+2} \circ g_{i+1} \circ h_i = f_{i+2} \circ f_{i+1} = 0 $ under the right arrow.  Since, $ \Hom(B_{i-1}[1], A_{i+1}) = 0 $, this forces $ f_{i+2} \circ g_{i+1} = 0 $ as desired.  Thus, we have produced $ g_{i+1} $ with the desired properties and hence a convolution exists by this construction.
\end{proof}

The behaviour of convolutions under functors is quite easy to understand.  Let $ \mathcal{D}_1, \mathcal{D}_2 $ be triangulated categories and let $ \F : \mathcal{D}_1 \rightarrow \mathcal{D}_2 $ be a triangulated functor.  Let $ (A_\bullet, f_\bullet) $ be a complex in $\mathcal{D}_1$.

\begin{Proposition} \label{th:pushcone}
If $ B $ is a convolution of $ (A_\bullet, f_\bullet ) $ then $\F(B) $ is a convolution of $ (\F(A_\bullet), \F(f_\bullet)) $.
\end{Proposition}

\begin{proof}
Since $ B $ is a convolution, it comes with a Postnikov system $ B_i, g_i, h_i $ as above.  Since $ \F $ is a functor and takes distinguished triangles to distinguished triangles, it follows that $ \F(B_i), \F(g_i), \F(h_i) $ is a Postnikov system, and so we see that $ \F(B) $ is a convolution of $ (\F(A_\bullet), \F(f_\bullet)) $.
\end{proof}

Two examples of unique convolutions worth highlighting. The first example occurs when $\mathcal{D}$ is the derived category of an abelian category $\mathcal{A}$. If every term in the complex $A_0 \rightarrow A_1 \rightarrow \dots \rightarrow A_n$ is an element of $\mathcal{A} \subset \mathcal{D}$ then the complex $ A_\bullet $ in $\mathcal{D} $ is the convolution. The second example is when $n=1$, in which case the usual cone of $A_0 \rightarrow A_1$ is the unique convolution of this 2-term complex. 

The first example can be generalized in the following way (though we may lose uniqueness).  Let $\mathcal{A} $ be an abelian category, and let $ A_0 \xrightarrow{\tilde{f}_1} A_1 \xrightarrow{\tilde{f}_2} \cdots \xrightarrow{\tilde{f}_n} A_n $ be a sequence of objects and morphisms in $ Kom(\mathcal{A}) $ (the category of complexes of objects in $ \mathcal{A} $) with $ \tilde{f}_{i+1} \circ \tilde{f}_i = 0 $.  We call such a collection a \textbf{complex} in $Kom(\mathcal{A}) $.  Since $Kom(\mathcal{A}) $ is not triangulated, this is not covered by our previous definition.  On the other hand, let $ f_i $ denote the homotopy class of $ \tilde{f}_i$.  Then $ (A_\bullet, f_\bullet) $ is a complex in $ K(\mathcal{A}) $ (the homotopy category of complexes of objects in $\mathcal{A}$).  

Now, consider the following construction in $ Kom(\mathcal{A}) $.  For each $ i $, we let $ B_i $ be the complex with $ B_i^k = A_0^{i+k} \oplus \cdots \oplus A_i^k $.  Equip $ B_i $ with a differential by defining matrix coefficients $ d_{j' j}^k : A_j^{i - j+k} \rightarrow A_{j'}^{i- j'+k+1} $ by
\begin{equation*}
d_{j' j}^k = \begin{cases}
(-1)^{i-j} d_j \quad \text{ if } j = j' \\
\tilde{f}_{j+1} \quad \text{ if } j' = j+1 \\
0 \quad \text{ otherwise}
\end{cases}
\end{equation*}
where $ d_j $ denotes the differential in the complex $ A_j $.

We also define maps $ \tilde{g}_i : B_{i-1} \rightarrow A_i $ by projecting $ B_{i-1} $ onto $ A_{i-1} $ and then using $ \tilde{f}_i $.  The projection is not a map of complexes, but its composition with $ \tilde{f}_i $ is (essentially because $\tilde{f}_i \circ \tilde{f}_{i-1} = 0 $).  Also the map $ \tilde{h}_i : A_i \rightarrow B_i $ is defined in the obvious way.

In this case, we call $ B = B_n $, \textbf{the convolution} of $ (A_\bullet, \tilde{f}_\bullet)$ and we write $ \Con(A_\bullet, \tilde{f}_\bullet) $ for $ B $. 

Let $ g_i, h_i $ be the homotopy classes of $ \tilde{g}_i, \tilde{h}_i$.  Then $ (B_\bullet, g_\bullet, h_\bullet) $ is a Postnikov system and so $ B $ is a convolution of $ (A_\bullet, f_\bullet) $ in $K(\mathcal{A}) $. In fact, $\Con(A_\bullet, \tilde{f}_\bullet)$ is nothing more than the total complex of the double complex $A_i^j$ with $B_i$ the total complex of the partial double complex $A_0 \xrightarrow{\tilde{f}_0} A_1 \xrightarrow{\tilde{f}_1} \cdots \xrightarrow{\tilde{f}_{i-1}} A_i$.

\subsection{Cones from tensor products} \label{se:coneandtensor}
Now, we will introduce our main example of complexes and their convolutions.  Let $ \mathcal{A} $ be an abelian category with a tensor product $ \otimes $.  Recall that $K(\mathcal{A})$ is a triangulated category where $ \otimes $ extends to a bifunctor $ K(\mathcal{A}) \times K(\mathcal{A}) \rightarrow K(\mathcal{A}) $.

Let $\{C_j \xrightarrow{r_j} D_j \}_{1 \le j \le n} $ be a collection of objects and morphisms of $ K(\mathcal{A}) $.  For each $ j $, let $ E_j $ be an object of $ K(\mathcal{A}) $ such that $ E_j $ is a cone of $ r_j $.  

Given this setup, for $ \delta \subset \{1, \dots, n\} $, we define the object $ A_\delta := X^\delta_1 \otimes \dots \otimes X^\delta_n $, where $ X^\delta_j = C_j $ if $ j \notin \delta $ and $ D_j $ if $ j \in \delta $.  Then, we define $ A_i := \oplus_{\delta, |\delta| = i} A_\delta $.  In other words $ A_i $ is the direct sum of all tensor products with $ n-i $ $C_j$s and $ i$ $D_j$s.  We define $ f_i : A_{i-1} \rightarrow A_i $ to be the morphism whose matrix entries $ f_{\delta' \delta} : A_\delta \rightarrow A_{\delta'}$ are 
\begin{equation} \label{eq:defmatent}
f_{\delta' \delta} = \begin{cases}
(-1)^{\# \{ k \in \delta : k \le j \}} 1 \otimes r_j \otimes 1 \quad \text{if } \delta' = \delta \sqcup j  \\
0 \quad \text{otherwise}
\end{cases}
\end{equation}
Here and later, when we write $ \delta' = \delta \sqcup j $, we mean that $ \delta' = \delta \cup \{j\} $ and that $ j \notin \delta $.

\begin{Proposition} \label{th:coneandtensor}
$ (A_\bullet, f_\bullet) $ is a complex and $ E_1 \otimes \cdots \otimes E_n $ is isomorphic to some convolution of $(A_\bullet, f_\bullet)$.
\end{Proposition}

\begin{proof}
To show that $ (A_\bullet, f_\bullet) $ is complex, we must show that $ f_{i+1} \circ f_i = 0 $.  For this it suffices to consider $ \delta $ with $ |\delta| = i-1 $.  Let $ j, k \notin \delta $ and consider the matrix element of $ f_{i+1} \circ f_i $ which is a map $ A_\delta \rightarrow A_{\delta \sqcup j,k} $.  By the definition, this matrix element is the sum 
\begin{equation*}
f_{\delta \sqcup j,k \, \delta \sqcup j} \circ f_{\delta \sqcup j \, \delta} + f_{\delta \sqcup j,k  \, \delta \sqcup k} \circ f_{\delta \sqcup k \, \delta}.
\end{equation*}
Examining the definition of $ f_{\delta \delta'} $ in (\ref{eq:defmatent}) we see that these two compositions are the same, except with opposite sign.  Thus $ f_{i+1} \circ f_i = 0 $ and so $ (A_\bullet, f_\bullet) $ is a complex.

Now it remains to show that $ E_1 \otimes \cdots \otimes E_n $ is isomorphic to a convolution of $ (A_\bullet, f_\bullet)$.  To do this, let us first choose a map of complexes $ \tilde{r}_j : C_j \rightarrow D_j $ which is a representative for $ r_i $.  Then we have a complex $ \Cone(\tilde{r}_j) $ with $ \Cone(\tilde{r}_j)^k = C_j^{k+1} \oplus D_j^k $ and with differential given by the matrix
\begin{equation*}
\begin{bmatrix}
-d_{C_j} & 0 \\ 
\tilde{r}_j & d_{D_j}  \\
\end{bmatrix}
\end{equation*}
where $ d_{C_j}, d_{D_j} $ denote the appropriate terms of the differentials for $ C_j, D_j $.  

By the definition of the triangulated structure in $ K(\mathcal{A}) $, we have $ \Cone(\tilde{r}_j) \cong E_j $.  Hence it suffices to show that $ \Cone(\tilde{r}_1) \otimes \cdots \otimes \Cone(\tilde{r}_n) $ is a convolution of $ (A_\bullet, f_\bullet) $.  Note that by definition 
\begin{equation*}
\big(\Cone(\tilde{r}_1) \otimes \cdots \otimes \Cone(\tilde{r}_n)\big)^k = \oplus_\delta A_\delta^{k+n-|\delta|}
\end{equation*}
with differential $ d $ having matrix coefficients $ d^k_{\delta' \delta} : A_\delta^{k+n-|\delta|} \rightarrow A_{\delta'}^{k+1+n-|\delta'|} $ defined by
\begin{equation*}
d^k_{\delta' \delta} = \begin{cases} 
(-1)^{n-|\delta|} d_{A_\delta} \quad \text{if $\delta' = \delta$} \\
(-1)^{\# \{l \in \delta : l \le j \}} 1 \otimes \tilde{r}_j \otimes 1 = f_{\delta' \delta} \quad \text{if $\delta' = \delta \sqcup j $} \\
0 \quad \text{otherwise}.
\end{cases}
\end{equation*}

The choice of $ \tilde{r}_i $ as a representative for $ r_i $ gives us natural choices for representatives $ \tilde{f}_i $ for the maps $ f_i $.  Note that $ \tilde{f}_{i+1} \circ \tilde{f}_i = 0 $.  Hence we have a complex $ (A_\bullet, \tilde{f}_\bullet) $, in $Kom(\mathcal{A}) $.  It is easily seen that $ \Con(A_\bullet, \tilde{f}_\bullet) = \Cone(\tilde{r}_1) \otimes \cdots \otimes \Cone(\tilde{r}_n) $ and hence $\Cone(\tilde{r}_1) \otimes \cdots \otimes \Cone(\tilde{r}_n)$ is a convolution of $ (A_\bullet, f_\bullet) $ as desired (see the remarks at the end of section \ref{se:cones}).
\end{proof}

\subsection{Cones from tangles}
Now we are ready to apply this theory of convolutions to the kernels coming from tangle diagrams.

Let $ T $ be an $(n,m) $ tangle diagram.  Recall from above that $ \Q(T) \in D(Y_n \times Y_m) $ is the kernel that associated to $T$ by the unoriented theory. Let $ S $ denote the set of crossings in $ T $.  Note that $ S $ acquires a total order coming from the order (in terms of height) in which the crossings appear in $ T $.  For any subset $ \delta $ of $ S $, let $ T_\delta $ denote the crossingless tangle created by resolving all the crossings of $ T $ according to the subset $ \delta $.  That is, a crossing $ s \in S $ is resolved into the ``0-resolution'' when $ s \notin \delta $ and into the ``1-resolution'' when $ s \in \delta $.  Let $ R_\delta := \Q(T_\delta)[|\delta|-|S|]\{-|\delta|\} $.  

If $\delta' = \delta \sqcup s $, then there is a morphism 
\begin{equation*}
(-1)^{\# \{ k \in \delta : k \le s \}} \alpha_U^s[|\delta'|-|S|]\{-|\delta|\} :  R_\delta \rightarrow R_{\delta \cup \{s\}}
\end{equation*}
where $ U $ is the one crossing tangle made by resolving all crossings except $ s $ according to $ \delta $ and leaving the crossing $ s $ (recall that $\alpha_U^s: \Q(T_\delta)[-1] \rightarrow \Q(T_{\delta'}) \{-1\}$). 

Now we let $ R_i = \oplus_{\delta, |\delta| = i} R_\delta $.  We get a complex $ (R_\bullet, h_\bullet) $ by using the above maps as the matrix entries.

\begin{Proposition} \label{th:Qiscone}
$\Q(T) $ is isomorphic to a convolution of $(R_\bullet, h_\bullet) $.  
\end{Proposition}

\begin{proof}
By definition, we have $ \Q(T) = \P_1 * \cdots * \P_k $ where each $ \P_j \in D(Y_{n_{j-1}} \times Y_{n_j}) $ is one of $ \sF_{n_j+2}^i, \sG_{n_j}^i, (\sT_{\sG_{n_j}^i})_L, \sT_{\sG_{n_j}^i}\{-1\} $ (here $ n_0 =n, n_k = m$).  

It follows from basic properties of the convolution product $ * $ that
\begin{equation*}
\Q(T) \cong {\pi_{0\, k}}_* \big( \pi_{0\, 1}^* \P_1 \otimes \cdots \otimes \pi_{k-1\, k}^* \P_k \big)
\end{equation*}
where the tensor products take place in $ D(Y_{n_0} \times \cdots \times Y_{n_k}) $ and where $ \pi_{i\, j} $ denotes the projection onto $ D(Y_{n_{i}} \times Y_{n_j}) $.  

Now, define objects and morphisms $ C_j \xrightarrow{r_j} D_j $ in $ D(Y_{n_0} \times \cdots \times Y_{n_k}) $ by the following rule
\begin{enumerate}
\item if $\P_j = \sF_n^i$, then $ C_j = 0, D_j = \pi_{j-1\, j}^* \sF_n^i $,
\item if $\P_j = \sG_n^i $, then $ C_j = 0, D_j = \pi_{j-1\, j}^* \sG_n^i $,
\item if $ \P_j = (\sT_{\sG_n^i})_L $, then $ C_j = \pi_{j-1\, j}^* \O_\Delta[-1],\ D_j = \pi_{j-1\, j}^* \sG_n^i * \sF_n^i \{-1\} $ and $ r_j $ comes from the natural map $ \O_\Delta \rightarrow \sG_n^i * (\sG_n^i)_L $ and
\item if $ \P_j = \sT_{\sG_n^i}\{-1\}$, then $ C_j = \pi_{j-1\, j}^* \sG_n^i * \sF_n^i [-1],\ D_j = \pi_{j-1\, j}^* \O_\Delta\{-1\} $ and $ r_j $ comes from the natural map $ \sG_n^i * (\sG_n^i)_R \rightarrow \O_\Delta $.
\end{enumerate}
Note that in each case, $ C_j \xrightarrow{r_j} D_j \rightarrow \pi_{j-1\, j}^* \P_j $ is a distinguished triangle.  

Now, replace each $ C_j, D_j , \pi_{j-1\, j}^* \P_j$ by quasi-isomorphic complexes $ C'_j, D'_j, E'_j $ of locally free sheaves and replace $ r_j $ by a map $r'_j $ in the homotopy category of locally free sheaves, such that we still have a distinguished triangle $ C'_j \xrightarrow{r'_j} D'_j \rightarrow E'_j $ and such that $ r'_j $ goes over to $ r_j $ under the quasi-isomorphisms (ie there is a commuting diagram in the derived category).  Since locally free sheaves are acyclic for $ \otimes $, the tensor product $ \pi_{0\,1}^*\P_1 \otimes \cdots \otimes \pi_{k-1\, k}^* \P_k  $ in the derived category is isomorphic to $ E'_1 \otimes \cdots \otimes E'_k $, where in the latter expression the tensor product is computed in the homotopy category of locally free sheaves.

Now, we are in the situation of Proposition \ref{th:coneandtensor}.  Adopting the notation of that proposition, for $ \delta \subset \{1, \dots, k \}$, define $ X_j^\delta $ to be $ C'_j $ if $ j \notin \delta $ or $ D'_j $ if $ j \in \delta $.  Let $ A_\delta = X_1^\delta \otimes \cdots \otimes X_k^\delta$, and let $ A_i = \oplus_{\delta, |\delta| = i} A_\delta $.  Let $ f_i : A_{i-1} \rightarrow A_i $ denote the maps defined in the beginning of section \ref{se:coneandtensor}.  

Proposition \ref{th:coneandtensor} shows that $ E'_1 \otimes \cdots \otimes E'_k $ is isomorphic in the homotopy category of locally free sheaves to a convolution of $ (A_\bullet, f_\bullet) $.  A convolution in the homotopy category of locally free sheaves is also a convolution in the derived category of coherent sheaves.  So now we pass back to the derived category and conclude that $\pi_{01}^*\P_1 \otimes \cdots \otimes \pi_{k-1 k}^* \P_k  $ is isomorphic to a convolution of $ (A_\bullet, f_\bullet) $.  Thus by Proposition \ref{th:pushcone}, $ {\pi_{0 k}}_* (\pi_{01}^*\P_1 \otimes \cdots \otimes \pi_{k-1 k}^* \P_k) $ is isomorphic to a convolution of $ ({\pi_{0 k}}_* A_\bullet, {\pi_{0 k}}_* f_\bullet) $.  Since $ \Q(T) $ is isomorphic to $ {\pi_{0 k}}_* (\pi_{01}^*\P_1 \otimes \cdots \otimes \pi_{k-1 k}^* \P_k)  $, we see that $ \Q(T) $ is isomorphic to a convolution of $ ({\pi_{0 k}}_* A_\bullet, {\pi_{0 k}}_* f_\bullet) $.

To complete the proof we need to relate $({\pi_{0 k}}_* A_\bullet, {\pi_{0 k}}_* f_\bullet) $ and $(R_\bullet, h_\bullet)$. Note that $ \{1, \dots, k\} $ is in bijection with the set of caps, cups, and crossings of $ T $ so if $S$ denotes the set of crossings we have $ S \subset \{1, \dots, k\} $.  Also, since $ C_j = 0 $ for $ j \notin S $, we see that $ A_\delta = 0 $ unless $\{1, \dots, k\} \smallsetminus S = S^c \subset \delta $.  Moreover, we see that if $ S^c \subset \delta$, then $ {\pi_{0\,k}}_* A_\delta $ is $\Q(T_{\delta \cap S})$ with some shift. From the definition of $\P_j$ above we see the shift is $[-1]$ for each 0-resolution and $\{-1\}$ for each 1-resolution. Since there are $|\delta \cap S|$ 0-resolutions and $|S| - |\delta \cap S|$ 1-resolutions we get ${\pi_{0\,k}}_* A_\delta \cong R_{\delta \cap S}$. 

Also if $ S^c \subset \delta $ and $ s \notin \delta$, then $ {\pi_{0\, k}}_* 1 \otimes r_s \otimes 1 $ which is a map $ {\pi_{0\, k}}_* A_\delta \rightarrow {\pi_{0\,k}}_\ast A_{\delta \sqcup s} $, agrees with the map $ \alpha_U^s $ (up to the appropriate shift) defined above. Hence we have that if $ i \ge k - |S| $, then 
\begin{equation*}
{\pi_{0\, k}}_* A_i = \bigoplus_{\delta \subset \{1, \dots, k\}, |\delta| = i} {\pi_{0\, k}}_* A_\delta = \bigoplus_{S^c \subset \delta, |\delta|=i} R_{\delta \cap S} = \bigoplus_{\delta' \subset S, |\delta'| = i-(k-|S|)} R_{\delta'} = R_{i-(k-|S|)}
\end{equation*}
while $ {\pi_{0\, k}}_* A_i = 0 $ if $ i < k-|S| $. Also, under the above isomorphism, the map $ {\pi_{0\,k}}_* (f_i: A_{i-1} \rightarrow A_i) $ goes over to the map $ h_i: R_{i-1-(k-|S|)} \rightarrow R_{i-(k-|S|)} $ defined above.

Since $ \Q(T) $ is isomorphic to the convolution of $ ({\pi_{0 k}}_* A_\bullet, {\pi_{0 k}}_* f_\bullet) $, we see that $ \Q(T) $ is the convolution of a sequence which is all 0 at the beginning and then is isomorphic to the sequence $ (R_\bullet, h_\bullet) $.  Hence $ \Q(T) $ is isomorphic to a convolution of $ (R_\bullet, h_\bullet) $.
\end{proof}

\subsection{Definition of Khovanov homology}
We will now recall the definition of Khovanov homology.  For our purposes, we will work with a ``sheared'' version of the Khovanov complex which we denote $ \sh $.  

We start with the vector space $ V = \C[1]\{-1\} \oplus \C[-1]\{1\} = \langle v^-, v^+ \rangle$. There are maps $ m : V \otimes V[-1] \rightarrow V \{-1\} $ and $ \Delta : V [-1] \rightarrow V \otimes V \{-1\} $ given by 
\begin{gather*}
m[1]: V^{\otimes 2} \rightarrow V[1]\{-1\} \\
v^- \otimes v^- \mapsto v^-[1]\{-1\} ,\ v^- \otimes v^+ \mapsto v^+[1]\{-1\} ,\ v^+ \otimes v^- \mapsto v^+[1]\{-1\} ,\ v^+ \otimes v^+ \mapsto 0 
\end{gather*}
and 
\begin{gather*}
\Delta[1]: V \rightarrow V^{\otimes 2}[1]\{-1\} \\
v^- \mapsto (v^- \otimes v^+ + v^+ \otimes v^-)[1]\{-1\} ,\ v^+ \mapsto v^+ \otimes v^+[1]\{-1\}.
\end{gather*}
For a crossingless link $ K $, define $\sh(K) := \otimes_{\text{circles of } K} V$, a tensor product of copies of $ V $ indexed by the circles of $ K $. 

Now $ K $ be an arbitrary link, let $ S $ be its set of crossings, and let $ n = |S| $.  As above, for $ \delta \subset S $, let $ K_\delta $ denote the result of resolving all the crossings according to $ \delta $.  Suppose that $ \delta' = \delta \sqcup s $.  Then switching the resolution of the crossing $ s $ either connects two circles of $ K_\delta $ into a circle of $ K_{\delta'} $ or it breaks a circle of $ K_\delta $ into two circles of $K_{\delta'}$.  Thus there is a natural correspondence between the circles of $ K_\delta $ and those of $ K_{\delta'} $ --- this correspondence is a bijection except for identifying a circle of $ K_\delta $ with two of $ K_{\delta'} $ or vice versa.  Using this correspondence and the maps $ m, \Delta $ above, we get a map $ \tilde{g}_{\delta' \delta} : \sh(K_\delta)[-1] \rightarrow \sh(K_{\delta'})\{-1\}$ in the category of complexes of graded vector spaces.

We define $ \sh(K)$ to be the complex with $ \sh(K)^k = \oplus_{\delta} \sh(K_\delta)^k \{-|\delta|\} $.  We define a differential on $ \sh(K) $ by giving it matrix entries $ d_{\delta' \delta} : \sh(K_\delta)^k  \{-|\delta|\} \rightarrow \sh(K_{\delta'})^{k+1}  \{-|\delta'|\} $ by
\begin{equation*}
d_{\delta' \delta} = \begin{cases}
(-1)^{\# \{ i \in \delta : i \le s \}} \tilde{g}_{\delta' \delta} \{-|\delta|\} \quad  \text{ if } \delta' = \delta \sqcup s \\
0 \quad \text{ otherwise} 
\end{cases}
\end{equation*}

Note that $ \sh(K) $ is the convolution of a sequence $(M_i(K), \tilde{g}_i) $ in the category of complexes of graded vector spaces, with $ M_i = \oplus_{\delta, |\delta| = i} \sh(K_\delta)[|\delta|-n]\{-|\delta|\} $ and $ \tilde{g}_i $ defined using shifts of the $ \tilde{g}_{\delta' \delta} $ as its matrix coefficients.

The definition of $ M(K) $ is just a sheared version of the definition of the ``Khovanov bracket'' $ [[K]]$ (notation from Bar-Natan \cite{BN}).  Hence it follows that $ H^{i,j}(M(K)) = H^{i+j,j}[[K]]$. 

\subsection{Proof of Theorem \ref{th:Khov}}

Let $ \triv^m $ denote the trivial link projection with $ m $ circles arranged into a vertical strip.  By Corollary \ref{cor:circle}, we have an isomorphism $ \Q(\triv^m) \cong V^{\otimes m} $.  Now, let $ A $ be an invertible combinatorial cobordism between $ \triv^m $ and itself.  So $ A $ is an invertible 2-morphism in $ \mathcal{C}' $ between $ \triv^m $ and itself.  The condition that it is invertible is equivalent to the condition that it is made by a composition of movie moves not including births, deaths and saddles.  So $ A $ consists of tubes joining the $ m $ circles with themselves and so gives rise to a permutation of the circles in $ \triv^m $ and hence an element $ \sigma \in S_m $.  The following result describes $ \Psi(A) $.

\begin{Lemma} \label{th:movingcircles}
Up to a non-zero scalar the following diagram commutes
\begin{equation*}
\begin{CD}
\Psi(\triv^m) @>>> V^{\otimes m} \\
@V\Psi(A)VV @VV\sigma V \\
\Psi(\triv^m) @>>> V^{\otimes m}
\end{CD}
\end{equation*}
\end{Lemma}

\begin{proof}
The set of such cobordisms forms a group called the ``loop braid group''.  Baez-Crans-Wise \cite{BCW} give a presentation of this group.  There are two types of generators, those which corresponds to exchanging two circles and those which correspond to threading one circle through another (see \cite{BCW} section 2 for nice pictures of these generators).  So it suffices to show that the above diagram commutes for these generators.  Each of these generators only involve two circles, so it suffices to consider $ m = 2 $.  The proof for each type of generator is similar, so we give it when $ A $ is ``threading'' generator (denoted $ \sigma_{ij} $ in \cite{BCW}).  In this case $ \sigma = id $ and so our goal is to prove that $ \Psi(A) $ is given by multiplication by a scalar.

We compose the birth cobordism $ b $ from $ \triv $ to $ \triv^2 $ with $ A $.  The composition $ A \circ b $ is equal (in $\mathcal{C'} $) to $ b $ (because we can move the cap along the tube).  Since $ A \circ b = b $, applying  $ \Psi $ gives $ \Psi(A \circ b) = \Psi(b) $ and so by Theorem \ref{thm:cobordism} we get the commutative diagram (up to scalar)
\begin{equation*}
\begin{CD}
V \otimes \C @>\Psi(b)>> V \otimes V \\
@| @VV\Psi(A)V \\
V \otimes \C @>>\Psi(b)> V \otimes V
\end{CD}
\end{equation*}
Let $ \lambda $ be the non-zero scalar involved.  By definition, $ \Psi(b)(v \otimes 1) = v \otimes v_- $.  Thus we see that $ \Psi(A)(v \otimes v_-) = \lambda v \otimes v_- $.

Now we apply the same argument again but using the birth of the other circle.  This shows that there exists a non-zero scalar $ \mu $ such that $ \Psi(A)(v_- \otimes v) = \mu v_- \otimes v $.  Because of the overlap of these two cases, we see that $ \lambda = \mu $.  

This almost determines $ \Psi(A)$ , except for $ \Psi(A)(v_+ \otimes v_+)$.  For this, we compose with the death cobordism $ d $.  Again $ d \circ A = d $ and so we get the commutative diagram (up to scalar)
\begin{equation*}
\begin{CD}
V \otimes V @>\Psi(d)>> V \otimes \C \\
@V\Psi(A)VV @| \\
V \otimes V @>>\Psi(d)> V \otimes \C
\end{CD}
\end{equation*}
Since $ \Psi(d)(v \otimes v_-) = 0 $ and $ \Psi(d) (v \otimes v_+) = v $, we see that there exists a non-zero scalar $ \mu $ and scalar $ \nu $, such that $ \Psi(A)(v_+ \otimes v_+) = \mu v_+ \otimes v_+ $ and $ \Psi(A)(v_- \otimes v_+) = \mu v_- \otimes v_+ + \nu v_+ \otimes v_- $.  By comparison with the above, we see that $ \mu = \lambda $ and $ \nu = 0 $.  Thus the result follows.
\end{proof}

For each $ \delta$, choose an invertible combinatorial cobordism of $ K_\delta $ with any fixed vertical strip of circles $S$ (with the same number of components as $ K_\delta $).  By Theorem \ref{thm:cobordism}, this cobordism gives rise to an isomorphism $ \Q(K_\delta) \rightarrow \Q(S) \cong V^{\otimes m} $.  By Lemma \ref{th:movingcircles}, this isomorphism gives rise to an isomorphism $ \Q(K_\delta) \cong \otimes_{\text{circles of $K_\delta$}} V = \sh(K_\delta)$, which (up to scalar) is independent of the choice of $ S $ and the cobordism between $ K_\delta $ and $ S $.

\begin{Lemma} \label{th:localKh}
Choose $ \delta, \delta' $ where $ \delta' = \delta \sqcup i $ for some $ i $.

Under the isomorphism $\Q(K_\delta) \cong \sh(K_\delta) $, the maps $ \alpha_U^s $ and $ g_{\delta' \delta} $ coincide up to a non-zero scalar, where, as before,  $ U $ denotes the one crossing link with a crossing at $ s $ and $g_{\delta' \delta} $ denotes the homotopy class of $ \tilde{g}_{\delta' \delta}$.
\end{Lemma}
\begin{proof}
We begin with the case that $ K $ is a one crossing link which consists of a vertical strip of circles and one figure 8 and we assume that the cobordisms of $ K_0 $ and $ K_1 $ with vertical strips are the identity cobordisms. There are two cases depending on whether the crossing is of type 1 or 2.

Suppose the crossing is of type 2. We can assume $K$ is just the figure 8 with $K_0$ and $K_1$ consisting of two and one circle respectively. For simplicity denote $\sF_2^1 = \sF$ and $\sG_2^1 = \sG$. Consider the induced map $\Q(K_0) = V^{\otimes 2}[-1] \xrightarrow{\alpha_2} V\{-1\} = \Q(K_1)$. More precisely it is $\sF \ast \left( \sG \ast \sF [-1] \right) \ast \sG \rightarrow \sF \ast \sG \{-1\}$ induced by the adjoint map $(\sG \ast \sG_R) \rightarrow \O_\Delta$ where we recall $\sF[-1] = \sG_R\{-1\}$. We need to show $\alpha_2$ is the same as $m$ (up to a non-zero multiple).

The identification $\sF \ast \sG = V$ is made using the maps from the standard exact triangle 
\begin{equation} \label{eq:Vses}
\langle v^- \rangle = \C[1]\{-1\} \rightarrow \sG_R \ast \sG [1]\{-1\} = \sF \ast \sG = \sG_L \ast \sG [-1]\{1\} \rightarrow \C[-1]\{1\} = \langle v^+ \rangle.
\end{equation}
If we apply $\ast \sF \ast \sG[-1]$ to the left side of (\ref{eq:Vses}) and compose with $\alpha_2$ we get
$$v^- \otimes \sF \ast \sG[-1] \rightarrow \sG_R \ast \sG [1]\{-1\} \ast \sF \ast \sG[-1] = \sF \ast \sG \ast \sF \ast \sG [-1] \xrightarrow{\alpha_2} \sF \ast \sG \{-1\}.$$ 
The composition of these maps is the identity since it's induced by the composition 
$$\sG_R \rightarrow (\sG_R \ast \sG) \ast \sG_R = \sG_R \ast (\sG \ast \sG_R) \rightarrow \sG_R$$
which is the identity by general theory of adjoints. This means that $\alpha_2(v^- \otimes v^-[-1]) = v^-\{-1\}$ and $\alpha_2(v^- \otimes v^+[-1]) = v^+\{-1\}$. Similarly, one can apply $\sF \ast \sG[-1] \ast$ to the left side of (\ref{eq:Vses}) and the same argument shows $\alpha_2(v^- \otimes v^-[-1]) = v^-\{-1\}$ and $\alpha_2(v^+ \otimes v^-[-1]) = v^+\{-1\}$. Finally, $\alpha_2(v^+ \otimes v^+[-1]) = 0$ since $v^+ \otimes v^+[-1]$ lies in degree $(1,-2)$ while $V\{-1\}$ is non-zero only in degrees $(1,0)$ and $(-1,2)$. Thus $\alpha_2$ coincides with $m$ in this case.

Similarly, for a type 1 crossing we get a map $\Q(K_0) = V[-1] \xrightarrow{\alpha_1} V^{\otimes 2}\{-1\} = \Q(K_1)$ which we must show to coincide with the map $\Delta$. To do this we imitate the argument above. Applying $\sF \ast \sG \{-1\} \ast$ to the right side of (\ref{eq:Vses}) and precomposing with $\alpha_1$ we get
$$\sF \ast \sG[-1] \xrightarrow{\alpha_1} \sF \ast \sG \ast \sF \ast \sG \{-1\} \xrightarrow{p} \sF \ast \sG [-1]$$
where the composition is the identity. Now degree considerations imply $\alpha_1(v^+) = a (v^+ \otimes v^+)$ for some $a \in \C$ while $p(v^+ \otimes v^+) = v^+$. Since the composition is the identity this means $a=1$. Similarly, degree considerations imply $\alpha_1(v^-) = a_1 (v^+ \otimes v^-) + a_2 (v^- \otimes v^+)$ for some $a_1,a_2 \in \C$. Since $p(v^- \otimes v^+)=0$ and $p(v^+ \otimes v^-) = v^+$ we find $a_1=1$. On the other hand, applying $\ast \sF \ast \sG \{-1\}$ to the right side of (\ref{eq:Vses}) also gives $a_2=1$. Hence we conclude that $\alpha_1 = \Delta$. 

Now, we proceed to the case of general $ K $. We choose a combinatorial cobordism of $ U $ with $ S_U $, which is a vertical strip with circles and one figure 8 (as above).  This cobordism induces a cobordism of $ K_\delta $ with  $ S_U^0 $ and $ K_{\delta'} $ with  $S_U^1 $ where $ S_U^0 $ and $ S_U^1 $ are vertical strips of $m$ and $ m'$ circles.  Moreover, we have a commutative diagram of 2-morphisms in $ \mathcal{C'} $ 
\begin{equation} \label{eq:toapplyQ}
\begin{CD}
K_\delta @>>> S_U^0  \\
@VVV @VVV \\
K_{\delta'}@>>> S_U^1 \\
\end{CD}
\end{equation}
where the horizontal arrows denote the above cobordisms and the vertical arrows denote the saddle maps.  We can then apply $ \Q $ and get a diagram which commutes up to scalar by Theorem \ref{thm:cobordism}. 

We now consider the following diagram.
\begin{equation*}
\begin{CD}
V^{\otimes m} @<<< \Q(S_U^0) @>>> \Q(K_\delta) @>>> \Q(S) @>>> V^{\otimes m} = \sh(S_\delta) \\
@V{m, \Delta}VV @V{\alpha_{S_U}}VV @V\alpha_UVV @. @V{g_{\delta \delta'}}VV \\
V^{\otimes m'} @<<< \Q(S_U^1) @>>> \Q(K_{\delta'}) @>>> \Q(S') @>>> V^{\otimes m'} = \sh(S_{\delta'})
\end{CD}
\end{equation*}
Here the leftmost vertical map is either $ m $ or $ \Delta $ depending on the type of the crossing.  The leftmost square commutes (up to scalar) by the first half of this proof.  The next square commutes (up to scalar) by applying $\Q $ to (\ref{eq:toapplyQ}).  

Now by Lemma \ref{th:movingcircles} the maps all the way across the top and across the bottom are just permutations (up to scalar).  Hence the right rectangle commutes up to scalar, which is what we needed to show.
\end{proof}

\begin{Remark} Khovanov describes the vector space $ V $ as the cohomology of $ \p $ and the map $ m $ as the cup product (see for example \cite{Kh5}).  This can be seen very naturally from our construction.

In Remark \ref{re:Dolb}, we explained that there were natural isomorphisms $$V [-1]\{1\} \cong \sF \ast \sG [-1]\{1\} \cong \Ext_{Y_2}(\O_{Y_1}, \O_{Y_1}) = \oplus_{i,j} \mathrm{H}^{i,j}(\p)[-i-j]\{2i\}.$$  The isomorphism $ \Ext_{Y_2}(\O_{Y_1}, \O_{Y_1}) = \oplus_{i,j} \mathrm{H}^{i,j}(\p)[-i-j]\{2i\}$ is an isomorphism of algebras where on the LHS we take the Ext multiplication and on the RHS we take the cup product.  Moreover, the natural adjunction map $\sF \ast (\sG \ast \sF[-1]) \ast \sG \rightarrow \sF \ast \sG \{-1\}$ induces the Ext multiplication by general principles.  Hence under the identification of $\sF \ast \sG  $ with the cohomology of $ \p $, the map $ \alpha $ goes over to cup product.  
\end{Remark}

\begin{Theorem} \label{th:Q=M}
If $ K $ is a link projection, then $ \Q(K) $ is isomorphic to $\sh(K)$.
\end{Theorem}

\begin{proof}
By Proposition \ref{th:Qiscone}, we know that $ \Q(K) $ is isomorphic to a convolution of $ (R_\bullet, h_\bullet) $.  Hence it suffices to show that $ \sh(K) $ is also a convolution of $ (R_\bullet, h_\bullet) $ and that $ (R_\bullet, h_\bullet) $ has a unique convolution.

For the first part, note that $ R_\delta \cong \sh(K_\delta)[|\delta|-n]\{-|\delta|\} $ and  by Lemma \ref{th:localKh} the matrix coefficients $h_{\delta \delta'}$ that define the $ h_j: R_{j-1} \rightarrow R_j $ are precisely the maps $ g_{\delta \delta'} $ up to a non-zero scalars $ \lambda_{\delta \delta'} $. Suppose we consider a square
\begin{equation*}
\begin{CD}
R_\delta @>>> R_{\delta \sqcup t} \\
@VVV @VVV \\
R_{\delta \sqcup s} @>>> R_{\delta \sqcup s,t}
\end{CD}
\end{equation*}
where the arrows are the various $ h $ maps.  Since these arrows are essentially adjunction maps, this square commutes.  Since the corresponding square involving the $ \sh(K_\delta) $ also commutes, we deduce that the scalars $ \lambda $ satisfy 
\begin{equation*}
\lambda_{\delta \, \delta \sqcup s} \lambda_{\delta \sqcup s \, \delta \sqcup s,t} = \lambda_{\delta \, \delta \sqcup t} \lambda_{\delta \sqcup t \, \delta \sqcup s,t}.
\end{equation*}
For each $ \delta = \{s_1, \dots, s_k \}$, we define a non-zero scalar $ \mu_\delta $ by $ \mu_\delta := \lambda_{\emptyset \, \{s_1 \}} \cdots \lambda_{\{s_1, \dots, s_{k-1}\} \, \{s_1, \dots, s_k \}} $.  By the above relation, $ \mu_\delta $ does not depend on the ordering of $ \delta $ used in its definition. 

Now, we can modify the isomorphism $ R_\delta \cong \sh(K_\delta)[|\delta|-n]\{-|\delta|\}$ by multiplying by the non-zero scalar $ \mu_\delta $.  Then, the $h_{\delta \delta'} $ and the $ g_{\delta \delta'} $ agree exactly and we deduce that the complex $ (M_\bullet(K), g_\bullet) $ is isomorphic to the complex $ (R_\bullet, h_\bullet) $.  Thus, by construction $ \sh(K) $ is a convolution of $ (R_\bullet, h_\bullet) $.

To show that $ (R_\bullet, h_\bullet) $ has a unique convolution, by Proposition \ref{th:uniquecone}, we see that it is sufficient to show that $ \Hom(R_i[k], R_{i+k+1}) $ vanishes for all $ k \ge 1 $ and $ i \ge 0 $.  To compute these hom spaces, we may pass to the category of bigraded vector spaces by taking the cohomologies of the complexes involved.  If $ U^{\bullet,\bullet}, W^{\bullet,\bullet} $ are bigraded vector spaces, then $ \Hom(U^{\bullet,\bullet}, W^{\bullet,\bullet}) $ vanishes unless there exists $ j, l$ with $ U^{j,l} \ne 0 $ and $ W^{j,l} \ne 0 $.

Because $ H^{j,l}(V) $ is supported along the diagonal $ j+l = 0 $, we see that $ H^{j,l}(R_\delta) $ is supported where $ n+ j+l = 0 $.  Hence $ H^{j,l}(R_i[k]) $ is supported where $ j+l+k +n = 0 $ while $H^{j,l}(R_{i+k+1}) $ is supported where $ j+l+n = 0 $.  Hence for $ k \ge 1$, we conclude that $\Hom(R_i[k], R_{i+k+1}) = 0 $.

\end{proof}

\begin{proof}[Proof of Theorem \ref{th:Khov}]
By the relation between our unoriented theory and our oriented theory, we see that $ \Psi(K)(\C) = \Q(K)[r]\{s\} $ where $ r, s $ are integers depending on $ K $ defined in \ref{se:uorientdef}.

On the other hand, from the definition of Khovanov homology, $ \Khh^{i,j}(K) = H^{i+r+s, j+s}[[K]] $ since $ r+s = k_2 +k_3 = n_- $  and $ s = 2(k_2+k_3) - (k_1 + k_4) = 2n_- - n_+ $ where $ n_+ , n_- $ are the number of right and left handed crossings in $ K $.  

Thus by Theorem \ref{th:Q=M}, we have that 
$$\Khh^{i+j,j}(K) = H^{i+j+r+s,j+s}[[K]] = H^{i+r,j+s}(M(K)) \cong H^{i+r,j+s}(\Q(K)) = \Kh^{i,j}(K)$$
as desired.
\end{proof}

\begin{Remark}
As Mikhail Khovanov pointed out to us, it should be possible to recover Khovanov homology with $ \mathbb{Z} $ coefficients, by working with varieties over $ \mathbb{Z} $ rather than over $ \C $ (note that $Y_n $ is defined over $ \mathbb{Z} $).  This will also reduce the ambiguity of the scalars from $ \C^\times $ to $\mathbb{Z}^\times = \{\pm 1 \}$ in the cobordism part of the theory.
\end{Remark}

\section{Reduced Homology} \label{se:reduced}

In this section we construct a reduced homology for marked links. To do this we first note that there is a standard bijection between oriented $(1,1)$ tangles and oriented links with a marked component. Namely, given such a link $K$ you can cut it anywhere along the marked component to obtain a $(1,1)$ tangle $K_{cut}$ and conversely you can glue the ends of a $(1,1)$ tangle to obtain a link with a distinguished component. We define the reduced homology $\rH_{\text{alg}}$ of the marked link $K$ by
$$\rH^{i,j}_{\text{alg}}(K) = \Hom_{Y_1}^{i,j}(\E_2, \Psi(K_{cut})(\E_1)[1])$$
where $Y_1 \cong \p$ and $\E_1 \cong \O_\p(-1) = \E_2^\vee$. 

\begin{Proposition} \label{th:reducedlong}
The homology of a marked link $K$ is related to the reduced homology of $K$ and its mirror $K^!$ by the long exact sequence 
$$\dots \rightarrow \rH^{-i,-j+1}_{\mathrm{alg}}(K^!)^\vee \rightarrow \rH^{i-1,j+1}_{\mathrm{alg}}(K) \rightarrow \Kh^{i,j}(K) \rightarrow \rH^{-i-1,-j+1}_{\mathrm{alg}}(K^!)^\vee \rightarrow \rH^{i,j+1}_{\mathrm{alg}}(K) \rightarrow \dots$$
\end{Proposition}
\begin{proof}
Consider the standard projection $\pi: Y_n \rightarrow Y_{n-1}$. Denote by $\O_{Y_n} \in D(Y_{n-1} \times Y_n)$ the FM kernel corresponding to $\pi_\ast$. Let $T$ be an $(n-1,m-1)$ tangle and denote by $TI$ the $(n,m)$ tangle obtained by adding an extra strand to the right of $T$. The following lemma relates the two worlds of tangles: those of type (odd,odd) and (even,even). 

\begin{Lemma} \label{th:evenodd} If $\P(T)$ and $\P(TI)$ denote the FM kernels associated to $T$ and $TI$ then 
\begin{equation*}
\P(TI) \ast \O_{Y_n} \cong \O_{Y_m} \ast \P(T) \ \text{ and }\ \P(T) \ast {\O_{Y_n}}_R \cong {\O_{Y_m}}_R \ast \P(TI)
\end{equation*}
In other words, $\pi^\ast \circ \Phi_{\P(T)} = \Phi_{\P(TI)} \circ \pi^\ast$ and $\Phi_{\P(T)} \circ \pi_\ast = \pi_\ast \circ \Phi_{\P(TI)}$. 
\end{Lemma}
\begin{proof}
It is enough to prove the isomorphisms when $T$ is either a cap, cup or crossing. For a cap the left isomorphism amounts to showing $\sG_{n+2}^i \ast \O_{Y_n} \cong \O_{Y_{n+2}} \ast \sG_{n+1}^i \in D(Y_{n-1} \times Y_{n+2})$. The proof follows much as in Proposition \ref{prop:capscups}. By Corollary \ref{cor:transverse} the intersection $W = \pi_{12}^{-1}(Y_n) \cap \pi_{23}^{-1}(X_{n+2}^i)$ is transverse so that 
$$\pi_{12}^\ast(\O_{Y_n}) \otimes \pi_{23}^\ast(\O_{X_{n+2}^i}) = \O_W \in D(Y_{n-1} \times Y_n \times Y_{n+2}).$$
A routine computation shows
\begin{equation*}
\begin{aligned}
W = \{ (L_\cdot, L'_\cdot, L''_\cdot) \in Y_{n-1} \times Y_n \times Y_{n+2} :\ &L_j = L'_j \text{ for } j \le n-1, \\
&L'_j = L''_j \text{ for } j \le i-1,\ L'_j = zL''_{j+2} \text{ for } j \ge i-1 \}.
\end{aligned}
\end{equation*}
so that the map $W \rightarrow \pi_{13}(W)$ is one-to-one and
\begin{equation*}
\begin{aligned}
\pi_{13}(W) = \{ (L_\cdot, L'_\cdot) \in Y_{n-1} \times Y_{n+2} :\ &L_j = L'_j \text{ for } j \le i-1, \\
 &L_j = zL'_{j+2} \text{ for } n-1 \ge j \ge i-1 \}.
\end{aligned}
\end{equation*}
Since $\sG_{n+2}^i = \O_{X_{n+2}^i} \otimes \E_i'\{-i+1\}$ we get 
\begin{eqnarray*} 
\sG_{n+2}^i \ast \O_{Y_n} &=& \pi_{13_\ast}(\O_W \otimes \E_i''\{-i+1\}) \\
&=& \O_{\pi_{13}(W)} \otimes \E_i'\{-i+1\}
\end{eqnarray*}
Similarly, the intersection $W' = \pi_{12}^{-1}(X_{n+1}^i) \cap \pi_{23}^{-1}(Y_{n+2})$ is transverse. We find
\begin{equation*}
\begin{aligned}
W' = \{ (L_\cdot, L'_\cdot, L''_\cdot) \in Y_{n-1} \times Y_{n+1} \times Y_{n+2} :\ &L_j = L'_j \text{ for } j \le i-1,\ L_j = zL'_{j+2} \text{ for } j \ge i-1 \\
&L'_j = L''_j \text{ for } j \le n+1 \}
\end{aligned}
\end{equation*}
so that $W' \rightarrow \pi_{13}(W')$ is one-to-one and
\begin{equation*}
\begin{aligned}
\pi_{13}(W') = \{ (L_\cdot, L'_\cdot) \in Y_{n-1} \times Y_{n+2} :\ &L_j = L'_j \text{ for } j \le i-1, \\
 &L_j = zL'_{j+2} \text{ for } n-1 \ge j \ge i-1 \}.
\end{aligned}
\end{equation*}
Thus $\pi_{13}(W') = \pi_{13}(W)$. Since on $W' \subset Y_{n-1} \times Y_{n+1} \times Y_{n+2}$ we have $\E_i' \cong \E_i''$ we find that 
$$\O_{Y_{n+2}} \ast \sG_{n+1}^i = \O_{\pi_{13}(W')} \otimes \E_i'\{-i+1\} \cong \sG_{n+2}^i \ast \O_{Y_n}.$$
A similar argument proves the left isomorphism when $T$ is a cup. A crossing corresponds to a kernel of the form $\Cone(\sG_n^i \ast \sF_n^i[-1]\{1\} \rightarrow \O_\Delta)$ (up to shifts). So the result for caps and cups implies the result for crossings. This proves the isomorphism on the left. Since $T$ is arbitrary, after taking right adjoints we also get the isomorphism on the right. 
\end{proof}

We will use the maps $Y_0 \xleftarrow{q} X_{2}^1 = Y_1 \xrightarrow{i} Y_2$. Since the link $K$ is obtained from $K_{cut}$ by gluing the two ends we have 
\begin{eqnarray*}
\Psi(K)(\cdot) &=& \Phi_{\sF_2^1} \circ \Psi(K_{cut}I) \circ \Phi_{\sG_2^1} (\cdot) \\
&=& q_\ast \left( i^\ast \left[ \Phi_{\P(K_{cut}I)} i_\ast(q^\ast(\cdot) \otimes \E_1) \right] \otimes \E_2^\vee \{1\} \right) \\
&=& q_\ast \left( i^\ast \sF \otimes i^\ast \E_2^\vee \{1\} \right)
\end{eqnarray*}
where $\sF = \Phi_{\P(K_{cut}I)} i_\ast (q^\ast(\cdot) \otimes \E_1)$. Denote by $\pi: Y_2 \rightarrow Y_1$ the standard projection map and notice that $\pi \circ i = id_{Y_1}$. Thus we can write $i^\ast \sF = \pi_\ast ( i_\ast i^\ast \sF)$. But $i_\ast i^\ast \sF = \sF \otimes \O_{Y_1}$ and $\O_{Y_1}$ has the standard resolution $\O_{Y_2}(-Y_1) \rightarrow \O_{Y_2}$. This gives us the distinguished triangle
$$ \sF \rightarrow i_\ast i^\ast \sF \rightarrow \sF \otimes \O_{Y_2}(-Y_1)[1].$$
Applying $q_\ast(\pi_\ast(\cdot) \otimes i^\ast \E_2^\vee\{1\})$ we obtain
$$q_\ast(\pi_\ast(\sF) \otimes i^\ast \E_2^\vee\{1\}) \rightarrow \Psi(K)(\cdot) \rightarrow q_\ast(\pi_\ast(\sF \otimes \O_{Y_2}(-Y_1)[1]) \otimes i^\ast \E_2^\vee \{1\}).$$
By Lemma \ref{th:evenodd} we have
$$\pi_\ast(\sF) = \pi_\ast \left( \Phi_{\P(K_{cut}I)} i_\ast (q^\ast(\cdot) \otimes \E_1) \right) = \Phi_{\P(K_{cut})} \pi_\ast i_\ast (q^\ast(\cdot) \otimes \E_1) = \Phi_{\P(K_{cut})}(q^\ast(\cdot) \otimes \E_1).$$ 
Thus 
$$q_\ast(\pi_\ast(\sF) \otimes i^\ast \E_2^\vee \{1\})  = \Hom_{Y_1}(\E_2, \Phi_{\P(K_{cut})}(q^\ast(\cdot) \otimes \E_1[1]))[-1]\{1\}.$$
On the other hand, $\O_{Y_2}(-Y_1) = \E_2 \otimes \E_1^\vee \{-2\}$. Hence
\begin{equation*}
q_\ast(\pi_\ast(\sF \otimes \O_{Y_2}(-Y_1)[1]) \otimes i^\ast \E_2^\vee\{1\}) = \Hom_{Y_1}(i^\ast \E_2 \{-1\}, \pi_\ast(\sF \otimes \E_2 \otimes \E_1^\vee \{-2\}[1])) 
\end{equation*}

But $\sF$ is supported on the Springer fibre $F_1 = Y_1 \subset Y_2$. To see this, note that $i_\ast(q^\ast(\cdot) \otimes \E_1)$ is supported on $F_1$.  Also, every basic functor associated to a cup, cap or crossing maps sheaves supported on the Springer fibre $F_n$ to sheaves supported on $F_{n-1}, F_{n+1}$ or $F_n$ respectively. Since $\Phi_{\P(K_{cut}I)}$ is a composition of functors associated to cups, caps and crossings $\sF$ must be supported on $F_1$.  

Since $\E_2$ and $\E_1^\vee$ are isomorphic in a neighbourhood of $Y_1 \subset Y_2$ we get $\sF \otimes \E_2 = \sF \otimes \E_1^\vee$ and so
\begin{equation*}
\begin{aligned}
\Hom_{Y_1}(i^\ast \E_2 \{-1\}, \pi_\ast(\sF \otimes \E_2 \otimes \E_1^\vee \{-2\}[1])) &= \Hom_{Y_1}(\E_1^\vee \{-1\}, \pi_\ast(\sF) \otimes \E_1^\vee \otimes \E_1^\vee \{-2\}[1])) \\
&=\Hom_{Y_1}(\E_1, \Phi_{\P(K_{cut})}(q^\ast(\cdot) \otimes \E_1) [1]\{-1\}
\end{aligned}
\end{equation*}
On the other hand, 
\begin{equation*}
\begin{aligned}
\Hom_{Y_1}(\E_1, \Phi_{\P(K_{cut})}(\E_1)) &= \Hom_{Y_1}({\Phi_{\P(K_{cut})}}_L(\E_1), \E_1)) \\
&= \Hom_{Y_1}(\Phi_{\P(K_{cut}^!)}(\E_1), \E_1) \\
&= \Hom_{Y_1}(\E_1, \Phi_{\P(K_{cut}^!)}(\E_1) \otimes \omega_{Y_1}[1])^\vee \\
&= \Hom_{Y_1}(\E_2, \Phi_{\P(K_{cut}^!)}(\E_1)[1])^\vee.
\end{aligned}
\end{equation*}
Here we used the fact that if $T$ is an $(n,n)$ tangle then the left adjoint of $\Phi_{\P(T)}$ is $\Phi_{\P(T')}$ where $T'$ is the mirror of $T$ read backwards. To see this notice that the left adjoint of the kernel corresponding to a crossing is the kernel associated to the reversed crossing. Moreover, the left adjoint of a cup or cap is a cap or cup with a $[-1]\{1\}$ or $[1]\{-1\}$ shift respectively. Since the number of cups and caps is equal the shifts cancel each other. This means 
$$\Hom_{Y_1}^{i+1,j-1}(\E_1, \Phi_{\P(K_{cut})}(\E_1)) = \Hom_{Y_1}^{-i-1,-j+1}(\E_2, \Phi_{\P(K_{cut}^!)}(\E_1)[1])^\vee = \rH^{-i-1,-j+1}_{\text{alg}}(K^!)^\vee$$
and the result follows. 
\end{proof}

\begin{Remark}
One can show that the graded Euler characteristic of $\rH^{i,j}_{\text{alg}}(K)$ is the Jones polynomial of $K$ normalized so that the unknot is assigned $1$ instead of $q+q^{-1}$. In \cite{Kpatterns} Khovanov defines a reduced homology $\rH_{\text{Kh}}$ whose graded Euler characteristic is also the normalized Jones polynomial. The proof in section \ref{se:unorient} showing the equivalence between $\Kh$ and $\Khh$ can be extended to show that $\rH^{i,j}_{\text{alg}}(K) \cong \rH^{i+j,j}_{\text{Kh}}(K)$. The only change is in lemma \ref{th:localKh} where one needs to check again that the maps coincide. 
The fact that reduced Khovanov homology fits into a long exact sequence as in Proposition \ref{th:reducedlong} is well-known to experts and was explained to us by Misha Khovanov.
\end{Remark}

\end{document}